\documentclass[reqno,12pt]{amsart}
\usepackage{a4wide}
\usepackage [frenchb]{babel}

\numberwithin{equation}{section}

\newcounter{saveeqn}
\newcommand{\alphaeqn}{\setcounter{saveeqn}{\value{equation}}%
\setcounter{equation}{0}%
\global\def\theequation{\mbox{\thesection.\arabic{saveeqn}\alph{equation}}}}
\newcommand{\reseteqn}{\setcounter{equation}{\value{saveeqn}}%
\global\def\theequation{\thesection.\arabic{equation}}}

\newtheorem{Theorem}{Th{\'e}or{\`e}me}
\newtheorem{conj}{Conjecture}
\newtheorem{Corollary}{Corollaire}
\newtheorem{Proposition}{Proposition}
\newtheorem{Lemma}{Lemme}
\theoremstyle{remark}
\newtheorem*{Remark}{Remarque}

\def\al{\alpha}
\def\be{\beta}

\def\de{\delta}
\def\ep{\varepsilon}

\def\et{\eta}

\def\Ga{\Gamma}

\def\({\left(}
\def\){\right)}
\def\[{\left[}
\def\]{\right]}
\def\ii{\infty}

\def\Li{\operatorname{Li}}
\def\dis{\displaystyle}
\def\dd{\textup{d}}

\begin{document}

\title[]{Hyperg{\'e}om{\'e}trie et fonction z{\^e}ta de Riemann}
\author[]{C. Krattenthaler$^\dagger$ et T. Rivoal}

\address{
Institut Girard Desargues,
Universit{\'e} Claude Ber\-nard Lyon-I,
21, avenue Claude Ber\-nard,
F-69622 Villeurbanne Cedex, France}
\email{kratt@igd.univ-lyon1.fr}
\address{
Laboratoire de Math{\'e}matiques Nicolas Oresme, CNRS UMR 6139,
Universit{\'e} de Caen,  BP 5186,
14032 Caen cedex,
France}
\email{rivoal@math.unicaen.fr}
\thanks{$^\dagger$ Recherche partiellement support{\'e}e par la Fondation
Autrichienne de la Recherche Scientifique FWF, contrat P12094-MAT,
et par le Programme {\og Accro{\^\i}tre le potentiel humain de recherche \fg}
de la Commission Europ{\'e}enne, contrat HPRN-CT-2001-00272,
``Algebraic Combinatorics in Europe"}



\begin{abstract}
Nous d{\'e}montrons la {\og conjecture des d{\'e}nominateurs\fg}
du deuxi{\`e}\-me auteur~\cite{ri3} sur le d{\'e}nominateur commun des coefficients
des combinaisons lin{\'e}aires en les valeurs de la fonction
z{\^e}ta de Riemann,  r{\'e}cemment construites
pour minorer la dimension de l'espace vectoriel engendr{\'e} sur $\mathbb Q$
par $1,\zeta(m),\zeta(m+2),\dots,\zeta(m+2h)$, o{\`u} $m$ et $h$ sont des
nombres entiers, $m\ge 2$ et $h\ge 0$.
En particulier, comme corollaires imm{\'e}diats, on obtient l'irrationalit{\'e}
d'au moins un des huit nombres
$\zeta(5),\zeta(7),\dots,\zeta(19)$ et l'existence d'un entier impair $j$ entre $5$ et $165$
tel que 1, $\zeta(3)$ et $\zeta(j)$ sont
lin{\'e}airement ind{\'e}pendants sur $\mathbb{Q}$,
ce qui am{\'e}liore des r{\'e}sultats de~\cite{ri2} et~\cite{br}, respectivement.
Nous prouvons {\'e}galement une conjecture connexe, due {\`a}
Vasilyev~\cite{va2},
ainsi qu'une conjecture de Zudilin~\cite{zu2} portant sur certaines
approximations rationnelles de $\zeta(4)$.
Les d{\'e}monstrations sont bas{\'e}es sur une identit{\'e} entre une somme
simple et une somme multiple, de nature hyperg{\'e}om{\'e}trique,
due \`a Andrews~\cite{Andr}.
Nous esp{\'e}rons que notre construction pourra aussi {\^e}tre appliqu{\'e}e
aux combinaisons lin{\'e}aires
plus g{\'e}n{\'e}rales construites par Zudilin~\cite{zu3}, afin d'am{\'e}liorer
son r{\'e}sultat sur l'irrationalit{\'e} d'au moins un des nombres
$\zeta(5),\zeta(7),\zeta(9),\zeta(11)$.
\medskip

\noindent
{\sc Abstract}.
We prove the second author's ``denominator
conjecture" 
\cite{ri3} concerning the common denominators of
coefficients of certain linear forms in zeta values. These forms were recently
constructed to obtain lower bounds for the dimension of the vector space over
$\mathbb Q$ spanned by $1,\zeta(m),\zeta(m+2),\dots,\zeta(m+2h)$, where $m$ and $h$
are integers such that $m\ge2$ and $h\ge0$. In particular, we
immediately get the following results as 
\hbox{corollaries\!\!:} at least one
of the eight numbers $\zeta(5),\zeta(7),\dots,\zeta(19)$ is irrational, and
there exists an odd integer $j$ between $5$ and $165$ such that $1$,
$\zeta(3)$ and $\zeta(j)$ are
linearly independent over $\mathbb{Q}$. This strengthens
some recent results in~\cite{ri2} and~\cite{br}, respectively. We also
prove a related conjecture, due to Vasilyev~\cite{va2},
and as well a conjecture, due to Zudilin~\cite{zu2}, on certain rational
approximations of $\zeta(4)$. The proofs are
based on a hypergeometric identity between a single sum and a multiple
sum due to Andrews~\cite{Andr}.
We hope that it will be possible to  apply our construction
to the more general linear forms constructed by Zudilin~\cite{zu3}, with the
ultimate goal of strengthening his result that one of
the numbers $\zeta(5),\zeta(7),\zeta(9),\zeta(11)$ is irrational.

\end{abstract}

\maketitle


\newpage

\begin{center}
\sc Table des mati{\`e}res
\end{center}
\contentsline {section}{\tocsection {}{1}{Introduction et plan de l'article}}{2}
\contentsline {section}{\tocsection {}{2}{Arri{\`e}re plan}}{4}
\contentsline {subsection}{\tocsubsection {}{2.1}{Le Th{\'e}or{\`e}me d'Ap{\'e}ry}}{4}
\contentsline {subsection}{\tocsubsection {}{2.2}{L'ind{\'e}pendance
lin{\'e}aire d'une infinit{\'e} de $\zeta $ impairs et la transcendance de~$\pi $}}{6}
\contentsline {subsection}{\tocsubsection {}{2.3}{{\`A} la recherche
d'un irrationnel parmi $\zeta (5)$, $\zeta (7)$, etc.}}{8}
\contentsline {subsection}{\tocsubsection {}{2.4}{La conjecture des d{\'e}nominateurs}}{9}
\contentsline {subsection}{\tocsubsection {}{2.5}{Les int{\'e}grales de Vasilyev}}{11}
\contentsline {section}{\tocsection {}{3}{Les r{\'e}sultats principaux}}{11}
\contentsline {section}{\tocsection {}{4}{Cons{\'e}quences diophantiennes du Th{\'e}or{\`e}me~\ref {thm:2}}}{13}
\contentsline {section}{\tocsection {}{5}{Le principe des
d{\'e}monstrations des Th{\'e}or{\`e}mes~\ref{thm:2} {\`a} \ref{thm:C=3}}}{14}
\contentsline {section}{\tocsection {}{6}{Deux identit{\'e}s entre une somme simple et une somme multiple}}{16}
\contentsline {section}{\tocsection {}{7}{Quelques explications}}{19}
\contentsline {section}{\tocsection {}{8}{Des identit{\'e}s hyperg{\'e}om{\'e}trico-harmoniques}}{22}
\contentsline {section}{\tocsection {}{9}{Corollaires au Th{\'e}or{\`e}me \ref {thm:A1}}}{31}
\contentsline {section}{\tocsection {}{10}{Corollaires au Th{\'e}or{\`e}me \ref {thm:1}}}{33}
\contentsline {section}{\tocsection {}{11}{Lemmes arithm{\'e}tiques}}{37}
\contentsline {section}{\tocsection {}{12}{D{\'e}monstration du Th{\'e}or{\`e}me \ref {thm:2}, partie i) }}{49}
\contentsline {section}{\tocsection {}{13}{D{\'e}monstration du Th{\'e}or{\`e}me \ref {thm:2}, partie ii) }}{51}
\contentsline {section}{\tocsection {}{14}{D{\'e}monstration du Th{\'e}or{\`e}me
\ref{thm:3}, partie {\rm i)}, et des Th{\'e}or{\`e}mes \ref{prop7gene} et \ref{prop:Phi}}}{55}
\contentsline {section}{\tocsection {}{15}{D{\'e}monstration du Th{\'e}or{\`e}me
\ref{thm:3}, partie {\rm ii)}, et du Th{\'e}or{\`e}me
\ref{thm:C=3}}}{57}
\contentsline {section}{\tocsection {}{16}{Encore un peu d'hyperg{\'e}om{\'e}trie}}{64}
\contentsline {section}{\tocsection {}{17}{Perspectives}}{66}
\contentsline {subsection}{\tocsubsection {}{17.1}{Les s{\'e}ries asym{\'e}triques de Zudilin}}{66}
\contentsline {subsection}{\tocsubsection {}{17.2}{La conjecture des d{\'e}nominateurs li{\'e}e aux
valeurs de la fonction beta}}{67}
\contentsline {subsection}{\tocsubsection {}{17.3}{La $q$-conjecture des d{\'e}nominateurs}}{68}
\contentsline {subsection}{\tocsubsection {}{17.4}{La version non-termin{\'e}e des identit{\'e}s gigantesques}}{68}
\contentsline {section}{\tocsection {}{}{Remerciements}}{70}
\contentsline {section}{\tocsection {}{}{Bibliographie}}{70}

\section{Introduction et plan de l'article}
La d{\'e}termination de la nature arithm{\'e}tique des
 valeurs aux entiers impairs $s\ge 3$ de la fonction z{\^e}ta de Riemann
$$
\zeta(s)=\dis\sum_{k=1}^{\ii}\frac{1}{k^s}
$$
est un des probl{\`e}mes parmi les plus difficiles de la th{\'e}orie des nombres.
Pour d{\'e}montrer les quelques r{\'e}sultats connus (voir les
paragraphes~\ref{ssec:hisapery},~\ref{ssec:hisinfzeta} et~\ref{sec:conj1}),
la seule m{\'e}thode d'attaque disponible
consiste {\`a} construire, gr{\^a}ce {\`a} divers proc{\'e}d{\'e}s hyperg{\'e}om{\'e}triques,
des suites de combinaisons lin{\'e}aires $(S_n)_{n\ge 0}$
en les valeurs de z{\^e}ta aux entiers $l\in\{2,\ldots, M\}$ 
et dont les coefficients
sont des rationnels, {\'e}ventuellement nuls~:
\begin{equation*}
S_n=p_{0,n}+\sum_{l=2}^M p_{l,n}\zeta(l).
\end{equation*}
Pour appliquer les crit{\`e}res d'irrationalit{\'e} ou d'ind{\'e}pendance
lin{\'e}aire, tel le crit{\`e}re de Nesterenko
\cite{ne}, il est n{\'e}cessaire de d{\'e}terminer un d{\'e}nominateur
commun aux $p_{l,n}$ qui soit
{\it le plus petit possible}.
Typiquement, ce d{\'e}nominateur est la puissance $M$-i{\`e}me du plus petit commun multiple des entiers
$1, 2, \ldots, n$, que nous notons $\textup{d}_n$ comme de coutume.
Or, dans certaines circonstances (li{\'e}es {\`a} la nature
hyperg{\'e}om{\'e}trique tr{\`e}s sp{\'e}ciale des constructions propos{\'e}es
et r{\'e}sum{\'e}es par les Conjectures~\ref{conj1},~\ref{conj2}
et~\ref{conj3} au paragraphe
\ref{sec:conj1}), on constate num{\'e}riquement que ce d{\'e}nominateur semble
{\^e}tre $\textup{d}_n^{M-1}$,
voire mieux,  ce qui permet d'am{\'e}liorer significativement certains r{\'e}sultats diophantiens
sur les valeurs de z{\^e}ta.

Nous prouvons ici toutes ces conjectures ({\`a} un facteur 2 pr{\`e}s),
ce qui constitue
nos Th{\'e}or{\`e}\-mes~\ref{thm:2},~\ref{thm:4} et~\ref{thm:3} au paragraphe \ref{resultats}.
De plus, les Th{\'e}or{\`e}mes~\ref{prop7gene}, \ref{prop:Phi} et
\ref{thm:C=3} dans le m{\^e}me paragraphe contiennent m{\^e}me
des am{\'e}liorations dans des cas sp{\'e}ciaux.
Les d{\'e}monstrations de ces th{\'e}or{\`e}mes, donn{\'e}es
aux paragraphes~\ref{sec:demothm1i} {\`a}
\ref{sec:Phi},  sont bas{\'e}es sur deux
identit{\'e}s hyperg{\'e}om{\'e}triques {\og gigantesques \fg}
(une due {\`a} Andrews \cite{Andr}, 
l'autre {\'e}tant une variante),
donn{\'e}es au
paragraphe~\ref{gigantesques}, et en fait surtout sur quelques unes de
leurs sp{\'e}cialisations, {\'e}nonc{\'e}es et d{\'e}montr{\'e}es 
aux paragraphes~\ref{CorollairesA} et \ref{Corollaires}.
Avant cela, nous donnons au paragraphe~\ref{sec:hyperharm} quelques
identit{\'e}s particuli{\`e}rement {\'e}l{\'e}gantes qui se d{\'e}duisent de
ces identit{\'e}s gigantesques et qui concernent les coefficients
{\og dominants
\fg} $p_{M,n}$ des combinaisons lin{\'e}aires les plus simples que l'on puisse construire.
En fait, les identit{\'e}s {\'e}nonc{\'e}es {\`a}
la Proposition~\ref{cor:A2} ont {\'e}t{\'e} le
point de d{\'e}part de ce travail, comme expliqu{\'e} au
paragraphe~\ref{sec:explications}, o{\`u} nous indiquons
comment nous est
venue l'id{\'e}e de leur improbable existence et comment elles 
nous ont men{\'e} {\`a} comprendre l'importance
des identit{\'e}s fondamentales du
paragraphe~\ref{gigantesques} 
(et qu'elles etaient, essentiellement, d{\'e}j{\`a} connues depuis trente 
ans). Tous les calculs hyperg{\'e}om{\'e}triques
des paragraphes~\ref{gigantesques}, \ref{sec:explications},
\ref{sec:hyperharm}, 
\ref{sec:Phi} et
\ref{sec:hyp} ont {\'e}t{\'e} faits
{\`a} l'aide du programme HYP, d{\'e}velopp{\'e} sous {\sl Mathematica} par le 
premier auteur \cite{hyp}%
\footnote{\label{foot1}Le programme HYP
\cite{hyp} permet de faire ais\'ement, et sans faute, des calculs routiniers
avec les s{\'e}ries hyperg{\'e}om{\'e}triques.
Cependant, tous ces calculs peuvent
aussi {\^e}tre faits sans probl{\`e}me {\`a} la
main, sauf que cela impliquerait peut-{\^e}tre plusieurs heures du travail.}. 
Pour aider le lecteur de comprendre
nos d{\'e}monstrations, nous avons
ins{\'e}r{\'e} le paragraphe~\ref{sec:idee}, o{\`u} nous expliquons {\`a}
grands traits l'id{\'e}e et la structure des d{\'e}monstrations de nos
th{\'e}or{\`e}mes principaux.

Les cons{\'e}quences diophantiennes de notre travail
sur les valeurs de z{\^e}ta aux entiers impairs sont formul{\'e}es
au paragraphe~\ref{sec:consdio}.
Notons que la preuve, qui nous {\'e}chappe encore,
des versions les plus g{\'e}n{\'e}rales
de ces conjectures (voir le paragraphe
\ref{ssec:zudconj})
pourrait impliquer
l'irrationalit{\'e}\footnote{certes bien connue, mais uniquement gr{\^a}ce au racourci
li{\'e} {\`a} la transcendance de $\pi$ et pas {\`a} la Ap{\'e}ry.
On obtiendrait ainsi la meilleure
mesure d'irrationalit{\'e} connue de $\pi^4$ (voir \cite{zu2}). 
Pour les meilleures
mesures d'irrationalit\'e connues pour $\pi$, 
respectivement $\zeta(2)$ et $\zeta(3)$,
voir Hata~\cite{hata}, 
respectivement Rhin et Viola~\cite{rv1, rv2}~: dans les trois cas, les
m\'ethodes utilis\'ees sont de nature hyperg\'eom\'etrique.}
de $\zeta(4)=\pi^4/90$, et surtout celle d'au moins un des trois nombres
$\zeta(5)$, $\zeta(7)$ et $\zeta(9)$, bien que nous n'ayons aucune
certitude {\`a} ce sujet.

Pour souligner davantage le don d'ubiquit\'e des
objets hyperg{\'e}om{\'e}triques dont nous nous servons, nous
ajoutons le paragraphe~\ref{sec:hyp}, o{\`u} nous montrons que
l'{\'e}quivalence de la s{\'e}rie de Beukers, Gutnik et Nesterenko
\eqref{eq:bgn} et de la s{\'e}rie de Ball \eqref{eq:Ball}, et
{\'e}galement l'{\'e}quivalence des s{\'e}ries \eqref{eq:Apery2} et
\eqref{eq:Ball2}, est une cons{\'e}quence directe de la transformation
de Whipple sous sa forme non termin{\'e}e due {\`a} Bailey~(une remarque qui
n'a apparemment pas {\'e}t{\'e} faite auparavant sous
cette forme).

Nous terminons notre article en {\'e}voquant au
paragraphe~\ref{sec:persp} certaines directions de
recherche que nous poursuivrons dans l'avenir.
Notamment, nous
discutons au sous-pa\-ra\-graphe~\ref{ssec:zudconj} des s{\'e}ries plus
g{\'e}n{\'e}rales de Zudilin et les raisons pour lesquelles
nous pensons que nos m{\'e}thodes devraient permettre d'attaquer les
conjectures associ{\'e}es {\`a} ces s{\'e}ries.
Aux sous-paragraphes~\ref{ssec:beta} et
\ref{ssec:q}, nous mentionnons bri{\`e}vement quelques conjectures
portant sur des combinaisons lin{\'e}aires en les valeurs de la fonction
beta et en les valeurs de $q$-z{\^e}ta, respectivement,
et leurs cons{\'e}quences possibles. Finalement, au
sous-paragraphe~\ref{ssec:nonterm}, nous
{\'e}voquons l'int{\'e}r{\^e}t de l'extension {\'e}ventuelle de nos
identit{\'e}s au cas de s{\'e}ries infinies.



\section{Arri{\`e}re plan}
\label{sec:histoire}

\subsection{Le Th{\'e}or{\`e}me d'Ap{\'e}ry}\label{ssec:hisapery}
La preuve de l'irrationalit{\'e} de $\zeta(3)$, due {\`a} Ap{\'e}ry
\cite{ap}, ne date que de 1978.
Sa d{\'e}monstration, qui fonctionne aussi pour $\zeta(2)=\pi^2/6$,
peut {\^e}tre synth{\'e}tis{\'e}e
ainsi\footnote{Voir le survol de Fischler~\cite{fi} pour
un expos{\'e}, apparemment exhaustif, des tr{\`e}s nombreuses preuves maintenant  disponibles de
l'irrationalit{\'e} de $\zeta(3)$.}~: il existe
deux suites $(a_n)_{n\ge 0}$ et
$(b_n)_{n\ge 0}$  telles que
$a_n\in \mathbb{Z}$, $\textup{d}_n^3b_n \in \mathbb{Z}$
 et
$$
\lim_{n\to +\ii}\vert 2a_n\zeta(3)-b_n\vert^{1/n} = (\sqrt{2}-1)^{4},
$$
o{\`u} $\textup{d}_n$ est le p.p.c.m des entiers $1, 2, \ldots, n$.
On conclut en remarquant que, en
vertu du th{\'e}or{\`e}me des nombres premiers,
$\textup{d}_n=e^{n+o(n)}$ et que $e^3(\sqrt{2}-1)^4<1$.
Il existe de nombreuses  fa{\c c}ons de produire ces suites, par
exemple au moyen de la s{\'e}rie suivante, due {\`a} Beukers, Gutnik et
Nesterenko \cite{beuk}, \cite{gu}, \cite{ne2}~:
\begin{equation} \label{eq:bgn}
-\sum_{k=1}^{\ii}\frac {\partial} {\partial k}
\bigg(\frac {(k-n)_n^2} {(k)_{n+1}^2}\bigg)
=
2a_n\zeta(3)-b_n,
\end{equation}
o{\`u} les symboles de Pochhammer sont d{\'e}finis par
$(\al)_k=\alpha(\alpha+1)\cdots(\alpha+k-1)$ si $k\ge 1$ et
$(\al)_0=1$. Ici et dans tout cet article, $n$
d{\'e}signe un entier positif, sauf mention contraire.

De
m{\^e}me, dans le cas de $\zeta(2)$, on construit
deux suites $(\al_n)_{n\ge 0}$ et
$(\be_n)_{n\ge 0}$  telles que
$\al_n\in \mathbb{Z}$, $\textup{d}_n^2\be_n \in \mathbb{Z}$
 et
$$
\lim_{n\to+\ii}\vert \al_n\zeta(2)-\be_n\vert^{1/n}
=\bigg(\frac{\sqrt{5}-1}{2}\bigg)^{5}.
$$
L{\`a} aussi, il existe
beaucoup de fa{\c c}ons de
g{\'e}n{\'e}rer ces suites, par exemple au moyen de la s{\'e}rie
\begin{equation} \label{eq:Apery2}
(-1)^n n!\sum_{k=1}^{\ii}\frac {(k-n)_n} {(k)_{n+1}^2}
=
\al_n\zeta(2)-\be_n.
\end{equation}
Les entiers  $a_n$ et $\al_n$ peuvent {\^e}tre
explicit{\'e}s sous forme
binomiale ou hyperg{\'e}om{\'e}trique~:
\begin{equation} \label{eq:a_n}
a_n=\sum_{j=0}^n\binom{n}{j}^2\binom{n+j}{n}^2
=
{}_{4}F_3
 \!\left [ \begin{matrix} {-n, -n, n+1, n+1}\\ {
    1,\;1,\;1}\end{matrix} ; {\displaystyle 1}\right ]
\end{equation}
et
\begin{equation} \label{eq:alpha_n}
\al_n=\sum_{j=0}^n\binom{n}{j}^2\binom{n+j}{n}
=
{}_{3}F_2\!\left [ \begin{matrix} {-n, -n, n+1}\\ {
      1,\; 1}\end{matrix} ; {\displaystyle 1}\right ].
\end{equation}
D'une
fa{\c c}on g{\'e}n{\'e}rale, les s{\'e}ries (ou fonctions)
hyperg{\'e}om{\'e}\-tri\-ques sont d{\'e}finies par
\begin{eqnarray*}
{}_{q+1}F_q\!\left [ \begin{matrix}
{ \alpha_0,\alpha_1,\ldots,\alpha_{q}}\\
{\beta_1,\ldots,\beta_q}\end{matrix} ;
{\displaystyle z}\right ]
=\sum_{k=0}^{\ii}
\frac{(\alpha_0)_k\,(\alpha_1)_k\cdots(\alpha_{q})_k}
{k!\,(\beta_1)_k\cdots(\beta_q)_k} z^k,
\label{eq:hyper}
\end{eqnarray*}
o{\`u}
$\alpha_j\in\mathbb{C}$ et
$\beta_j\in\mathbb{C}\setminus\mathbb{Z}_{\le 0}$.
La s{\'e}rie converge pour tout $z\in\mathbb{C}$ tel que $\vert z\vert <1$,
et aussi pour $z=\pm1$ lorsque
$\text{Re}(\be_1+\dots+\be_q)>\text{Re}(\al_0+\al_1+\dots+\al_q)$.
Dans les  ouvrages traitant de ces fonctions (par exemple
\cite{AAR}, \cite{BailAA}, \cite{GaRaAA},  \cite{SlatAC}), on
trouve les d{\'e}finitions suivantes~:
la s{\'e}rie hyperg{\'e}om{\'e}trique ${}_{q+1}F_q$ est dite

\bigskip

\begin{itemize}
\item  \textit{balanc{\'e}e} (balanced) si
$\al_0+\cdots+\al_q+1=\be_1+\cdots +\beta_q$~;

\medskip
\item \textit{quasi {\'e}quilibr{\'e}e de premi{\`e}re
esp{\`e}ce} (nearly-poised of the first kind) si
$\al_1+\beta_1=\cdots=\alpha_{q}+\beta_q$~;

\medskip
\item \textit{bien {\'e}quilibr{\'e}e} (well-poised) si
$\alpha_0+1=\al_1+\beta_1=\cdots=\alpha_{q}+\beta_q$~;

\medskip
\item \textit{tr{\`e}s bien {\'e}quilibr{\'e}e}
(very-well-poised) si
elle est bien {\'e}quilibr{\'e}e et
$\alpha_1=\frac12\,\alpha_0+1$.
\end{itemize}

\bigskip

\noindent Ces s{\'e}ries v{\'e}rifient d'innombrables identit{\'e}s
recens{\'e}es dans les livres cit{\'e}s ci-dessus%
\footnote{%
Le manuel du logiciel
HYP, d{\'e}j{\`a} mentionn{\'e} {\`a} la note \ref{foot1} de bas de page,
contient la plus grande liste
d'identit{\'e}s hyperg{\'e}om{\'e}triques actuellement disponible.}.
Les s{\'e}ries (tr{\`e}s) bien {\'e}quilibr{\'e}es
y sont abondamment repr{\'e}sent{\'e}es, {\`a} la mesure de leur {\'e}norme influence
sur le d{\'e}veloppement de la th{\'e}orie hyperg{\'e}om{\'e}trique
au cours du XXi{\`e}me si{\`e}cle
(voir le survol d'Andrews \cite{andr2}
et le livre \cite{AAR} d'Andrews, Askey et Roy {\`a} ce sujet).
Le pr{\'e}sent article n'{\'e}chappe pas {\`a} cette influence.

Dans \cite{ne2}, Nesterenko a pos{\'e} le probl{\`e}me de trouver une
preuve de l'irrationalit{\'e} de
$\zeta(3)$ aussi {\'e}l{\'e}mentaire que celle du nombre
$e=\sum_{n\ge 0} 1/n!$, due {\`a} Fourier. Pour attaquer ce probl{\`e}me,
Ball a  introduit la
s{\'e}rie hyperg{\'e}om{\'e}trique tr{\`e}s bien
{\'e}quilibr{\'e}e suivante (voir l'introduction de~\cite{ri1})~:
\begin{align} \notag
{\bf B}_n&=n!^2\sum_{k=1}^{\ii}\(k+\frac n2\)
\frac{(k-n)_n(k+n+1)_n}{(k)_{n+1}^4}
\\&=\frac{n!^7(3n+2)!}{2(2n+1)!^5}\;  _{7}F_6\!\left [ \begin{matrix}
{ 3n+2, \frac32 n+2 , n+1    , \ldots,   n+1}\\
{\frac32 n+1 ,  2n+2, \ldots,  2n+2}\end{matrix} ;
{\displaystyle 1}\right ].
\label{eq:Ball}
\end{align}
Il a alors observ{\'e} le fait remarquable que
${\bf B}_n={\bf a}_n\zeta(3)-{\bf b}_n$, alors que
l'on s'attend aussi {\`a} voir appara{\^\i}tre
$\zeta(4)$ et $\zeta(2)$.
En d{\'e}finissant le $m$-i{\`e}me nombre harmonique par $H_m=
\sum _{j=1} ^{m}\frac {1} {j}$ si $m\ge 1$ et $H_0=0$,
on a en particulier
\begin{multline} \label{eq:bolda_n}
{\bf a}_n=(-1)^{n+1}\sum_{j=0}^n \(\frac n2-j\)
 \binom{n}{j}^4\binom{n+j}{n}\binom{2n-j}{n}\\
\cdot
\left(5H_{n-j}-5H_j+H_{n+j}-H_{2n-j}-\frac {1} {\frac {n} {2}-j}\right),
\end{multline}
et on montre que $\textup{d}_n{\bf a}_n$
et $\textup{d}_n^4{\bf b}_n$ sont entiers.
Un deuxi{\`e}me point remarquable de
cette s{\'e}rie r{\'e}side dans  la propri{\'e}t{\'e},
initialement conjecturale, que
${\bf a}_n$ et $\textup{d}_n^3{\bf b}_n$ sont en fait des
entiers, et
${\bf a}_n$ et ${\bf b}_n$ co{\"\i}ncident avec
les nombres d'Ap{\'e}ry $a_n$ et $b_n/2$
pour $\zeta(3)$. L'{\'e}galit{\'e} de la s{\'e}rie \eqref{eq:Ball} et de la moiti{\'e} de la
s{\'e}rie \eqref{eq:bgn}
a ensuite {\'e}t{\'e} prouv{{\'e}e} par
Zudilin \cite{zu1} et le deuxi{\`e}me auteur ind{\'e}pendamment,
alors que l'{\'e}galit{\'e} ${\bf a}_n=a_n$ a
{\'e}t{\'e} prouv{{\'e}e} par le premier auteur, dans les deux cas
par une utilisation subtile
de l'algorithme de calcul de
r{\'e}currences lin{\'e}aires
de Gosper--Zeilberger (voir \cite{ek}, \cite{PeWZAA}, \cite{ZeilAM},
\cite{ZeilAV}), ce qui implique l'{\'e}galit{\'e} ${\bf b}_n=b_n/2$.
Voir le paragraphe~\ref{sec:hyp} pour une d{\'e}monstration
ind{\'e}pendante, utilisant
des identit{\'e}s hyperg{\'e}om{\'e}triques classiques.
Au passage, on obtient bien une
nouvelle d{\'e}monstration
de l'irrationalit{\'e} de $\zeta(3)$~: bien qu'{\'e}l{\'e}mentaire,
elle n'est cependant pas aussi simple que celle de $e$.


\subsection{L'ind{\'e}pendance lin{\'e}aire d'une infinit{\'e} de
$\zeta$ impairs et la transcendance de~$\pi$}\label{ssec:hisinfzeta}

La disparition de la moiti{\'e}  des valeurs de $\zeta$ attendues
n'est pas un miracle isol{\'e}~: elle s'explique
par la nature
(tr{\`e}s) bien {\'e}quilibr{\'e}e
de ${\bf B}_n$, alors qu'une s{\'e}rie seulement
quasi {\'e}quilibr{\'e}e ne la produit pas.
Ce pr{\'e}cieux ph{\'e}nom{\`e}ne  a  {\'e}t{\'e} g{\'e}n{\'e}ralis{\'e}
dans \cite{ri1} et \cite{br}
au moyen, essentiellement\footnote{Plus
pr{\'e}cis{\'e}ment~:
sans le facteur {\og tr{\`e}s bien {\'e}quilibrant\fg} $k+n/2$, qui ne joue aucun
r{\^o}le pour obtenir le
r{\'e}sultat vis{\'e},
mais qui est la myst{\'e}rieuse
raison d'{\^e}tre  des identit{\'e}s prouv{\'e}es dans cet article.}, de
la s{\'e}rie
\begin{align*}
\bar{\bf S}_{n,A,r}(z)&=n!^{A-2r}\sum_{k=1}^{\ii}
\left(k+\frac n2\right)\frac{(k-rn)_{rn}(k+n+1)_{rn}}{(k)_{n+1}^A}z^{-k}
\\
&=z^{-rn-1}n!^{A-2r}\frac{(rn)!^{A+1}((2r+1)n+2)!}{2((r+1)n+1)!^{A+1}}
\\
&\kern3cm
\times {}_{A+3}F_{A+2}
\!\left [ \begin{matrix}
{(2r+1)n+2, \frac{2r+1}{2}n+2, rn+1, \ldots,   rn+1}\\
{\frac{2r+1}{2}n+1, (r+1)n+2, \ldots,  (r+1)n+2}\end{matrix} ;
{\displaystyle z^{-1}}\right ].
\end{align*}
avec $\vert z\vert \ge 1$ et $A$ et $r$ des entiers tels
que $0\le r<A/2$, ce qui a permis
de montrer qu'une infinit{\'e} des valeurs de la fonction
z{\^e}ta aux entiers impairs sont lin{\'e}airement ind{\'e}pendantes sur $\mathbb{Q}$.
Esquissons rapidement
la preuve. On d{\'e}finit d'abord les fonctions polylogarithmes, pour tout
$s\ge 1$ et $z$ complexe v{\'e}rifiant $\vert z\vert \le 1$ et
$(s,z)\not=(1,1)$, 
par
$$
\dis \Li_s(z)=\sum_{n=1}^{\ii}\frac{z^n}{n^s}\,.
$$
En d{\'e}veloppant en {\'e}l{\'e}ments simples le sommande de
$\bar{\bf S}_{n,A,r}(z)$, on v{\'e}rifie
qu'il existe des polyn{\^o}mes
$\bar{\bf p}_{l,n}(X)$
(d{\'e}pendant aussi de $A$ et $r$) tels que
$\textup{d}_n^{A-l}\bar{\bf p}_{l,n}(X)\in \mathbb{Z}[X]$ et
$$
\bar{\bf S}_{n,A,r}(z)=\bar{\bf p}_{0,n}(z)+\sum_{l=1}^A
\bar{\bf p}_{l,n}(z)\Li_l(1/z).
$$
Le {\og tr{\`e}s  bon {\'e}quilibrage\fg} de
$\bar{\bf S}_{n,A,r}(z)$ se traduit par la relation de
r{\'e}ciprocit{\'e}
\begin{equation}\label{eq:recipari}
z^n\bar{\bf p}_{l,n}(1/z)=(-1)^{A(n+1)+l+1}\bar{\bf p}_{l,n}(z),
\end{equation}
dont on d{\'e}duit que pour tout  $A$ pair et
tout $n\ge 0$,
on a\footnote{\label{foot}Lorsque $z$ tend vers 1, la s{\'e}rie
$\bar{\bf S}_{n,A,r}(z)$ converge mais
 $\Li_1(1/z)=-\log(1-1/z)$ diverge~: on a donc n{\'e}cessairement
 $\bar{\bf p}_{1,n}(1)=0$, ce qui {\'e}limine la valeur divergente
 $\zeta(1)$ de la combinaison lin{\'e}aire.}
$$
\bar{\bf S}_{n,A,r}(1)=\bar{\bf p}_{0,n}(1)+
\underset{l \,\textup{impair}}{\sum_{l=3, \ldots,
A-1}} \bar{\bf p}_{l,n}(1)\zeta(l).
$$
On conclut en utilisant un crit{\`e}re d'ind{\'e}pendance
lin{\'e}aire d{\^u} {\`a} Nesterenko \cite{ne} et en
optimisant le param{\`e}tre $r$ en fonction de $A$.
Comme pour ${\bf B}_n$, on constate
exp{\'e}rimentalement que le d{\'e}nominateur
$\textup{d}_n^{A-l}$ est trop g{\'e}n{\'e}reux
lorsque $z=1$~: il semble en effet que\break
$\textup{d}_n^{A-l-1}\bar{\bf p}_{l,n}(1)$ soit
d{\'e}j{\`a} entier pour tout
$l\in\{0,\ldots, A-1\}$.
L'int{\'e}r{\^e}t est que cela permettrait  d'appliquer plus finement le
crit{\`e}re de Nesterenko (voir
la partie~ii) du Th{\'e}or{\`e}me~\ref{thm:ameliodio}
au paragraphe~\ref{sec:consdio}).

Lorsque $A$
est impair et $z=-1$, la s{\'e}rie
$\bar{\bf S}_{n,A,r}(-1)$ est une combinaison
lin{\'e}aire rationnelle en les valeurs de $\tilde\zeta$ aux entiers pairs,
o{\`u}, par d{\'e}finition,
$$
\tilde\zeta(s)=\sum_{k=1}^{\ii}\frac{(-1)^k}{k^s}=(2^{1-s}-1)\zeta(s).
$$
Puisque pour tout entier $k\ge 1$, on a
$$
\zeta(2k)=(-1)^{k-1}\frac{2^{2k-1}B_{2k}}{(2k)!}\pi^{2k},
$$
o{\`u} $B_k$ est le $k$-i{\`e}me nombre de Bernoulli, on obtient en fait une suite de
combinaisons
lin{\'e}aires rationnelles de puissances de $\pi$~: la
transcendance de~$\pi$ en d{\'e}coule par un argument
simple de th{\'e}orie des nombres alg{\'e}briques
(voir \cite{rey} pour un cas analogue).
Le m{\^e}me ph{\'e}nom{\`e}ne arithm{\'e}tique que celui mis {\`a}
jour pour $A$ pair et $z=1$ semble l{\`a} aussi se
produire~: si l'on consid{\`e}re par exemple la s{\'e}rie
\begin{equation} \label{eq:Ball2}
\bar{\bf S}_{n,3,1}(-1)=n!\sum_{k=1}^{\ii}(-1)^k\(k+\frac{n}{2}\)
\frac{(k-n)_n(k+n+1)_n}{(k)_{n+1}^3}= p_n\tilde\zeta(2)-q_n,
\end{equation}
le rationnel
\begin{multline} \label{eq:q_n}
p_n=(-1)^{n+1}\sum_{j=0}^n \(\frac n2-j\)
\binom{n}{j}^3\binom{n+j}{n}\binom{2n-j}{n}\\
\cdot
\left(4H_{n-j}-4H_j+H_{n+j}-H_{2n-j}-\frac {1} {\frac {n} {2}-j}\right),
\end{multline}
joue le m{\^e}me r{\^o}le pour $\zeta(2)$ que ${\bf a}_n$
pour $\zeta(3)$.
En effet, on a {\it a priori} que
$\textup{d}_np_n$ et $\textup{d}_n^3q_n$ sont entiers, mais
on montre par les  m{\'e}thodes
de \cite{zu1} que
$p_n=\al_n$ et $q_n=-\beta_n/2$,
o{\`u} $\al_n$ et $\be_n$ sont les nombres d'Ap{\'e}ry
pour $\zeta(2)$ 
d{\'e}finis au sous-paragraphe~\ref{ssec:hisapery}.
Voir le paragraphe~\ref{sec:hyp} pour une d{\'e}monstration utilisant
des identit{\'e}s hyperg{\'e}om{\'e}triques classiques.


\subsection{{\`A} la recherche d'un irrationnel parmi $\zeta(5)$,
$\zeta(7)$, etc.}\label{sec:conj1}
Cette recherche a {\'e}t{\'e} initi{\'e}e dans
\cite{ri2} au moyen d'un autre type de s{\'e}ries, que nous
appellerons improprement {\og s{\'e}ries d{\'e}riv{\'e}es\fg},
qui ne sont plus formellement
hyperg{\'e}om{\'e}triques mais en sont tr{\`e}s
proches\footnote{De fa{\c c}on plus pr{\'e}cise, une s{\'e}rie
d{\'e}riv{\'e}e appara{\^\i}t naturellement par la
 m{\'e}thode de Frobenius dans le calcul des
solutions des {\'e}quations
diff{\'e}rentielles hyperg{\'e}om{\'e}triques. L'{\'e}quation
diff{\'e}rentielle satisfaite par une s{\'e}rie
hyperg{\'e}om{\'e}trique (tr{\`e}s) bien {\'e}quilibr{\'e}e est
invariante par le changement de variable
$z\mapsto 1/z$, ce qui {\og explique\fg} la r{\'e}ciprocit{\'e} des
polyn{\^o}mes ${\bf p}_{l,n}(z)$. Voir \cite{no} pour une tr{\`e}s belle
exposition du calcul des solutions des
{\'e}quations diff{\'e}rentielles hyperg{\'e}om{\'e}triques.}.
Consid{\'e}rons, pour $A$ pair $\ge 6$ et $\vert z\vert \ge1$, la s{\'e}rie
$$
\tilde{\bf S}_{n,A}(z)=n!^{A-6}\sum_{k=1}^{\ii}\frac{1}{2}\frac{\partial^2}{\partial k^2} \(
\(k+\frac n2\)\frac{(k-n)_n^3(k+n+1)_n^3}{(k)_{n+1}^A}\)z^{-k}.
$$
Par rapport \`a l'expression {\og d\'evelopp\'ee \fg} de $\bar{\bf S}_{n,A,r}(1)$ 
du paragraphe pr{\'e}c{\'e}dent, 
l'introduction d'une d{\'e}rivation d'ordre $2$ permet de remplacer
$\zeta(l)$ par $\zeta(l+2)$, ce qui, avec le 
tr{\`e}s bon {\'e}quilibrage de $\tilde{\bf S}_{n,A}(z)$, montre qu'il
existe des polyn{\^o}mes
$\tilde{\bf p}_{l,n}(X)$, d{\'e}pendant de $A$, tels que
$\textup{d}_n^{A+2}\tilde{\bf p}_{0,n}(1)$ et
 $\textup{d}_n^{A-l}\tilde{\bf p}_{l,n}(1)$ ($l\ge 1$) sont entiers et
$$
\tilde{\bf S}_{n,A}(1)
=\tilde{\bf p}_{0,n}(1)+\underset{l\equiv 1
\,(\textup{mod}\,2)}{\sum_{l=3, \,\ldots, \,A-1}}
\tilde{\bf p}_{l,n}(1)\zeta(l+2).
$$
Ainsi,  on fait non seulement dispara{\^\i}tre
les valeurs de $\zeta$ aux entiers pairs, mais aussi $\zeta(3)$.
L'entier $A=20$ est le plus petit entier pair tel que
$0<\liminf_{n\to+\ii}\vert
\textup{d}_n^{A+2}\tilde{\bf S}_{n,A}(1)\vert^{1/n}<1$,
ce qui se traduit par  l'irrationalit{\'e} d'au moins un des neuf
nombres $\zeta(5), \zeta(7), \ldots, \zeta(21)$.
Sans surprise, on constate
exp{\'e}rimentalement\footnote{Ce gain arithm{\'e}tique n'est pas
anodin~: voir 
la partie~i) du Th{\'e}or{\`e}me~\ref{thm:ameliodio}
au paragraphe~\ref{sec:consdio} 
pour son utilisation diophantienne.}
que l'on pourrait prendre $\textup{d}_n^{21}$ {\`a} la place de
$\textup{d}_n^{22}$.

\subsection{La conjecture des d{\'e}nominateurs}\label{ssec:conjdenom}
Toutes ces donn{\'e}es ont
naturellement conduit {\`a} formuler une conjecture
g{\'e}n{\'e}rale sur les d{\'e}nominateurs des combinaisons lin{\'e}aires
rationnelles construites sur les s{\'e}ries
d{\'e}riv{\'e}es suivantes~:
\begin{equation}
\label{eq:seriesderiveesformeslineaires}
{\bf S}_{n,A,B,C,r}(z)=n!^{A-2Br}\sum_{k=1}^{\ii}\frac{1}{C!}
\frac{\partial^{C}}{\partial k^{C}} \(
\(k+\frac n2\)\frac{(k-rn)_{rn}^{B}(k+n+1)_{rn}^B}{(k)_{n+1}^A}\)z^{-k},
\end{equation}
o{\`u} $\vert z\vert \ge1$ et $A$, $B$, $C$, $r$ sont des entiers positifs
v{\'e}rifiant\footnote{
Cette condition, qui sert uniquement {\`a} faire converger la s{\'e}rie en
$z=\pm 1$,
sera assouplie
aux paragraphes~\ref{resultats}, \ref{sec:hyperharm} 
{\`a}
\ref{Corollaires}, et
\ref{sec:demothm1i}
{\`a} \ref{sec:demothmvaszud}.}
$0\le 2Br<A$. Il existe alors des polyn{\^o}mes
${\bf p}_{0, C, n}(X)$ et ${\bf p}_{l,n}(X)$ 
pour $l\in\{1, \ldots, A\}$,
d{\'e}pendant de $A$, $B$ et $r$ mais pas de $C$, 
tels que
\begin{equation} \label{eq:poldef}
\textup{d}_n^{A+C}{\bf p}_{0, C, n}(X) \in\mathbb{Z}[X],\quad
\textup{d}_n^{A-l}{\bf p}_{l,n}(X)\in \mathbb{Z}[X]
\end{equation}
et
\begin{equation} \label{eq:Sdef}
{\bf S}_{n,A,B,C,r}(z)={\bf p}_{0,C,n}(z)+
(-1)^C\sum_{l=1}^A\binom{C+l-1}{l-1} {\bf p}_{l,n}(z)\Li_{C+l}(1/z).
\end{equation}

Posons pour simplifier
$$
R_{n,A,B,r}(k)=n!^{A-2Br}\(k+\frac n2\)\frac{(k-rn)_{rn}^{B}(k+n+1)_{rn}^B}{(k)_{n+1}^A}.
$$
En suivant la d{\'e}marche classique (voir \cite{br}), on
a pour tout
$l\in\{1,\dots, A\}$~:
\begin{equation}
{{\bf p}_{l,n}\(X\)}=
\sum_{j=0}^n\frac{1}{(A-l)!}\frac{\partial^{A-l}}
{\partial k^{A-l}}\big(R_{n,A,B,r}(k)(k+j)^A\big)\bigg\vert_{k=-j}\,X^j
\label{eq:expressionp_mn}
\end{equation}
et
\begin{align} \notag
{\bf p}_{0,C,n}(X)&=-\sum _{j=0} ^{n}
\sum _{e=1} ^{A}(-1)^C \binom {C+e-1}{e-1}\\
&\kern2cm\cdot
\bigg(\frac {1} {(A-e)!}
\frac {\partial^{A-e}}
{\partial k^{A-e}}\big(R_{n,A,B,r}(k)(k+j)^A\big)\Big\vert_{k=-j}\bigg)
\sum _{i=1} ^{j}\frac {1} {i^{e+C}}\,X^{j-i}.
\label{eq:p0}
\end{align}
Nous aurons besoin de l'expression plus explicite suivante :
\begin{multline}
\frac {\partial^{h}}
{\partial k^{h}}\big(R_{n,A,B,r}(k)(k+j)^A\big)\Big\vert_{k=-j}\\
=(-1)^{Aj+Brn}\frac {(rn)!^{2B}} {n!^{2rB}}
\frac {\partial^h} {\partial\ep^h}\(\frac {n} {2}-j+\ep\)
\(\frac {n!} { (1-\ep)_j \, (1+\ep)_{n-j}}\)^{A}\\
\cdot
{\binom {  r n+j-\ep}{rn}}^{B}
{\binom { (r+1) n-j+\ep}{rn}}^{B}\Bigg\vert_{\ep=0}.
\label{eq:expressionp_mnplusexplicite}
\end{multline}

Quand on sp{\'e}cialise \eqref{eq:Sdef} en $z=(-1)^A$, on
obtient
\begin{multline} \label{eq:SnABCr}
{\bf S}_{n,A,B,C,r}\((-1)^A\)\\={\bf p}_{0,C,n}\((-1)^A\)+
(-1)^C\sum_{l=1}^A\binom{C+l-1}{l-1} {\bf p}_{l,n}\((-1)^A\)\Li_{C+l}\((-1)^A\).
\end{multline}
De nouveau, la relation de r{\'e}ciprocit{\'e}~\eqref{eq:recipari} vaut
avec ${\bf p}_{l,n}(z)$ {\`a} la place de $\bar{\bf p}_{l,n}(z)$, et,
de m{\^e}me, la note \ref{foot} de bas de page 
(avec ${\bf S}_{n,A,B,C,r}(z)$ {\`a} la place de $\bar{\bf S}_{n,A,r}(z)$)
s'applique dans la
situation plus g{\'e}n{\'e}rale que l'on consid{\`e}re ici.
Par cons{\'e}quent, si $A$ est pair, la s{\'e}rie
${\bf S}_{n,A,B,C,r}\((-1)^A\)$
est une combinaison lin{\'e}aire en 1, $\zeta(C+3), \zeta(C+5),
\ldots, \zeta(C+A-1)$, alors que si  $A$ impair, c'est une combinaison
lin{\'e}aire en 1, $\tilde\zeta(C+2), \tilde\zeta(C+4),
\ldots, \tilde\zeta(C+A-1)$.
Dans les deux cas, chacun des coefficients est
un rationnel dont un d{\'e}nominateur est donn{\'e} par
\eqref{eq:poldef}.
En particulier, si on multiplie ${\bf S}_{n,A,B,C,r}\((-1)^A\)$ par
$\textup{d}_n^{A+C}$, on obtient une combinaison lin{\'e}aire {\`a} coefficients
entiers en
les valeurs de z{\^e}ta
(respectivement de z{\^et}a {\og altern{\'e}e\fg} $\tilde\zeta$)
aux entiers impairs, respectivement aux entiers pairs.
Cependant, num{\'e}riquement, il appara{\^\i}t que l'on peut
esp{\'e}rer mieux.

\begin{conj}\label{conj1} Dans les conditions ci-dessus, pour tout
$A\ge 2$ et pour tout $l\in\{1, \ldots, A\}$, les nombres
$\textup{d}_n^{A+C-1}{\bf p}_{0,C,n}\((-1)^A\)$ et
 $\textup{d}_n^{A-l-1}{\bf p}_{l,n}\((-1)^A\)$ sont entiers.
\end{conj}

\begin{Remark}
Nous rappelons que l'int{\'e}r{\^e}t de cette conjecture
est qu'elle permettrait d'appliquer plus finement le
crit{\`e}re de Nesterenko
(voir 
le Th{\'e}or{\`e}me~\ref{thm:ameliodio} au paragraphe~\ref{sec:consdio}).
Pour $A$ pair, la conjecture est formul{\'e}e dans \cite[p.~51]{ri3}, et
pour $A$ impair, dans \cite[Remarque~2.14]{fi}.
\end{Remark}

\bigskip

Zudilin \cite[paragraphe~2]{zu2}
a consid{\'e}r{\'e} plus en d{\'e}tail la s{\'e}rie d{\'e}riv{\'e}e
$$
{\bf S}_{n,4,2,1,1}(1)=\sum_{k=1}^{\ii}\frac{\partial}{\partial k} \(
\(k+\frac n2\)\frac{(k-n)_{n}^{2}(k+n+1)_{n}^2}{(k)_{n+1}^4}\)
={\bf u}_n\zeta(4)-{\bf v}_n,
$$
o{\`u} $\textup{d}_n {\bf u}_n$ et $\textup{d}_n^5 {\bf v}_n$ sont entiers.
En utilisant l'algorithme de Gosper--Zeilberger
\cite{ek}, \cite{PeWZAA}, \cite{ZeilAM}, \cite{ZeilAV}, il a constat{\'e} que
les suites $({\bf u}_n)_{n\ge 0}$
et  $({\bf v}_n)_{n\ge 0}$ v{\'e}rifient la r{\'e}currence
lin{\'e}aire d'ordre deux suivante\footnote{Cette r{\'e}currence
a aussi {\'e}t{\'e} obtenue par Cohen et Rhin \cite{cr}
par la m{\'e}thode d'Ap{\'e}ry en 1981
et par Sorokin \cite{so} en 2001~; il est remarquable que
les trois m{\'e}thodes soient {\it a priori} totalement ind{\'e}pendantes.
}~:
\begin{small}
\begin{multline*}
(n+1)^5Y_{n+1}=3(2n+1)(3n^2+3n+1)(15n^2+15n+4)Y_n+3n^2(3n-1)(3n+1)Y_{n-1}.
\end{multline*}
\end{small}%
En s'aidant de cette r{\'e}currence
(pour calculer un grand nombre de valeurs de
${\bf u}_n$ et ${\bf v}_n$),
Zudilin a affin{\'e} la Conjecture~\ref{conj1} dans ce cas.
Posons avec lui (voir la fin du paragraphe~2 dans \cite{zu2})
\begin{equation}
\label{eq:Phi}
\Phi_n=\underset{\{n/p\}\in[2/3,1[}{\prod_{p\text{ premier}}} p\ ,
\end{equation}
o{\`u} $\{n/p\}$ est la partie fractionnaire de $n/p$.
\begin{conj}\label{conj2}
Pour tout entier $n\ge 0$, les nombres
$\Phi_n^{-1}{\bf u}_n$ et $\Phi_n^{-1}\textup{d}_n^4{\bf v}_n$ sont entiers.
\end{conj}


\subsection{Les int{\'e}grales de Vasilyev}\label{sec:conjvasilyev}
Une des nombreuses d{\'e}monstrations de l'irrationalit{\'e} de
$\zeta(2)$ et $\zeta(3)$ utilise les
c{\'e}l{\`e}bres int{\'e}grales de Beukers \cite{be}~:
$$
\int_0^1\int_0^1\frac{x^n(1-x)^ny^n(1-y)^n}{(1-(1-x)y)^{n+1}}
\,\dd x\,\dd y 
=\al_n\zeta(2)-\be_n
$$
et
$$
\int_0^1\int_0^1\int_0^1\frac{x^n(1-x)^ny^n(1-y)^nz^n(1-z)^n}
{(1-(1-(1-x)y)z)^{n+1}}
\,\dd x\,\dd y\,\dd z
=2a_n\zeta(3)-b_n.
$$
{\`A} la suite de Vasilenko \cite{vas}, Vasilyev \cite{va1}, \cite{va2} a
consid{\'e}r{\'e} des int{\'e}grales qui
g{\'e}n{\'e}ralisent naturellement celles de Beukers~:
$$
J_{E,n}=\int_{[0,1]^E}\frac{\prod_{j=1}^E x_j^n(1-x_j)^n}
{Q_E(x_1,x_2,\ldots, x_E)^{n+1}}
\,\dd x_1\,\dd x_2\,\cdots\, \dd x_E,
$$
o{\`u} $Q_E(x_1,x_2,\ldots, x_E)=
1-(\cdots (1-(1-x_E)x_{E-1})\cdots)x_1$.
Il a alors formul{\'e} la conjecture
suivante, qu'il a prouv{\'e}e pour $E=4$ et $5$, et qui est aussi
vraie pour $E=2$ et $3$ depuis Beukers.
\begin{conj}
\label{conj3}\leavevmode\newline
{\em i)} Pour tous entiers $E\ge 2$ et $n\ge 0$, il existe
des rationnels $p_{l,E,n}$ tels que
\begin{equation}\label{eq:vasint}
J_{E,n}=p_{0,E,n}+\underset{l \equiv E \,
(\textup{mod}\, 2)}{\sum_{l=2, \ldots, E}} p_{l,E,n} \zeta(l).
\end{equation}
{\em ii)} De plus, $\textup{d}_n^{E}p_{l,E,n}$ est un entier pour
tout $l\in\{0, 1, \ldots, E\}$.
\end{conj}
La partie i) de cette conjecture a {\'e}t{\'e} d{\'e}montr{\'e}e
par Zudilin \cite[paragraphe~8]{zu2} au moyen
d'une identit{\'e}
inattendue entre des int{\'e}grales g{\'e}n{\'e}ralisant celles
de Vasilyev et certaines  s{\'e}ries hyperg{\'e}om{\'e}triques
tr{\`e}s bien {\'e}quilibr{\'e}es\footnote{Voir aussi~\cite{KR} pour une
nouvelle d{\'e}monstration de cette identit{\'e}, bas{\'e}e sur
celle d'Andrews~\cite{Andr}~: nous la rappelons au
Th{\'e}or{\`e}me~\ref{thm:A1}
parce qu'elle est aussi centrale dans le pr{\'e}sent article.}.
Mais le d{\'e}nominateur alors obtenu est $\textup{d}_n^{E+1}$ et il reste
donc une puissance de $\textup{d}_n$ {\`a} {\'e}liminer.
En effet, l'identit{\'e} de \cite[Theorem~5]{zu2}, se lit dans ce cas~:
$$
J_{E,n}=\frac{n!^{2E+1}(3n+2)!}{(2n+1)!^{E+2}}
\;{}_{E+4}F_{E+3}
\!\left [ \begin{matrix}
{ 3n+2, \frac32 n+2 , n+1    , \ldots,   n+1}\\
{\frac32 n+1 ,  2n+2, \ldots,  2n+2}\end{matrix} ;
{\displaystyle (-1)^{E+1}}\right ],
$$
c'est-{\`a}-dire $J_{E,n}={\bf S}_{n,E+1,1,0,1}((-1)^{E+1})$ avec
la notation employ{\'e}e au paragraphe pr{\'e}\-c{\'e}\-dent.
Par cons{\'e}quence, 
la Conjecture~\ref{conj3}, 
partie~ii) d{\'e}coule de la
Conjecture~\ref{conj1}.


\section{Les r{\'e}sultats principaux}\label{resultats}

Dans ce paragraphe,
nous pr\'esentons les th\'eor\`emes centraux de cet article. 
Les Th{\'e}o\-r{\`e}mes~\ref{thm:2} {\`a} \ref{thm:3} prouvent les
Conjectures~\ref{conj1} 
{\`a} \ref{conj3}, {\`a} un facteur 2 pr{\`e}s, et
les Th{\'e}or{\`e}mes~\ref{prop7gene} {\`a} \ref{thm:C=3} les affinent
dans des cas sp{\'e}ciaux. Notons que la
restriction {\og analytique\fg} $0\le 2Br< A$
n'intervient pas dans ces th{\'e}or{\`e}mes~:
les {\'e}nonc{\'e}s sont valables sans
aucune hypoth{\`e}se de ce type, ni sur $r\ge 0$
qui est quelconque.
Les preuves des th{\'e}or{\`e}mes sont don\-n{\'e}es aux
paragraphes~\ref{sec:demothm1i} {\`a} \ref{sec:Phi}~;
elles sont bas{\'e}es sur un
certain nombre de corollaires, que nous {\'e}noncerons
aux paragraphes~\ref{CorollairesA} et \ref{Corollaires}, des
Th{\'e}or{\`e}mes~\ref{thm:A1} et \ref{thm:1} au
paragraphe~\ref{gigantesques}.

\begin{Theorem} \label{thm:2}\leavevmode\newline
\hangindent15pt\hangafter2
{\em i)} La Conjecture~{\em\ref{conj1}} est vraie quels que soient
  $A\ge 2$, $B\ge 1$, $C\ge 0$ et $r\ge 0$ pour tous les coefficients
${\bf p}_{l,n}\((-1)^A\)$, $l\in\{1,\ldots, A\}$, c'est-{\`a}-dire que
$\dd_n^{A-l-1}{\bf p}_{l,n}\((-1)^A\)$ est un nombre entier.

\noindent
\leavevmode\newline
\hangindent15pt\hangafter2
{\em ii)}
De plus, dans les m{\^e}mes conditions, les coefficients
$2\textup{d}_n^{A+C-1}{\bf p}_{0,C,n}\((-1)^A\)$
sont des nombres entiers.
\end{Theorem}

\begin{Theorem} \label{thm:4}\leavevmode\newline
\hangindent15pt\hangafter2
{\em i)} La Conjecture~{\em\ref{conj3}}, {\em ii)} est vraie
pour tous les coefficients $p_{l,E,n}$, $l\in\{1,\ldots, E\}$,
c'est-{\`a}-dire que $\textup{d}_n^{E-l}p_{l,E,n}$ est un nombre entier.

\noindent
\leavevmode\newline
\hangindent15pt\hangafter2
{\em ii)}
De plus, les coefficients
$2\textup{d}_n^{E}p_{0,E,n}$ sont des nombres entiers.
\end{Theorem}

Comme not{\'e} {\`a} la fin du
paragraphe~\ref{sec:conjvasilyev},
la Conjecture~\ref{conj3} est un cas particulier de la
Conjecture~\ref{conj1}~: le Th{\'e}or{\`e}me~\ref{thm:4} est donc
une cons{\'e}quence directe du Th{\'e}or{\`e}me~\ref{thm:2}.

\begin{Theorem} \label{thm:3}
\leavevmode\newline
\hangindent15pt\hangafter2
{\em i)} La Conjecture~{\em\ref{conj2}} est vraie pour le coefficient
${\bf u}_n$, c'est-{\`a}-dire que $\Phi_n^{-1}{\bf u}_n$ est un nombre
entier,
o{\`u} $\Phi_n$ est le nombre d{\'e}fini dans \eqref{eq:Phi}.

\noindent
{\em ii)} De plus, $2\Phi_n^{-1}\textup{d}_n^4{\bf v}_n$ est un nombre
entier, c'est-\`a-dire, qu'au facteur $2$ pr{\`es}, la
Conjecture~{\em\ref{conj2}} est aussi vraie pour le coefficient
${\bf v}_n$.
\end{Theorem}

Ce th{\'e}or{\`e}me est {\'e}videmment une am{\'e}lioration du
Th{\'e}or{\`e}me~\ref{thm:2} dans le cas o{\`u} $r=1$, $A=4$, $B=2$ et
$C=1$.
Il est alors naturel d'attendre des am{\'e}liorations similaires pour les
coefficients les plus g\'en\'eraux.
Par exemple, le r{\'e}sultat suivant est
une am{\'e}lioration du cas $r=1$ du
Th{\'e}or{\`e}me~\ref{thm:2} pour le coefficient dominant
${\bf p}_{A-1,n}\((-1)^A\)$. Il contient en m{\^e}me temps
la partie~i) du Th{\'e}or{\`e}me~\ref{thm:3}.

\begin{Theorem}\label{prop7gene}
Pour $r=1$, $A\ge2$ et $B\ge 1$,
le nombre $\Phi_n^{-B+1}{\bf p}_{A-1,n}\((-1)^A\)$ est entier,
\end{Theorem}

Pour les autres coefficients (toujours dans le cas o{\`u} $r=1$), nous
pouvons d{\'e}montrer un r{\'e}sultat un peu plus faible.
Au lieu de $\Phi_n$, consid{\'e}rons la
quantit{\'e} inf{\'e}rieure
\begin{equation*}
\tilde\Phi_n=\underset{\{n/p\}\in[2/3,1[}{\prod_{p\text{
premier},\,\,p<n}} p\ .
\end{equation*}
Nous avons alors le th{\'e}or{\`e}me suivant.

\begin{Theorem}\label{prop:Phi} Pour $r=1$, $A\ge2$, $B\ge 1$, $C\ge0$, et
pour tout $l\in\{1, \ldots, A\}$, les nombres\break
$\tilde\Phi_n^{-B+1}\textup{d}_n^{A-l-1}{\bf p}_{l,n}\((-1)^A\)$ et
$2\tilde\Phi_n^{-B+1}\textup{d}_n^{A+C-1}{\bf p}_{0,C,n}\((-1)^A\)$
sont entiers.
\end{Theorem}

Le cas $r=1$, $A=6$, $B=3$, $C=2$ de ce th{\'e}or{\`e}me r{\'e}pond {\`a} la
question sur la s{\'e}rie $\tilde F_n$ {\`a} la fin du paragraphe~7
dans \cite{zu2}, {\`a} un facteur 2 pr{\`e}s.

En principe, on pourrait aussi s'attendre {\`a} remplacer
$\tilde\Phi_n^{B-1}$ dans le Th{\'e}or{\`e}me~\ref{prop:Phi}
par $\Phi_n^{B-1}$.
Des calculs num{\'e}riques sugg{\`e}rent cependant qu'une divisibilit{\'e}
par $\Phi_n^{B-1}$
du d{\'e}nominateur commun des coefficients n'a lieu que pour
$A=4$, $B=2$ et $C=1$ ou $3$~:
le cas $A=4$, $B=2$ et $C=1$ est couvert
par le Th{\'e}or{\`e}me~\ref{thm:3} et
notre dernier r{\'e}sultat confirme
les observations num{\'e}riques pour $A=4$, $B=2$ et $C=3$.

\begin{Theorem} \label{thm:C=3}
Pour $r=1$, $A=4$, $B=2$, le nombre
$2\Phi_n^{-1}\textup{d}_n^{6}{\bf p}_{0,3,n}\(1\)$ est entier.
\end{Theorem}


\section{Cons{\'e}quences diophantiennes du Th{\'e}or{\`e}me~\ref{thm:2}}\label{sec:consdio}
Nous mentionnons maintenant deux applications imm{\'e}diates de
la conjecture des d{\'e}\-no\-mi\-na\-teurs, dont la preuve nous permet
d'am{\'e}liorer les r{\'e}sultats suivants~:
{\og {\it au moins un des neuf nombres $\zeta(5), \zeta(7), \ldots,
\zeta(21)$ est irrationnel} (\cite{ri2})\fg} et {\og {\it
il existe un entier impair $j$ entre $5$ et $169$,  tel que $1$, $\zeta(3)$ et $\zeta(j)$
soient lin{\'e}airement ind{\'e}pendants sur $\mathbb{Q}$} (\cite{br})\fg}.

\begin{Theorem}\label{thm:ameliodio}\leavevmode \newline
\noindent
{\em i)} Au moins un des huit nombres $\zeta(5), \zeta(7), \ldots,
\zeta(19)$ est irrationnel.

\noindent
{\em ii)} Il existe un entier impair $j$ entre $5$ et $165$,  tel que $1$,
$\zeta(3)$ et $\zeta(j)$
soient lin{\'e}airement ind{\'e}\-pen\-dants sur $\mathbb{Q}$.
\end{Theorem}

\begin{proof} Nous ne faisons que l'esquisser car elle suit les lignes de
celles de~\cite{ri2} et~\cite{br}.

i) Le th{\'e}or{\`e}me de \cite{ri2} est prouv{\'e} par l'utilisation de la s{\'e}rie
$\tilde{\bf S}_{n,20}(1)$ introduite au paragraphe~\ref{sec:conj1}, avec un {\og mauvais\fg}
d{\'e}nominateur $\textup{d}_n^{22}$
pour les combinaisons lin{\'e}aires rationnelles
en 1, $\zeta(5), \zeta(7), \ldots,\zeta(21)$ construites.
Le Th{\'e}or{\`e}me~\ref{thm:2}
nous permet de faire le m{\^e}me travail avec la s{\'e}rie $\tilde{\bf S}_{n,18}(1)$ et un {\og bon\fg}
d{\'e}nominateur $2\,\textup{d}_n^{19}$
pour les combinaisons lin{\'e}aires rationnelles
en 1, $\zeta(5), \zeta(7), \ldots, \zeta(19)$.

ii) Le th{\'e}or{\`e}me de \cite{br} est d{\'e}montr{\'e} {\`a} l'aide d'une
s{\'e}rie bien {\'e}quilibr{\'e}e qui est
$\bar{\mathbf S}_{n,169,10}(1)$ du paragraphe~\ref{ssec:hisinfzeta}, sans le facteur
$k+n/2$. Notons qu'avec la s{\'e}rie
$\bar{\mathbf S}_{n,168,10}(1)$, on peut d{\'e}j{\`a} d{\'e}montrer le r{\'e}sultat
avec 167 {\`a} la place de 169, en construisant
des combinaisons lin{\'e}aires rationnelles en 1, $\zeta(3), \zeta(5), \ldots, \zeta(167)$, avec un
{\og mauvais\fg}
d{\'e}nominateur $\textup{d}_n^{168}$. Notre am{\'e}lioration repose sur
la s{\'e}rie $\bar{\mathbf S}_{n,166,10}(1)$,
qui permet de construire des combinaisons lin{\'e}aires rationnelles
en 1, $\zeta(3), \zeta(5), \ldots, \zeta(165)$, avec un {\og bon\fg}
d{\'e}nominateur $2\,\textup{d}_n^{165}$ comme cons{\'e}quence du  Th{\'e}or{\`e}me~\ref{thm:2}.
\end{proof}

\begin{Remark} Si l'on essaie de tirer avantage de la divisibilit{\'e} des
combinaisons lin{\'e}aires
$\tilde{\bf S}_{n,A}(1)$ par $\tilde{\Phi}_n^2$ (Th{\'e}or{\`e}me~\ref{prop:Phi}),
alors il s'en faut d'extr{\^e}mement peu que l'on parvienne {\`a} montrer
l'irrationalit{\'e} de l'un des sept
nombres 1, $\zeta(5), \zeta(7), \ldots, \zeta(17)$~:
on a en effet $\liminf_{n\to\ii}
\tilde{\Phi}_n^{-2}\dd_n^{17}\tilde{\bf S}_{n,16}(1)\approx 1.007$.
\end{Remark}

Ces am{\'e}liorations ne sont pas n{\'e}gligeables  mais on sait
obtenir beaucoup mieux. En effet,
Zudilin a montr{\'e} qu'{\og {\it au moins un des quatre
nombres $\zeta(5), \zeta(7), \zeta(9), \zeta(11)$ est irrationnel}\fg}
(\cite{zu3})
et qu'{\og{\it il existe un entier impair $j$ entre $5$ et $69$,
tel que $1$, $\zeta(3)$ et $\zeta(j)$
soient lin{\'e}airement ind{\'e}pendants sur $\mathbb{Q}$}\fg}%
\footnote{Communication personnelle de W.~Zudilin~; 
voir aussi \cite{zu4} pour un r{\'e}sultat un peu plus faible.}.
Ces r{\'e}sultats sont bas{\'e}s sur des combinaisons lin{\'e}aires en les valeurs
de z{\^e}ta construites
{\`a} l'aide de s{\'e}ries beaucoup plus g{\'e}n{\'e}rales que les n{\^o}tres.
Le gain obtenu r{\'e}sulte d'une {\'e}tude
$p$-adique tr{\`e}s fine des coefficients des combinaisons lin{\'e}aires,
ce qui permet d'{\'e}liminer de {\og gros\fg} facteurs
communs {\`a} ces coefficients, {\`a} la mani{\`e}re de la Conjecture~\ref{conj2}.
Zudilin a {\'e}galement formul{\'e} une conjecture des d{\'e}nominateurs pour ses
s{\'e}ries, qui pourrait peut-{\^e}tre
permettre de montrer l'irrationalit{\'e} d'au moins un des trois nombres
$\zeta(5), \zeta(7)$ et $\zeta(9)$.
Mais il n'est pas {\'e}vident que l'on puisse aborder cette conjecture
avec nos m{\'e}thodes, et encore moins {\'e}vident
que cela permette de prouver le r{\'e}sultat envisag{\'e},
qui n'appara{\^\i}tra {\'e}ventuellement
qu'au bout des calculs~: voir le paragraphe~\ref{ssec:zudconj}
pour plus de d{\'e}tails sur cette conjecture.

\section{Le principe des d{\'e}monstrations des
Th{\'e}or{\`e}mes~\ref{thm:2} {\`a}~\ref{thm:C=3}}\label{sec:idee}

L'id{\'e}e des d{\'e}monstrations des
Th{\'e}or{\`e}mes~\ref{thm:2} {\`a}
\ref{thm:C=3} consiste {\`a} ne pas {\'e}tudier
les expressions des coefficients ${{\bf p}_{l,n}\((-1)^A\)}$ et
${\bf p}_{0,C,n}\((-1)^A\)$ donn{\'e}es par
\eqref{eq:expressionp_mn}--\eqref{eq:expressionp_mnplusexplicite}
dir{\'e}ctement, mais {\`a} chercher des expressions {\'e}quivalentes, dont
il sera alors possible d'extraire le comportement arithm{\'e}tique {\'e}nonc{\'e}
dans les th{\'e}or{\`e}mes. Pour illustration, et pour {\^e}tre plus concret,
consid{\'e}rons des coefficients sp{\'e}ciaux et commen{\c c}ons avec les
coefficients {\og dominants\fg}
${{\bf p}_{A-1,n}\((-1)^A\)}$. Par exemple,
pour $r=1$, $A=3$, $B=1$,
$C=0$ (ce choix correspond {\`a} la s{\'e}rie
$\bar{\bf S}_{n,3,1}(-1)$ dans \eqref{eq:Ball2}), le coefficient
${{\bf p}_{2,n}\(-1\)}$ est {\'e}gal {\`a}
\begin{align} \notag
 {\bf p}_{2,n}\(-1\)
&=(-1)^n
\sum _{j=0} ^{n}\frac {\partial} {\partial\ep}\Bigg( \(\frac {n} {2}-j+\ep\)
\(\frac {n!} { (1-\ep)_j \, (1+\ep)_{n-j}}\)^{3}\\
\notag
&\kern4cm\cdot{\binom {
n+j-\ep}{rn}}
{\binom { 2 n-j+\ep}{rn}}\Bigg)\Bigg\vert_{\ep=0}\\
\notag
&=(-1)^n\sum _{j=0} ^{n}\(\frac {n}
{2}-j\){\binom nj}^3{\binom {n+j}n}{\binom {2n-j}n}\\
&\kern4cm
\cdot
\(4H_{j}-4H_{n-j}+H_{2n-j}-H_{n+j}+\frac {1}
{\frac {n} {2}-j}\),
\label{eq:apery2}
\end{align}
ou pour $r=1$, $A=4$, $B=1$,
$C=0$ (ce choix correspond {\`a} la s{\'e}rie
$\mathbf B_n$ dans \eqref{eq:Ball}), le coefficient
${{\bf p}_{3,n}\(1\)}$ est {\'e}gal {\`a}
\begin{multline}
 {{\bf p}_{3,n}\(1\)}=(-1)^n\sum _{j=0} ^{n}\(\frac {n}
{2}-j\){\binom nj}^4{\binom {n+j}n}{\binom {2n-j}n}\\
\cdot
\(5H_{j}-5H_{n-j}+H_{2n-j}-H_{n+j}+\frac {1}
{\frac {n} {2}-j}\).
\label{eq:apery3}
\end{multline}
Selon le Th{\'e}or{\`e}me~\ref{thm:2}, partie~i),
il faudrait d{\'e}montrer que ces
quantit{\'e}s sont des nombres entiers, ce qui n'est pas du tout
{\'e}vident {\`a} partir
de la forme des expressions {\`a} droite dans \eqref{eq:apery2} et
\eqref{eq:apery3}, {\`a} cause de la pr{\'e}sence des nombres harmoniques
dans les sommandes. En fait, on peut se convaincre que
les sommandes sont rarement entiers. Or, au
paragraphe~\ref{ssec:hisapery}, on a d{\'e}j{\`a} remarqu{\'e} l'identit{\'e}
$a_n=\mathbf a_n$, c'est-{\`a}-dire, explicitement,
\begin{multline*}
(-1)^n\sum _{j=0} ^{n}\(\frac {n}
{2}-j\){\binom nj}^3{\binom {n+j}n}{\binom {2n-j}n}\\
\cdot
\(4H_{j}-4H_{n-j}+H_{2n-j}-H_{n+j}+\frac {1}
{\frac {n} {2}-j}\)=
\sum _{j=0} ^{n}{\binom nj}^2\binom {n+j}n,
\end{multline*}
et celle-ci d{\'e}montre imm{\'e}diatement que ${{\bf p}_{2,n}\(-1\)}$ dans
\eqref{eq:apery2} est bien un nombre entier.

Pour le coefficient ${{\bf p}_{3,n}\(1\)}$ dans \eqref{eq:apery3} on a
aussi constat{\'e} une co{\"\i}ncidence avec une somme simple dont les
sommandes sont des coefficients binomiaux, {\`a} savoir la somme dans
\eqref{eq:a_n}, mais ce type d'identit{\'e} ne se g{\'e}n{\'e}ralise pas~:
il est trop optimiste de s'attendre {\`a} exprimer ces coefficients
sous forme d'une somme {\it simple} dont tous les sommandes sont entiers.
La cl{\'e} de nos d{\'e}monstrations est de chercher des expressions sous
forme d'une somme {\it multiple}. Par exemple, nous allons d{\'e}montrer
dans la Proposition~\ref{cor:A2}
au paragraphe~\ref{sec:hyperharm} (voir aussi le
paragraphe~\ref{sec:explications}) que
\begin{multline*}
(-1)^n\sum _{j=0} ^{n}\(\frac {n}
{2}-j\){\binom nj}^4{\binom {n+j}n}{\binom {2n-j}n}
\(5H_{j}-5H_{n-j}+H_{2n-j}-H_{n+j}+\frac {1}
{\frac {n} {2}-j}\)\\
=-\sum_{0\le i\le j\le n} (-1)^j\binom{n}{j}
\binom{n}{i}^2\binom{n+j}{n}\binom{n+j-i}{n}.
\end{multline*}
{\'E}videmment cette identit{\'e} d{\'e}montre que le coefficient
${{\bf p}_{3,n}\(1\)}$ dans \eqref{eq:apery3} est un nombre entier.
Plus g{\'e}n{\'e}ralement, la Proposition~\ref{cor:A2} donne la
d{\'e}monstration du Th{\'e}or{\`e}me~\ref{thm:2} pour le coefficient
dominant (dans le cas o{\`u} $r=1$).

Malheureusement, les identit{\'e}s de la Proposition~\ref{cor:A2} ne
suffisent pas pour d{\'e}montrer le
Th{\'e}or{\`e}me~\ref{thm:2} pour les autres coefficients. Pour cela,
nous avons cherch{\'e} des identit{\'e}s plus g{\'e}n{\'e}rales entre une
somme simple et une somme multiple, le r{\'e}sultat {\'e}tant donn{\'e} par
les deux identit{\'e}s des Th{\'e}or{\`e}mes~\ref{thm:A1} et
\ref{thm:1} au paragraphe~\ref{gigantesques}.
(Voir le paragraphe~\ref{sec:explications} pour savoir comment nous avons
{\'e}t{\'e} conduits {\`a} ces identit{\'e}s.)

Comme cons{\'e}quence (voir les Corollaires~\ref{cor:1} {\`a} \ref{cor:2a} au
paragraphe~\ref{Corollaires}), nous obtenons une expression
alternative pour la somme
\begin{equation} \label{eq:somme}
\sum _{j=0} ^{n}
 \left( \frac{n}{2} - j + \ep\right)
   \left(
      \frac{n!}{( 1-\ep )_j\,
         (1+\ep )_{n-j}} \right)^A
   { \binom {  n+j-\ep} { n} }^{B}\,
   { \binom {2n-j+\ep } { n} }^{B},
\end{equation}
sous la forme $\ep\cdot \Sigma$, o{\`u} $\Sigma$ est
une somme multiple dont tous les sommandes
sont (essentiellement) des produits de coefficients binomiaux. Selon
\eqref{eq:expressionp_mn} et \eqref{eq:expressionp_mnplusexplicite}
(pour $r=1$), le coefficient ${{\bf p}_{l,n}\((-1)^A\)}$ r{\'e}sulte de
\eqref{eq:somme} en appliquant l'op{\'e}rateur
$\frac{1}{(A-l)!}\frac{\partial^{A-l}}
{\partial \ep^{A-l}}$ {\`a} l'expression \eqref{eq:somme} (au signe
pr{\`e}s), et en posant ensuite $\ep=0$. Gr{\^a}ce {\`a} l'expression
{\'e}quivalente $\ep\cdot \Sigma$, le coefficient ${{\bf
p}_{l,n}\((-1)^A\)}$ s'obtient donc aussi en appliquant l'op{\'e}rateur
$\frac{1}{(A-l-1)!}\frac{\partial^{A-l-1}}{\partial \ep^{A-l-1}}$
{\`a} la somme multiple
$\Sigma$, et
en posant ensuite {\'e}galement $\ep=0$. (C'est le contenu de \eqref{eq:3}
dans un contexte un peu plus g{\'e}n{\'e}ral.) Le point important de cet argument
est que ce ne sont que $A-l-1$ d{\'e}riv{\'e}es qui sont appliqu{\'e}es {\`a}
$\Sigma$, au lieu des $A-l$ d{\'e}riv{\'e}es qui ont {\'e}t{\'e} appliqu{\'e}es au
coefficient originel ${{\bf p}_{l,n}\((-1)^A\)}$. En utilisant plusieurs
lemmes concernant des propri{\'e}t{\'e}s arithm{\'e}tiques de
{\og briques\fg} (voir les
Lemmes~\ref{lem:briques} et \ref{lem:Faktor} au
paragraphe~\ref{sec:LemmesArithm}), il devient maintenant
possible de montrer qu'il suffit de multiplier ce coefficient par
$\dd_n^{A-l-1}$ pour obtenir un nombre entier.
(On peut adopter le {\og principe \fg} que chaque fois que l'on applique une
d{\'e}riv{\'e}e {\`a} une expression qui est un produit des coefficients
binomiaux, il faut multiplier par $\dd_n$ pour obtenir un
nombre entier.)
Cette m{\^e}me approche fonctionne aussi pour d{\'e}montrer le
Th{\'e}or{\`e}me~\ref{thm:2} dans le cas o{\`u} $r$ est quelconque~: voir le
paragraphe~\ref{sec:demothm1i}.

Pour le coefficient ${\bf
p}_{0,C,n}\((-1)^A\)$, il faut choisir une autre strat{\'e}gie
puisque la forme de l'expression qui d{\'e}finit ${\bf
p}_{0,C,n}\((-1)^A\)$ est d'une nature diff{\'e}rente.
Plus pr{\'e}cis{\'e}ment, selon \eqref{eq:p0} avec $r=1$,
le coefficient ${\bf p}_{0,C,n}\((-1)^A\)$ s'exprime sous forme d'une
somme (sur $e$ et $i$) des expressions (voir \eqref{p0alter})
\begin{equation} \label{eq:sommep0}
\frac {1} {i^{e+C}}\sum _{j=i} ^{n}
 \left( \frac{n}{2} - j + \ep\right)
   \left(
      \frac{n!}{( 1-\ep )_j\,
         (1+\ep )_{n-j}} \right)^A
   { \binom {  n+j-\ep} { n} }^{B}\,
   { \binom {2n-j+\ep } { n} }^{B},
\end{equation}
auxquelles on applique l'op{\'e}rateur
$\frac{1}{(A-e)!}\frac{\partial^{A-e}}{\partial \ep^{A-e}}$
et on pose ensuite {\'e}galement $\ep=0$. (C'est le contenu de \eqref{p0alter}
dans un contexte un peu plus g{\'e}n{\'e}ral.)
Contrairement au cas des autres coefficients, l'identit{\'e}
du Th{\'e}or{\`e}me~\ref{thm:A1} joue ici le r\^ole principal
en nous permettant d'{\'e}crire
\eqref{eq:sommep0} d'une fa{\c c}on 
diff{\'e}rente~:
cette fois-ci sous la forme $\frac {i-\ep} {2i^{e+C}}\tilde\Sigma$,
o{\`u} $\tilde\Sigma$ est
une somme multiple dont tous les sommandes
sont essentiellement des produits de coefficients binomiaux
(voir la Proposition~\ref{conj:4} au paragraphe~\ref{sec:demothm1ii}). On
applique alors les m{\^e}mes arguments que ci-dessus {\`a} la somme
multiple $\tilde\Sigma$.

Pour d{\'e}montrer les Th{\'e}or{\`e}mes~\ref{thm:3} {\`a} \ref{thm:C=3}, nous
suivons la m{\^e}me approche, mais il faut raffiner
l'analyse arithm{\'e}tique des sommandes de $\Sigma$ et $\tilde\Sigma$ (voir les
paragraphes~\ref{sec:demothmvaszud} et \ref{sec:Phi}).


\section{Deux identit{\'e}s entre une somme simple et une somme
multiple}\label{gigantesques}

Nous pr{\'e}sentons ici deux identit{\'e}s de type
hyperg{\'e}om{\'e}trique entre une s{\'e}rie tr{\`e}s bien {\'e}quilibr{\'e}e et
une somme multiple. Ces deux identit{\'e}s sont centrales dans les
d{\'e}monstrations des nos r{\'e}sultats principaux au
paragraphe~\ref{resultats}.
Plus pr{\'e}cis{\'e}ment, ce sont leurs cas sp{\'e}ciaux,
{\'e}nonc{\'e}s aux paragraphes~\ref{sec:hyperharm} {\`a}
\ref{Corollaires}, qui sont utilis{\'e}s dans les d{\'e}monstrations
aux paragraphes~\ref{sec:demothm1i} {\`a} \ref{sec:Phi}.

La premi{\`e}re identit{\'e} est un cas sp{\'e}cial
d'une identit{\'e} de nature
hyperg{\'e}om{\'e}trique basique, due \`a Andrews \cite[Theorem~4]{Andr}~:
si l'on remplace dans le th{\'e}or{\`e}me d'Andrews $a$ par $q^a$,
$b_j$ par $q^{b_j}$, $c_j$ par $q^{c_j}$, $k$ par $m+1$, $N$ par $n$,
$m_j$ par $i_j-i_{j-1}$, pour chaque~$j$ (avec la convention $i_0=0$),
et que l'on fait tendre $q$ vers $1$,
on obtient l'identit{\'e} {\'e}nonc{\'e}e dans le th{\'e}or{\`e}me suivant.

\begin{Theorem}[\sc Andrews] \label{thm:A1}
Pour tous entiers $m,n\ge0$, on a
\begin{multline}
\label{eq:Amain}
{}_{2m+5}F_{2m+4}\left[
\begin{array}{c}
a,\frac{a}{2}+1, b_1,  c_1, \ldots,  b_{m+1},  c_{m+1}, -n \\
     \frac{a}{2}, 1+a-b_1, 1+a-c_1, \ldots,  1+a-b_{m+1}, 1+a-c_{m+1}, 1+a+n
\end{array}
\,;\, 1 \right]
\\
=
\frac{(1+a)_n\,(1+a-b_{m+1}-c_{m+1})_n}{(1+a-b_{m+1})_n\,(1+a-c_{m+1})_n}
\sum _{0\le i_1\le i_2\le\dots\le i_{m}\le n} ^{}
\frac{(-n)_{i_m}}{(b_{m+1}+c_{m+1}-a-n)_{i_m}}\kern2cm\\
\cdot
\Bigg(\prod_{k=1}^m
\frac{(1+a-b_{k}-c_k)_{i_k-i_{k-1}}\,(b_{k+1})_{i_k}\,(c_{k+1})_{i_k}}
{(i_k-i_{k-1})!\,(1+a-b_k)_{i_k}\,(1+a-c_k)_{i_k}}\Bigg),
\end{multline}
o{\`u}, par d{\'e}finition, $i_0=0$ et o{\`u} dans le cas $m=0$,
la somme vide doit {\^e}tre interpr{\'e}t{\'e}e
comme valant 1.
\end{Theorem}

La d{\'e}monstration de ce th{\'e}or{\`e}me dans \cite{Andr} utilise
la transformation de Whipple entre une s{\'e}rie $_4F_3$ balanc{\'e}e
et une s{\'e}rie $_7F_6$ tr{\`e}s bien {\'e}quilibr{\'e}e
(voir \cite[(2.4.1.1)]{SlatAC})~:
\begin{multline} \label{Whipple}
{} _{4} F _{3} \!\left [ \begin{matrix} { a, b, c, -N}\\
{ e, f, 1 + a + b + c - e -
   f - N}\end{matrix} ; {\displaystyle 1}\right ] =
\frac {( -a - b + e + f)_N\,( -a - c + e
+ f) _{N}} {( -a + e + f)_N\, (-a - b - c + e + f) _{N}} \\
\times
  {} _{7} F _{6} \!\left [ \begin{matrix} { -1 - a + e + f, \frac{1}{2} -
\frac{a}{2} + \frac{e}{2} + \frac{f}{2}, -a + f, -a + e, b, c, -N}\\ {
- \frac{1}{2}  - \frac{a}{2} + \frac{e}{2} + \frac{f}{2}, e,
f, -a - b + e + f, -a - c + e + f, -a + e + f + N}\end{matrix} ;
{\displaystyle 1}\right ],
\end{multline}
o{\`u} $N$ est un entier positif, et la formule de
Pfaff--Saalsch\"utz (voir
\cite[(2.3.1.3), Appendix~(III.2)]{SlatAC})~:
\begin{equation}
\label{saal}
{} _{3} F _{2} \!\left [ \begin{matrix} { a, b, -N}\\
{ c, 1 + a + b - c -
   N}\end{matrix} ; {\displaystyle 1}\right ]  =
  {\frac{({ \textstyle c-a}) _{N} \,({ \textstyle c-b}) _{N} }
    {({ \textstyle c}) _{N} \,({ \textstyle c-a - b}) _{N} }},
\end{equation}
o{\`u} $N$ est un entier positif, de fa{\c c}on it{\'e}rative.
En particulier, l'identit{\'e} \eqref{eq:Amain} se r{\'e}duit {\`a}
l'equation \eqref{Whipple} pour $m=1$.

La deuxi{\`e}me identit{\'e} peut {\^e}tre consid{\'e}r{\'ee} comme
une variation de l'identit{\'e} d'Andrews. Elle
g{\'e}n{\'e}ralise une combinaison de la transformation de Whipple et
une transformation entre deux s{\'e}ries $_4F_3$ balanc{\'e}es (voir
\eqref{Sears}).

\begin{Theorem} \label{thm:1}
Pour tous entiers $m,n\ge1$, on a
\begin{multline}
\label{eq:main}
{}_{2m+5}F_{2m+4}\left[
\begin{array}{c}
a,\frac{a}{2}+1, b_1,  c_1, \ldots,  b_{m+1},  c_{m+1}, -n \\
     \frac{a}{2}, 1+a-b_1, 1+a-c_1, \ldots,  1+a-b_{m+1}, 1+a-c_{m+1}, 1+a+n
\end{array}
\,;\, 1 \right]
\\
=
\frac {({ \textstyle 1 + a}) _{n} \,
     ({ \textstyle 1 + a - b_{m} - c_{m+1}}) _{n} \,
     ({ \textstyle 1 + a - b_{m+1} - c_{m+1}}) _{n} \,
     ({ \textstyle 1 + a - c_{m} - c_{m+1}}) _{n} }
{   ({ \textstyle 1 + a - b_{m}}) _{n} \,
     ({ \textstyle 1 + a - b_{m+1}}) _{n} \,
     ({ \textstyle 1 + a - c_{m}}) _{n} \,
     ({ \textstyle 1 + a - c_{m+1}}) _{n} }\kern1cm\\
\times
\sum _{0\le i_1\le i_2\le\dots\le i_{m}\le n} ^{}
\frac {({ \textstyle -n}) _{i_{m}} \,
     ({ \textstyle c_{m+1}}) _{i_{m}}
}
{      ({ \textstyle -a - n + b_{m} + c_{m+1}}) _{i_{m}} \,
     ({ \textstyle -a - n + b_{m+1} + c_{m+1}}) _{i_{m}} \,
}\kern2cm\\
\cdot
{
 \frac {
  ({ \textstyle -1 - 2 a - n + b_{m} + b_{m+1} + c_{m} + c_{m+1}}) _{i_{m}}\,
    ({ \textstyle -a - n + c_{m+1}}) _{i_{m}-i_{m-1} } \,
     ({ \textstyle b_{m+1}}) _{i_{m-1}} \,
       }
   {
     ({ \textstyle -a - n + c_{m} + c_{m+1}}) _{i_{m}} \,
({i_{m}-i_{m-1} })! \,
     ({ \textstyle -1 - 2 a - n + b_{m} + b_{m+1} + c_{m} + c_{m+1}})
   _{i_{m-1}} }}\\
\cdot
\left( \prod_{k = 1}^{m-1}
        {\frac {            ({ \textstyle 1 + a - b_{k} - c_{k}}) _{i_{k}-i_{k-1}} \,
({ \textstyle b_{k+1}}) _{i_{k}} \,
            ({ \textstyle c_{k+1}}) _{i_{k}} }
          {({ i_{k}-i_{k-1} })! \,
            ({ \textstyle 1 + a - b_{k}}) _{i_{k}} \,
            ({ \textstyle 1 + a - c_{k}}) _{i_{k}} }} \right),
\end{multline}
o{\`u}, par d{\'e}finition, $i_0=0$ et o{\`u}, dans le cas $m=1$,
le produit vide doit {\^e}tre interpr{\'e}t{\'e}
comme valant 1.
\end{Theorem}

\begin{proof}
On commence avec l'identit{\'e} d'Andrews \eqref{eq:Amain}.
On reformule le membre de droite de cette identit{\'e}
en {\'e}crivant la somme
sur $i_m$ en notation hyperg{\'e}om{\'e}trique~:
\begin{multline}
\label{eq:4F3}
\frac{(1+a)_n\,(1+a-b_{m+1}-c_{m+1})_n}{(1+a-b_{m+1})_n\,(1+a-c_{m+1})_n}
\sum _{0\le i_1\le i_2\le\dots\le i_{m-1}\le n} ^{}\Bigg(
\frac{(-n)_{i_{m-1}}}{(b_{m+1}+c_{m+1}-a-n)_{i_{m-1}}}\kern1.7cm\\
\cdot
\frac{(b_{m+1})_{i_{m-1}}\,(c_{m+1})_{i_{m-1}}}
{(1+a-b_m)_{i_{m-1}}\,(1+a-c_m)_{i_{m-1}}}
\Bigg(\prod_{k=1}^{m-1}
\frac{(1+a-b_{k}-c_k)_{i_k-i_{k-1}}\,(b_{k+1})_{i_k}\,(c_{k+1})_{i_k}}
{(i_k-i_{k-1})!\,(1+a-b_k)_{i_k}\,(1+a-c_k)_{i_k}}\Bigg)\\
\cdot
{} _{4} F _{3} \!\left [ \begin{matrix} { c_{m+1}+i_{m-1},
b_{m+1}+i_{m-1},1+a-b_{m}-c_m,-n+i_{m-1}}\\
{1+a-b_m+i_{m-1},1+a-c_m+i_{m-1}, b_{m+1}+c_{m+1}-a-n+i_{m-1}
}\end{matrix} ; {\displaystyle 1}\right ] \Bigg).
\end{multline}
On applique maintenant une transformation liant deux s{\'e}ries $_4F_3$
balanc{\'e}es (voir\break \cite[(4.3.5.1)]{SlatAC})
\begin{multline}
\label{Sears}
{} _{4} F _{3} \!\left [ \begin{matrix}
{ a, b, c, -N}\\ { e, f, 1 + a + b + c - e -
   f - N}\end{matrix} ; {\displaystyle 1}\right ]  =
{\frac {({ \textstyle  e-a})_N\,({ f-a}) _{N}}  {({ \textstyle
e})_N\,({ f}) _{N}}}\\
\times
  {} _{4} F _{3} \!\left [ \begin{matrix} { -N, a,
 1 + a + c - e - f - N, 1 + a + b -
    e - f - N}\\ { 1 + a + b + c - e - f - N, 1 + a - e - N, 1 + a - f -
    N}\end{matrix} ; {\displaystyle 1}\right ]  ,
\end{multline}
o{\`u} $N$ est un entier positif, {\`a} la s{\'e}rie $_4F_3$ dans
\eqref{eq:4F3}. Ensuite, on {\'e}crit la (nouvelle) s{\'e}rie $_4F_3$
comme somme sur $i_m$. 
Enfin, apr{\`e}s quelques modifications, on arrive au membre
de droite de \eqref{eq:main}, ce qui d{\'e}montre le
th{\'e}or{\`e}me.
\end{proof}


\section{Quelques explications}\label{sec:explications}
Avant de continuer, il nous semble
b{\'e}n{\'e}fique d'indiquer
comment nous avons 
{\'e}t{\'e} men{\'e}s aux identit{\'e}s des
Th{\'e}or{\`e}mes~\ref{thm:A1}
et~\ref{thm:1}, et, en tout premier lieu, {\`a} certaines de leurs
sp{\'e}\-cia\-li\-sa\-tions les plus simples, donn{\'e}es
dans les diverses propositions du paragraphe \ref{sec:hyperharm}.
L'impulsion initiale a {\'e}t{\'e} donn{\'e}e par l'identit{\'e} suivante~:
\begin{multline}
\sum_{j=0}^n \frac{\dd}{\dd j}\(\frac n2-j\)\binom{n}{j}^4
\binom{n+j}{n}\binom{2n-j}{n}
\\=-\sum_{0\le i\le j\le n} (-1)^j\binom{n}{j}
\binom{n}{i}^2\binom{n+j}{n}\binom{n+j-i}{n}
=(-1)^{n+1}\sum_{j=0}^{n}\binom{n}{j}^2\binom{n+j}{n}^2,
\label{apmult}
\end{multline}
o{\`u} la 
{\og d{\'e}riv{\'e}e $\frac{\dd}{\dd j}$ avec $j$ entier\fg} doit
{\^e}tre entendue au sens suivant~:
\begin{align} \notag
\frac{\dd}{\dd j}\binom{n}{j}&=
\frac{\partial}{\partial \et}\(\frac{n!}
{\Ga(\et+1)\,\Ga(n-\et+1)}\)\bigg\vert_{\et=j}\\
&=\binom{n}{j}\(H_{n-j}-H_{j}\)
=\frac{\partial}{\partial \ep}\(\frac{n!}
{(1+\ep)_j\,(1-\ep)_{n-j}}\)\bigg\vert_{\ep=0},
\label{eq:d/dj}
\end{align}
o{\`u} $\Ga(x)$ d{\'e}signe la fonction gamma de $x$,
avec des identit{\'e}s similaires pour les autres
coefficients binomiaux.
On reconna{\^\i}t {\`a} gauche les nombres ${\bf a}_n$  et {\`a} droite les
nombres $a_n$, d{\'e}finis en
\eqref{eq:bolda_n} et \eqref{eq:a_n} respectivement~: nous avons d{\'e}j{\`a}
expliqu{\'e} la d{\'e}couverte de l'{\'e}galit{\'e} ${\bf a}_n=a_n$.
La somme double
au centre de~\eqref{apmult}
(qui a {\'e}t{\'e} d{\'e}j{\`a} mentionn{\'e} au paragraphe~\ref{sec:idee})
est apparue de fa{\c c}on tr{\`e}s
diff{\'e}rente~: dans~\cite{so1}, Sorokin a
montr{\'e} par des techniques d'approximants de Pad{\'e} que l'int{\'e}grale
$$
s_n=\int_0^1\int_0^1\int_0^1\frac{x^n(1-x)^ny^n(1-y)^nz^n(1-z)^n}
{(1-xy)^{n+1}(1-xyz)^{n+1}}
\,\dd x\,\dd y\,\dd z
$$
est {\'e}gale, pour tout entier $n\ge 0$,
{\`a} la combinaison lin{\'e}aire d'Ap{\'e}ry  $2a_n\zeta(3)-b_n$. En
d{\'e}veloppant le d{\'e}nominateur de l'int{\'e}grande de $s_n$ {\`a} l'aide du
th{\'e}or{\`e}me binomial et apr{\`e}s quelques manipulations {\'e}l{\'e}mentaires, on obtient
la s{\'e}rie double infinie
$$
s_n=n!\sum_{1\le \ell\le k}\frac{(k-\ell+1)_n(\ell-n)_n}{(k)^2_{n+1}(\ell)_{n+1}}.
$$
On d{\'e}compose alors en {\'e}l{\'e}ments simples le
sommande de $s_n$, ce qui permet de montrer que
$$
s_n=\bigg((-1)^n\sum_{0\le i\le j\le n} (-1)^j\binom{n}{j}
\binom{n}{i}^2\binom{n+j}{n}\binom{n+j-i}{n}\bigg)2\zeta(2,1)+\tilde{r}_n\zeta(2)+r_n,
$$
o{\`u}, par d{\'e}finition,
$\zeta(2,1)=\sum_{1\le \ell <  k}1/(k^{2}\ell)$ est un polyz{\^e}ta et o{\`u}
$r_n$,  $\tilde{r}_n$ sont des rationnels dont les expressions sont tr{\`e}s
compliqu{\'e}es\footnote{Il semble difficile de {\og voir \fg} directement
l'identit{\'e} combinatoire $\tilde{r}_n=0$ sur le d{\'e}veloppement
en {\'e}l{\'e}ments simples de $s_n$.
Voir le paragraphe \ref{ssec:nonterm} 
pour plus de d{\'e}tails concernant les s{\'e}ries multiples.}.
Or Euler a prouv{\'e}
que $\zeta(2,1)=\zeta(3)$ (voir les
{\it Opera  Omnia} \cite[p. 228]{euler})~: comme, conjecturalement, les
nombres 1,
$\zeta(2)$ et $\zeta(3)$ sont lin{\'e}airement ind{\'e}pendants sur
$\mathbb{Q}$, et compte-tenu du r{\'e}sultat de Sorokin,
il est l{\'e}gitime d'affirmer que $\tilde{r}_n=0$, 
que $r_n=-b_n$, et que
la deuxi{\`e}me {\'e}galit{\'e} de~\eqref{apmult} a lieu.
Pour prouver
ces affirmations, nous allons montrer directement que l'on a bien
l'identit{\'e}
\begin{equation*}
a_n=(-1)^n
\sum_{0\le i\le j\le n}
(-1)^j\binom{n}{j}\binom{n}{i}^2\binom{n+j}{n}\binom{n+j-i}{n},
\end{equation*}
(et on invoque ensuite l'irrationalit{\'e} de $\zeta(2)$).
Rappelons que les nombres $a_n$
s'expriment sous la forme
$$
a_n
=
{}_{4}F_3\!\left [ \begin{matrix}
{-n, -n, n+1,n+1}\\
{1,\;1,\; 1}\end{matrix} ;
{\displaystyle 1}\right ].
$$
Compte-tenu des binomiaux $\binom{n}{j}\binom{n+j-i}{n}$,
on peut {\'e}tendre {\`a} l'infini la somme double de~\eqref{apmult},
not{\'e}e ici $A_n$,
sans changer sa valeur, puis  mettre sous forme
hyperg{\'e}om{\'e}trique la somme int{\'e}rieure sur $i$~:
\begin{eqnarray*}
A_n
&=&\sum_{j=0}^{\ii} \sum_{i=0}^{\ii}(-1)^j
\binom{n}{j}\binom{n}{i}^2\binom{n+j}{n}\binom{n+j-i}{n}\\
&=&\sum_{j=0}^{\ii}\frac{(-n)_j\,(n+1)_j^2}{j!^3} \,
{}_{3}F_2
\!\left [ \begin{matrix}
{-n, -n, -j}\\
{1, -j-n}\end{matrix} ;
{\displaystyle 1}\right ].
\end{eqnarray*}
On applique ensuite {\`a} la s{\'e}rie ${}_{3}F_2$ l'une des transformations
de Thomae (voir  \cite[(3.1.1)]{GaRaAA})~:
$$
{}_{3}F_2
\!\left [ \begin{matrix}
{-n, a, b}\\
{d ,e}\end{matrix} ;
{\displaystyle 1}\right ]
=\frac{(e-b)_n}{(e)_n}\,
{}_{3}F_2
\!\left [ \begin{matrix}
{-n, b, d-a}\\
{d ,1+b-e-n}\end{matrix} ;
{\displaystyle 1}\right ].
$$
On obtient ainsi apr{\`e}s quelques manipulations~:
\begin{align*}
A_n &=
\sum_{j=0}^{\ii}\frac{(-j)_j\,(-n)_j\,(n+1)_j^2}{j!^3\,(-j-n)_j}\,
{}_{3}F_2
\!\left [ \begin{matrix}
{-n, -j, n+1}\\
{1,\; 1}\end{matrix} ;
{\displaystyle 1}\right ]\\
&=
\sum_{j=0}^{\ii}\frac{(-j)_j\,(-n)_j\,(n+1)_j^2}{j!^3\,(-j-n)_j}
\sum_{k=0}^j\frac{(-n)_k\,(-j)_k\,(n+1)_k}{k!^3}\\
&= \sum_{k=0}^{\ii}\sum_{j=k}^{\ii}
\frac{(-j)_j\,(-j)_k\,(-n)_j\,(-n)_k\,(n+1)_j^2\,(n+1)_k}
{j!^3\,k!^3\,(-j-n)_j}\\
&= \sum_{k=0}^{\ii}(-1)^k\frac{(-n)_k^2\,(n+1)_k^2}{k!^4}
\,{}_{2}F_1
\!\left [ \begin{matrix}
{k-n, n+k+1}\\
{k+1}\end{matrix} ;
{\displaystyle 1}\right ].
\end{align*}
La s{\'e}rie ${}_{2}F_1$ peut {\^e}tre somm{\'e}e par la formule
de Chu--Vandermonde (sous  forme
hyper\-g{\'e}o\-m{\'e}\-trique~; voir \cite[(1.7.7),
Appendix~(III.4)]{SlatAC})~:
\begin{equation} \label{vand}
{} _{2} F _{1} \!\left [ \begin{matrix} { a, -N}\\ { c}\end{matrix} ;
   {\displaystyle 1}\right ]  =
  {\frac{({ \textstyle c-a}) _{N} }{({ \textstyle c}) _{N} }},
\end{equation}
o{\`u} $N$ est un entier positif. D'o{\`u} finalement~:
\begin{multline*}
A_n=\sum_{k=0}^{\ii}(-1)^k
\frac{(-n)_k^2\,(-n)_{n-k}\,(n+1)_k^2}{k!^4\,(k+1)_{n-k}}
\\
=\;\frac{(-n)_n}{n!}\,
{}_{4}F_3
\!\left [ \begin{matrix}
{-n, -n, n+1,n+1}\\
{1,\;1,\;1}\end{matrix} ;
{\displaystyle 1}\right ]
=(-1)^n a_n.
\end{multline*}

\medskip
L'argument {\og divinatoire\fg} utilis{\'e} plus haut peut se
g{\'e}n{\'e}raliser. En effet,
au moyen de certains changements de variables birationnels,
Fischler a montr{\'e} dans \cite{fiscvar}
que les int{\'e}grales de Vasilyev \eqref{eq:vasint} prennent
la forme d'int{\'e}grales similaires {\`a} celle de  Sorokin,
ce qui est aussi une cons{\'e}quence d'un th{\'e}or{\`e}me de Zlobin \cite{zlo}.
Plus pr{\'e}cis{\'e}ment : si $E=2e\ge 2$ est pair, alors
$$
J_{E,n}=\int_{[0,1]^{E}}
\prod_{j=1}^e\frac{x_j^ny_j^n(1-x_j)^n(1-y_j)^n}
{(1-x_1y_1\cdots x_jy_j)^{n+1}}
\,\dd x_j\,\dd y_j
$$
et si $E=2e+1\ge 1$ est impair, alors
$$
J_{E,n}=\int_{[0,1]^{E}}
\prod_{j=1}^e\frac{x_j^ny_j^n(1-x_j)^n(1-y_j)^n}
{(1-x_1y_1\cdots x_jy_j)^{n+1}}
\,\dd x_j\,\dd y_j 
\cdot
\frac{t^n(1-t)^n}{(1-x_1y_1\cdots x_ey_et)^{n+1}}
\,\dd t.
$$
En effectuant un d{\'e}veloppement en s{\'e}rie des int{\'e}grandes
et en invoquant le th{\'e}or{\`e}me de Zudilin \cite{zu2} (liant
les int{\'e}grales de Vasilyev {\`a} certaines
s{\'e}ries tr{\`e}s bien {\'e}quilibr{\'e}es), on parvient {\`a} deviner
les identit{\'e}s de la Proposition~\ref{cor:A2} au paragraphe \ref{sec:hyperharm},
entre une somme finie tr{\`e}s bien {\'e}quilibr{\'e}e et une somme finie multiple.
Un argument similaire
(bas{\'e} sur la comparaison
du Th{\'e}or{\`e}me~1 de \cite{firi} et
d'une int{\'e}grale de
type Sorokin, due {\`a} Amoroso \cite{am}, au dernier paragraphe de \cite{firi})
nous a  permis de deviner les
identit{\'e}s de la Proposition~\ref{cor:A1} au
paragraphe~\ref{sec:hyperharm}.

Nous avons ensuite r{\'e}ussi
{\`a} trouver la d{\'e}monstration de ces deux
identit{\'e}s~: dans une version ant{\'e}rieure
\cite{KR1} du pr\'esent article,
il s'agissait essentiellement de la d{\'e}monstration du
Th{\'e}or{\`e}me~7 (au
paragraphe~9 de~\cite{KR1}) avec les valeurs des param{\`e}tres
sp{\'e}cialis{\'e}es comme dans la d{\'e}monstration de la
Proposition~\ref{cor:A2}, respectivement de la d{\'e}monstration du
Th{\'e}or{\`e}me~8 (au
paragraphe~10 de~\cite{KR1}) avec les valeurs des param{\`e}tres
sp{\'e}cialis{\'e}es comme dans la d{\'e}monstration de la
Proposition~\ref{cor:A1}. Or dans ces d{\'e}monstrations, on
n'utilisait que des identit{\'e}s hyperg{\'e}om{\'e}triques
tr{\`e}s classiques qui autorisent une tr{\`e}s grande libert{\'e}
dans le choix des param{\`e}tres~:
nous en avons donc mis le plus possible,
le r{\'e}sultat
{\'e}tant donn{\'e} par les tr\`es g\'en\'eraux 
Th{\'e}or{\`e}mes~7 et 8 de cette version
ant{\'e}rieure.

Nous n'avons r{\'e}alis{\'e} que beaucoup plus tard
que ce 
Th{\'e}or{\`e}me~8 \'etait en fait {\'e}quivalent {\`a} l'identit{\'e}
d'Andrews, vieille de trente ans et cit{\'e}e ici au
paragraphe~\ref{gigantesques},
tandis que ce Th{\'e}or{\`e}me~7 {\'e}tait simplement le r{\'e}sultat d'une
combinaison de l'identit{\'e} d'Andrews et de la transformation
\eqref{Sears} comme expliqu{\'e} ici dans la d{\'e}monstration du
Th{\'e}or{\`e}me~\ref{thm:1}. 
(Dans les deux cas, ces remarques valent 
{\`a} une inversion de l'ordre de sommation
dans la somme simple pr{\`e}s.)

Il est amusant de constater
que la cl{\'e} de la conjecture des d{\'e}nominateurs {\'e}tait cach{\'e}e
dans la recherche d'une preuve {\it directe} (cas $A=4$ de la Proposition~\ref{cor:A2})
de  l'{\'e}galit{\'e} ${\bf a}_n=(-1)^nA_n$,
et pas dans celle de  ${\bf a}_n=a_n$,
qui a pourtant lanc{\'e} toute cette recherche.


\section{Des identit{\'e}s hyperg{\'e}om{\'e}trico-harmoniques}
\label{sec:hyperharm}
Dans ce paragraphe, nous nous int{\'e}ressons au coefficient dominant\footnote{Ce n'est pas
${\bf  p}_{A,n}\((-1)^A\)$ qui est toujours nul!}
de la combinaison lin{\'e}aire issue de ${\bf S}_{n,A,B,C,1}\((-1)^A\)$,
c'est-{\`a}-dire
au coefficient devant $\zeta(A+C-1)$
(respectivement $\tilde\zeta(A+C-1)$) dans le d{\'e}veloppement \eqref{eq:SnABCr}.
Ce coefficient est
$$(-1)^C \binom {A+C-2}{A-2}
{\bf p}_{A-1,n}\((-1)^A\)$$
au paragraphe~\ref{ssec:conjdenom}.
En ignorant le signe et le coefficient binomial,
nous noterons\break
$P_n(A,B)=(-1)^{Bn+1}\mathbf  p_{A-1,n}\((-1)^A\)$, c'est-{\`a}-dire
explicitement, avec le symbole $\frac{\dd}{\dd j}$ entendu selon \eqref{eq:d/dj},
\begin{align*}
P_{n}(A,B)&=
\sum_{j=0}^n\frac{\dd}
{\dd j}\(\frac {n} {2}-j\)
{\binom nj}^{A}
{\binom {  n+j}{n}}^{B}
{\binom { 2 n-j}{n}}^{B}\\
&=\sum _{j=0} ^{n}\(\frac {n}
{2}-j\){\binom nj}^A{\binom {n+j}n}^B{\binom {2n-j}n}^B\\
&\kern2cm
\cdot\((A+B)H_{n-j}-(A+B)H_j+BH_{n+j}-BH_{2n-j}-\frac {1}
{\frac {n} {2}-j}\).
\end{align*}
Les propositions\footnote{Elles sont apparemment
nouvelles. Il est curieux que, ind{\'e}pendamment de nous,
Ahlgren a aussi d{\'e}couvert certaines identit{\'e}s
{\'e}quivalentes aux n{\^o}tres
pour $P_n(A,0)$ avec $A\in\{1, \ldots, 5\}$~: il
{\'e}tait motiv{\'e} par l'{\'e}tude \cite{AhOnAA} de congruences
satisfaites par des nombres d'Ap{\'e}ry g{\'e}n{\'e}ralis{\'e}s.
Paule et Schneider \cite{ps} ont remarqu{\'e} que l'on peut d{\'e}montrer
ces cinq identit{\'e}s en utilisant l'algorithme de Gosper--Zeilberger
\cite{ek,PeWZAA,ZeilAM,ZeilAV}. 
Depuis, Chu et De Donno \cite{ChuDonno} ont, comme nous, 
montr{\'e}
que ces cinq identit{\'e}s (et beaucoup d'autres) peuvent {\^e}tre trouv{\'e}es
comme cas limite de sommations classiques pour les s{\'e}ries
hyperg{\'e}om{\'e}triques.} de ce paragraphe montrent
que les nombres $P_n(A,B)$ sont entiers, ce qui est surprenant en
raison de la pr{\'e}sence des nombres harmoniques.
Sans le facteur $n/2-j$, qui
correspond au facteur tr{\`e}s bien {\'e}quilibrant $k+n/2$ de la s{\'e}rie
${\bf S}_{n,A,B,C,1}\((-1)^A\)$, aucun ne serait entier.
Remarquons que cette m{\^e}me conclusion d{\'e}coulera aussi des
Corollaires~\ref{cor:1} {\`a} \ref{cor:2a} (sauf pour de petites valeurs de
$A$ ou $B$) en faisant tendre $\ep$ vers~0. De plus, ces
corollaires nous permettront aussi de d{\'e}montrer la
Conjecture~\ref{cor:1} pour les autres coefficients ${\bf
  p}_{l,n}\((-1)^A\)$, alors que c'est impossible avec l'approche qui
m{\`e}ne aux propositions ci-dessous. Cependant, les identit{\'e}s exprim{\'e}es par ces
propositions sont beaucoup plus {\'e}l{\'e}gantes que celles que l'on d{\'e}duit
des corollaires cit{\'e}s, et elles ont donc un int{\'e}r{\^e}t
intrins{\`e}que. Enfin, elles seront essentielles dans
la d{\'e}monstration du Th{\'e}or{\`e}me~\ref{prop7gene}
au paragraphe~\ref{sec:demothmvaszud}.

La Proposition~\ref{cor:A1} se d{\'e}duit du Th{\'e}or{\`e}me~\ref{thm:A1},
alors que les Propositions~\ref{cor:A2} {\`a} \ref{cor:A4} r{\'e}sultent
du Th{\'e}or{\`e}me~\ref{thm:1},
sauf quelques cas particuliers
que l'on doit traiter {\og {\`a} la main\fg}.
Notons que la restriction $0\le 2Br<A$
(avec ici $r=1$) dispara{\^\i}t, puisque
nous donnons des expressions de $P_n(A,B)$
pour tous les entiers $A\ge 0$ et $B\ge 0$. Cependant, sans cette restriction,
la s{\'e}rie  ${\bf S}_{n,A,B,C,1}\((-1)^A\)$ peut {\^e}tre
divergente
et donc ne pas produire de combinaisons lin{\'e}aires de valeurs de
z{\^e}ta~:  c'est en particulier le cas pour $A=0$.

Dans tous les {\'e}nonc{\'e}s ci-dessous, le produit vide (par
exemple,
$\prod _{k=1} ^{m-1}$ pour $m=1$) doit {\^e}tre interpr{\'e}t{\'e} comme {\'e}gal
{\`a} 1.

\begin{Remark} On pourrait aussi d{\'e}montrer des identit{\'e}s
similaires pour un
entier $r$ quelconque, et pas seulement pour $r=1$.
\end{Remark}

\begin{Proposition} \label{cor:A1}
Soit
\begin{equation*}
P_n(A,0)=\sum_{j=0}^n\frac{\dd}{\dd j}
\(\frac {n}{2}-j\){\binom nj}^A.
\end{equation*}
Pour $A=2m+3\ge5$ impair, soit
$$p_n(A,0)=
\sum _{0\le i_1\le i_2\le\dots\le i_{m}\le n} ^{}
\prod _{k=1} ^{m}{\binom n{i_k}}^2\binom {n+i_{k+1}-i_k}n,$$
o{\`u}, par d{\'e}finition, $i_{m+1}=n$, et soient
$p_n(3,0)=1$ et $p_n(1,0)=(-1)^n$,
et pour $A=2m+2\ge 4$ pair, soient
$$p_n(A,0)=
\sum _{0\le i_1\le i_2\le\dots\le i_{m}\le n}
{\binom n{i_{m}}}^2\,
\prod _{k=1} ^{m-1}{\binom n{i_k}}^2\binom {n+i_{k+1}-i_k}n,
$$
$p_n(0,0)=-(n+1)$ et $p_n(2,0)=0$.
Alors pour tous entiers $A\ge0$ et $n\ge0$, on a
$P_n(A,0)=(-1)^{n+1}p_n(A,0)$.
\end{Proposition}

\begin{Remark}
Pour $A=4$, on a aussi $P_n(4,0)=(-1)^{n+1}\binom {2n}n$ puisque
$$p_n(4,0)=\sum _{0\le i_1\le n} ^{}
{\binom n{i_{1}}}^2=\binom {2n}n,
$$
en vertu de l'identit{\'e} de Chu--Vandermonde sous forme
binomiale (voir  par exemple \cite[paragraphe~5.1, (5.27)]{GrKPAA}).
\end{Remark}

\begin{Proposition} \label{cor:A2}
Soit
\begin{equation*}
P_n(A,1)=
\sum_{j=0}^{n}\frac{\dd}{\dd j}
\(\frac {n}{2}-j\){\binom nj}^A\binom{n+j}{n}\binom{2n-j}{n}.
\end{equation*}
Pour $A=2m+1\ge 3$ impair, soit
$$
p_n(A,1)=\sum_{0\le i_1\le i_2\le \cdots\le i_m\le n}
\binom{n}{i_{m}}^2\binom{n+i_{m}}{n}
\prod_{k=1}^{m-1}\binom{n}{i_k}^2 \binom{n+i_{k+1}-i_{k}}{n},
$$
et $p_n(1,1)=(-1)^n$, et pour $A=2m\ge 2$ pair, soit
$$
p_n(A,1)=\sum_{0\le i_1\le i_2\le \cdots\le i_{m}\le n} (-1)^{i_{m}}
\binom{n}{i_{m}}\binom{n+i_{m}}{n}
\prod_{k=1}^{m-1}\binom{n}{i_k}^2 \binom{n+i_{k+1}-i_{k}}{n},
$$
et $p_n(0,1)=-\binom{2n+1}{n+1}$.
Alors pour tous entiers
$A\ge0$ et $n\ge0$, on a $P_n(A,1)=\break
(-1)^{An+1}p_n(A,1)$.
\end{Proposition}

\begin{Remark}
Pour $A=2$, on a aussi $P_n(2,1)=(-1)^{n+1}$ puisque
$$p_n(2,1)=
\sum _{i_1=0} ^{n}(-1)^{i_1}\binom n{i_1}
\binom {n+i_1}n,$$
soit en termes hyperg{\'e}om{\'e}triques~:
$${} _{2} F _{1} \!\left [ \begin{matrix} { n+1,-n}\\ { 1}\end{matrix} ; \
{\displaystyle 1}\right ],$$
ce qui se simplifie {\`a} l'aide de l'identit{\'e} de Chu--Vandermonde
\eqref{vand} en $(-1)^n$.
\end{Remark}


\begin{Proposition} \label{cor:B1}
Soit
\begin{equation*}
P_n(0,B)=\sum_{j=0}^n\frac{\dd}{\dd j}
\(\frac {n}{2}-j\)\binom{n+j}{n}^B\binom{2n-j}{n}^B.
\end{equation*}
Pour $B\ge2$, soit
\begin{multline*}
p_n(0,B)=
\sum _{0\le i_1\le i_2\le\dots\le i_{B-1}\le n} ^{}
(-1)^{i_{B-1}+i_{B-2}}\binom {n+i_{B-1}-i_{B-2}} {i_{B-1}-i_{B-2}}
\binom {2n+1} {n-i_{B-1}}\\
\cdot\binom {n+i_{B-1}}n
\binom {3n+1} {n-i_{B-2}}
\prod _{k=1} ^{B-2}{\binom {2n-i_k}{n-i_{k-1}}}\binom {n+i_{k-1}}n
\binom {n+i_{k}}{n+i_{k-1}},
\end{multline*}
o{\`u}, par d{\'e}finition, $i_{0}=0$.
Alors pour tous entiers $B\ge2$ et $n\ge0$, on a
$P_n(0,B)=(-1)^{n+1}p_n(0,B)$.
\end{Proposition}

\begin{Proposition} \label{cor:B2}
Soit
\begin{equation*}
P_n(1,B)=\sum_{j=0}^n\frac{\dd}{\dd j}
\(\frac {n}{2}-j\)\binom nj\binom{n+j}{n}^B\binom{2n-j}{n}^B.
\end{equation*}
Pour $B\ge2$, soit
\begin{multline*}
p_n(1,B)=
\sum _{0\le i_1\le i_2\le\dots\le i_{B-1}\le n} ^{}
(-1)^{i_{B-1}+i_{B-2}}\binom {3n+1} {n-i_{B-1}}\\
\cdot
\binom {n+i_{B-1}-i_{B-2}} n
\binom {n+i_{B-1}}n
\prod _{k=1} ^{B-2}{\binom {n}{i_k-i_{k-1}}}\binom {2n-i_{k}}n
\binom {n+i_{k}}{n},
\end{multline*}
o{\`u}, par d{\'e}finition, $i_{0}=0$.
Alors pour tous entiers $B\ge2$ et $n\ge0$, on a
$P_n(1,B)=(-1)^{n+1}p_n(1,B)$.
\end{Proposition}

\begin{Proposition} \label{cor:A4}
Soit
\begin{equation*}
P_n(A,B)=
\sum_{j=0}^{n}\frac{\dd}{\dd j}
\(\frac {n}{2}-j\){\binom nj}^A{\binom {n+j}n}^B{\binom {2n-j}n}^B.
\end{equation*}
Pour $B\ge2$ et $A=2m+1\ge 3$ impair, soit
\begin{multline*}
p_n(A,B)=
\sum _{0\le i_1\le i_2\le\dots\le i_{m+B-1}\le n} ^{}
\kern-16pt(-1)^{i_{m+B-1}}
{\binom n{i_{m+B-1}}}\binom {n+i_{m+B-1}}n\binom {n+i_{m+B-1}-i_{m+B-2}}n\\
\cdot
\Bigg(\prod _{k=B} ^{m+B-2}{\binom n{i_k}}^2\binom
{n+i_{k}-i_{k-1}}n\Bigg)
\binom n{i_{B-1}}\binom {2n-i_{B-1}}n\binom n{i_{B-1}-i_{B-2}}\\
\cdot
\Bigg(\prod _{k=1} ^{B-2}{\binom {n+i_k}{n}}\binom {2n-i_k}n
\binom n{i_{k}-i_{k-1}}\Bigg),
\end{multline*}
o{\`u}, par d{\'e}finition, $i_0=0$,
et pour $B\ge2$ et $A=2m+2\ge 2$ pair, soit
\begin{multline*}
p_n(A,B)=
\sum _{0\le i_1\le i_2\le\dots\le i_{m+B-1}\le n} ^{}
\binom {n+i_{m+B-1}}n
\Bigg(\prod _{k=B} ^{m+B-1}{\binom n{i_k}}^2\binom
{n+i_{k}-i_{k-1}}n\Bigg)\\
\cdot
\binom n{i_{B-1}}\binom {2n-i_{B-1}}n\binom n{i_{B-1}-i_{B-2}}
\Bigg(\prod _{k=1} ^{B-2}{\binom {n+i_k}{n}}\binom {2n-i_k}n
\binom n{i_{k}-i_{k-1}}\Bigg).
\end{multline*}
Alors pour tous entiers $A\ge2$ et $B\ge2$, on a
$P_n(A,B)=(-1)^{(A+1)n+1} p_n(A,B)$.
\end{Proposition}

En sp{\'e}cialisant la Proposition~\ref{cor:A1} en $A=5$,
on obtient l'identit{\'e}
suivante (aussi d{\'e}montr{\'e}e par Paule et Schneider \cite{ps}) qui fait
intervenir les nombres $\al_n$ d{\'e}finis par~\eqref{eq:alpha_n}~:
\begin{equation*}
\sum_{j=0}^n \frac{\dd}{\dd j}\(\frac n2-j\) \binom{n}{j}^5
=(-1)^{n+1}\sum_{j=0}^{n}\binom{n}{j}^2\binom{n+j}{n}.
\label{eq:apubiq}
\end{equation*}
L'ubiquit{\'e} de ces nombres
est remarquable puisque, en vertu du cas $A=3$
de la Proposition~\ref{cor:A2}, on a {\'e}galement~:
\begin{equation*}
\sum_{j=0}^n \frac{\dd}{\dd j}\(\frac n2-j\)
\binom{n}{j}^3\binom{n+j}{n}\binom{2n-j}{n}
=(-1)^{n+1}\sum_{j=0}^{n}\binom{n}{j}^2\binom{n+j}{n}.
\label{q=al}
\end{equation*}
Cette {\'e}galit{\'e}  prouve que l'on a bien
$p_n=\al_n$, o{\`u} les nombres $p_n$ sont ceux d{\'e}finis par~\eqref{eq:q_n}.

\begin{proof}[D{\'e}monstration de la Proposition \ref{cor:A1}]
Le cas $A=0$ est trivial  et
pour les cas $A=1,2,3$, voir les
Lemmes~\ref{lem:1} {\`a}~\ref{lem:3} ci-dessous.
On suppose donc d{\'e}sormais que $A\ge4$.

Consid{\'e}rons d'abord le cas o{\`u} $A$ est impair, $A=2m+3$.
Dans le Th{\'e}or{\`e}me~\ref{thm:A1} on sp{\'e}cialise
$a=2\ep-n$,
$b_1=b_2=\dots=b_{m+1}=c_2=\dots=c_{m+1}=\ep-n$,
et enfin on fait tendre $c_1$ vers $\infty$.
On simplifie ensuite les deux membres
de \eqref{eq:Amain} et on fait tendre $\ep$ vers 0, ce
qui produit les nombres
harmoniques dans le sommande de la s{\'e}rie {\`a} gauche de
\eqref{eq:Amain}.

Si $A$ est pair, $A=2m+2$,
dans le Th{\'e}or{\`e}me~\ref{thm:A1} on fait tendre  $b_{m+1}$
et $c_1$ vers $\infty$,
puis on sp{\'e}cialise
$a=2\ep-n$ et
$b_1=b_2=\dots=b_{m}=c_2=\dots=c_{m+1}=\ep-n$.
On simplifie ensuite les deux membres de \eqref{eq:Amain} et
on fait finalement tendre $\ep$ vers 0.
\end{proof}

\begin{Lemma} \label{lem:1}
On a $P_n(1,0)=-1$.
\end{Lemma}
\begin{proof}[D{\'e}monstration]
Consid{\'e}rons l'expression
$$-\frac {1} {\ep}\sum _{j=0} ^{n}\(\frac {n}{2}-j+\ep\)
\frac {n!} {(1-\ep)_j\,(1+\ep)_{n-j}}.$$
La limite de cette expression quand $\ep\to0$ vaut
$P_n(1,0)$. D'autre part, on peut r{\'e}{\'e}crire cette
somme finie comme une diff{\'e}rence de deux sommes infinies, {\`a}
savoir
$$-\frac {1} {\ep}\sum _{j=0} ^{\infty}\(\frac {n} {2}-j+\ep\)
\frac {n!} {(1-\ep)_j\,(1+\ep)_{n-j}}
+
\frac {1} {\ep}\sum _{j=n+1} ^{\infty}\(\frac {n} {2}-j+\ep\)
\frac {n!} {(1-\ep)_j\,(1+\ep)_{n-j}},
$$
soit en notation hyperg{\'e}om{\'e}trique~:
\begin{multline*}
-\frac {\left( n+2\ep  \right) } {2\ep} \frac{
    n! }
     {({ \textstyle 1 + \ep}) _{n} }
  {} _{4} F _{3} \!\left [ \begin{matrix} { -2\,\ep - n, 1 -
      \ep - \frac{n}{2}, -\ep - n, 1}\\ {
      -\ep - \frac{n}{2}, 1 - \ep, -2\,\ep -
      n}\end{matrix} ; {\displaystyle -1}\right ]\\
+ \frac{\left( n - 2\ep +2 \right) \,
             n! }{2\,
      ({ \textstyle 1 - \ep}) _{n+1} }
  {} _{4} F _{3} \!\left [ \begin{matrix} { 2 - 2\,\ep + n,
         2 - \ep + \frac{n}{2}, 1 - \ep, 1}\\ { 1 -
         \ep + \frac{n}{2}, 2 - \ep + n, 2 -
         2\,\ep + n}\end{matrix} ; {\displaystyle -1}\right ].
\end{multline*}
Ces deux s{\'e}ries $_4F_3$ sont tr{\`e}s bien {\'e}quilibr{\'e}es et
on peut donc appliquer la sommation suivante (voir
\cite[(2.3.4.6); Appendix
(III.10)]{SlatAC})~:
\begin{equation} \label{4F3[-1]}
 {} _{4} F _{3} \!\left [ \begin{matrix} { a, 1 + \frac{a}{2}, b, c}\\ {
    \frac{a}{2}, 1 + a - b, 1 + a - c}\end{matrix} ; {\displaystyle -1}\right
    ]  = \frac{\Gamma({ \textstyle 1 + a - b}) \,
      \Gamma({ \textstyle 1 + a - c}) }{\Gamma({ \textstyle 1 + a}) \,
      \Gamma({ \textstyle 1 + a - b - c}) }.
\end{equation}
Apr{\`e}s quelques simplifications, on obtient ainsi
$$- \frac{n!}
    {2({ \textstyle 1 + \ep}) _{n} } -\frac{n! }
    {2({ \textstyle 1 - \ep}) _{n} } .
  $$
Apr{\`e}s avoir fait tendre $\ep$ vers 0, cela vaut $-1$.
\end{proof}

\begin{Lemma} \label{lem:2}
On a $P_n(2,0)=0$.
\end{Lemma}
\begin{proof}[D{\'e}monstration]
Consid{\'e}rons l'expression
$$-\frac {1} {\ep}\sum _{j=0} ^{n}\(\frac {n}
{2}-j+\ep\)\(\frac {n!}
{(1-\ep)_j\,(1+\ep)_{n-j}}\)^2.$$
La limite de cette expression quand $\ep$ tend vers 0 est
$P_n(2,0)$. D'autre part, on peut r{\'e}{\'e}crire cette
somme finie comme une diff{\'e}rence de deux sommes infinies, {\`a} savoir
\begin{equation*}
-\frac {1} {\ep}\sum _{j=0} ^{\infty}\(\frac {n} {2}-j+\ep\)
\(\frac {n!} {(1-\ep)_j\,(1+\ep)_{n-j}}\)^2
+
\frac {1} {\ep}\sum _{j=n+1} ^{\infty}\(\frac {n} {2}-j+\ep\)
\(\frac {n!} {(1-\ep)_j\,(1+\ep)_{n-j}}\)^2,
\end{equation*}
soit en notation hyperg{\'e}om{\'e}trique~:
\begin{multline*}
- \frac{\left( n+2\ep \right) \,
    n!^2}{2\ep\,
      ({ \textstyle 1 + \ep}) _{n} ^2}
  {} _{5} F _{4} \!\left [ \begin{matrix} { -2\,\ep - n, 1 -
       \ep - \frac{n}{2}, -\ep - n, -\ep - n,
       1}\\ { -\ep - \frac{n}{2}, 1 - \ep, 1 -
       \ep, -2\,\ep - n}\end{matrix} ; {\displaystyle
       1}\right ]\\
 -
  \frac{ \ep\left( n - 2\ep +2 \right) \,
              n!^2 }{2\,
      ({ \textstyle 1 - \ep}) _{n+1} ^2}
 {} _{5} F _{4} \!\left [ \begin{matrix} { 2 - 2\,\ep + n,
         2 - \ep + \frac{n}{2}, 1 - \ep, 1 -
         \ep, 1}\\ { 1 - \ep + \frac{n}{2}, 2 -
         \ep + n, 2 - \ep + n, 2 - 2\,\ep +
         n}\end{matrix} ; {\displaystyle 1}\right ] .
\end{multline*}
La deuxi{\`e}me expression s'annule quand on fait
tendre $\ep$ vers 0. La premi{\`e}re s{\'e}rie $_5F_4$ se traite
{\`a} l'aide de la sommation suivante (voir \cite[(2.3.4.5); Appendix (III.12)]{SlatAC})~:
\begin{multline} \label{eq:5F4}
 {} _{5} F _{4} \!\left [ \begin{matrix} { a, 1 + \frac{a}{2}, b, c, d}\\ {
    \frac{a}{2}, 1 + a - b, 1 + a - c, 1 + a - d}\end{matrix} ; {\displaystyle
    1}\right ] \\
 = \frac{\Gamma({ \textstyle 1 + a - b}) \,
      \Gamma({ \textstyle 1 + a - c}) \,\Gamma({ \textstyle 1 + a - d}) \,
      \Gamma({ \textstyle 1 + a - b - c - d}) }{\Gamma({ \textstyle 1 + a}) \,
      \Gamma({ \textstyle 1 + a - b - c}) \,
      \Gamma({ \textstyle 1 + a - b - d}) \,
      \Gamma({ \textstyle 1 + a - c - d}) }.
\end{multline}
Apr{\`e}s quelques simplifications, on obtient
$$ \frac{\ep\,n!\,
     \Gamma({ \textstyle n}) }{2\,
     ({ \textstyle 1 + \ep}) _{n} ^2},$$
qui tend vers 0 quand $\ep$ tend vers 0.
\end{proof}

\begin{Lemma} \label{lem:3}
On a $P_n(3,0)=(-1)^{n+1}$.
\end{Lemma}
\begin{proof}[D{\'e}monstration]
Consid{\'e}rons l'expression
$$-\frac {1} {\ep}\sum _{j=0} ^{n}\(\frac {n}
{2}-j+\ep\)\(\frac {n!}
{(1-\ep)_j\,(1+\ep)_{n-j}}\)^3.$$
La limite de cette expression quand $\ep$ tend vers 0 est
$P_n(3,0)$. D'autre part, on peut r{\'e}{\'e}crire cette
somme finie comme une diff{\'e}rence de deux sommes infinies, {\`a} savoir
\begin{equation*}
-\frac {1} {\ep}\sum _{j=0} ^{\infty}\(\frac {n} {2}-j+\ep\)
\(\frac {n!} {(1-\ep)_j\,(1+\ep)_{n-j}}\)^3
+
\frac {1} {\ep}\sum _{j=n+1} ^{\infty}\(\frac {n} {2}-j+\ep\)
\(\frac {n!} {(1-\ep)_j\,(1+\ep)_{n-j}}\)^3,
\end{equation*}
soit en notation hyperg{\'e}om{\'e}trique~:
\begin{multline*}
 -\frac{\left( n+2\ep  \right) \,
          n!^3}{2\ep\,
      ({ \textstyle 1 + \ep}) _{n} ^3}
  {} _{6} F _{5} \!\left [ \begin{matrix} { -2\ep - n, 1 -
       \ep - \frac{n}{2}, -\ep - n, -\ep - n,
       -\ep - n, 1}\\ { -\ep - \frac{n}{2}, 1 -
       \ep, 1 - \ep, 1 - \ep,
       -2\ep - n}\end{matrix} ; {\displaystyle -1}\right ]  \\
 -
  \frac{{\ep}^2\left( n - 2\ep + 2 \right) \,
          n!^3}{2\,
      ({ \textstyle 1 - \ep}) _{n+1} ^3}
 {} _{6} F _{5} \!\left [ \begin{matrix} { 2 - 2\ep + n, 2
       - \ep + \frac{n}{2}, 1 - \ep, 1 - \ep,
       1 - \ep, 1}\\ { 1 - \ep + \frac{n}{2}, 2 -
       \ep + n, 2 - \ep + n, 2 - \ep + n, 2 -
       2\ep + n}\end{matrix} ; {\displaystyle -1}\right ] .
\end{multline*}
La deuxi{\`e}me expression s'annule quand
$\ep$ tend vers 0. On applique cette fois-ci une transformation de la s{\'e}rie
$_6F_5$ restante en une s{\'e}rie $_3F_2$ (voir \cite[4.4(2)]{BailAA})~:
\begin{multline} \label{eq:6F5}
 {} _{6} F _{5} \!\left [ \begin{matrix} { a, 1 + \frac{a}{2}, b, c, d,
    e}\\ { \frac{a}{2}, 1 + a - b, 1 + a - c, 1 + a - d, 1 + a -
    e}\end{matrix} ; {\displaystyle -1}\right ]\\
  =
   \frac{\Gamma({ \textstyle 1 + a - d}) \,\Gamma({ \textstyle 1 + a - e}) }
      {\Gamma({ \textstyle 1 + a}) \,\Gamma({ \textstyle 1 + a - d -
e}) }
{} _{3} F _{2} \!\left [ \begin{matrix} { 1 + a - b - c, d, e}\\ {
       1 + a - b, 1 + a - c}\end{matrix} ; {\displaystyle 1}\right ].
\end{multline}
On obtient ainsi
$$-\frac{ n!^3}{2\,
     ({ \textstyle 1 +\ep}) _{n} ^3}
{} _{3} F _{2} \!\left [ \begin{matrix} { -\ep -
      n, 1, n+1}\\ { 1 -\ep, 1 -\ep}\end{matrix} ;
      {\displaystyle 1}\right ].
$$
Ici, on doit {\^e}tre prudent et  restreindre le
param{\`e}tre $\ep$ aux valeurs
n{\'e}gatives puisque, sinon, la s{\'e}rie $_3F_2$ ne converge pas. On
applique alors la transformation (voir \cite[Ex.~7, p.~98]{BailAA})~:
$$
{} _{3} F _{2} \!\left [ \begin{matrix} { a, b, c}\\ { d, e}\end{matrix} ;
   {\displaystyle 1}\right ]  =
  \frac {\Gamma( e)\, \Gamma(  d + e-a - b - c)}
{\Gamma( e -a)\, \Gamma( d +
    e -b - c}
  {} _{3} F _{2} \!\left [ \begin{matrix} { a,  d-b,  d-c}\\ { d,  d +
    e-b - c}\end{matrix} ; {\displaystyle 1}\right ].
$$
D'o{\`u}
$${\left( -1 \right) }^{n+1} \frac{
     \Gamma({ \textstyle 1 - \ep}) ^2\,
     n! ^3\,
     ({ \textstyle 2\ep+1}) _{n} }{
     \Gamma({ \textstyle 1 - 2\,\ep}) \,
     \Gamma({ \textstyle n+1}) \,
     ({ \textstyle 1 + \ep}) _{n} ^3}
{} _{3} F _{2} \!\left [ \begin{matrix} { -\ep - n,
      -\ep, -\ep - n}\\ { 1 - \ep,
      -2\ep - n}\end{matrix} ; {\displaystyle 1}\right ].
$$
On peut finalement faire tendre $\ep$ vers 0, ce qui donne
$(-1)^{n+1}$.
\end{proof}

\begin{proof}[D{\'e}monstration de la Proposition \ref{cor:A2}]
Pour les cas $A=0,1,2,3$, voir les
Lemmes~\ref{lem:6} {\`a}~\ref{lem:9} ci-dessous.
On suppose donc d{\'e}sormais que $A\ge4$.

Consid{\'e}rons d'abord le cas o{\`u} $A$ est impair, $A=2m+1$.
Dans le Th{\'e}or{\`e}me~\ref{thm:1} on sp{\'e}cialise
$a=2\ep-n$,
$b_1=b_2=\dots=b_{m}=c_2=\dots=c_{m}=\ep-n$,
$b_{m+1}=\ep-n$, $c_{m+1}=n+\ep+1$,
et enfin on fait tendre $c_1$ vers $\infty$.
On simplifie ensuite les deux membres
de \eqref{eq:main} et on fait tendre  $\ep$ vers 0, ce
qui produit les nombres
harmoniques dans le sommande de la s{\'e}rie {\`a} gauche de
\eqref{eq:main}. On trouve
l'identit{\'e} annonc{\'e}e.

Si $A$ est pair, $A=2m$,
dans le Th{\'e}or{\`e}me~\ref{thm:1},
on fait tendre $b_{m+1}$
et $c_1$ vers $\infty$,
puis on sp{\'e}cialise comme ci-dessus
$a=2\ep-n$,
$b_1=b_2=\dots=b_{m}=c_2=\dots=c_{m}=\ep-n$, $c_{m+1}=n+\ep+1$.
On simplifie ensuite les deux membres de \eqref{eq:main} et
on fait finalement tendre $\ep$ vers 0.
\end{proof}

\begin{Lemma} \label{lem:6}
On a $P_n(0,1)=-p_n(0,1)$.
\end{Lemma}

\begin{proof}[D{\'e}monstration]
La s{\'e}rie $P_n(0,1)$ est {\'e}gale {\`a}
$$\lim_{\ep\to0}\frac {1} {(-\ep)}
\sum _{j=0} ^{n}\(\frac n2-j+\ep\)
\binom {n+j-2\ep} {n}\binom {2n-j} {n}.$$
En termes hyperg{\'e}om{\'e}triques, cette
somme s'{\'e}crit
\begin{equation*}
\lim_{\ep\to0}\frac {1} {(-\ep)}
{\frac {\left( n+2 \ep  \right)
\,  ({ \textstyle 1 - 2 \ep}) _{n} \,
     ({ \textstyle n+1}) _{n} } {2 \,{n!^2}}}
     {} _{5} F _{4} \!\left [ \begin{matrix} { -2 \ep - n, 1 - \ep
      - {\frac n 2}, 1 - 2 \ep + n, 1, -n}\\ { -\ep - {\frac n 2},
      -2 n, -2 \ep - n, 1 - 2 \ep}\end{matrix} ; {\displaystyle
      1}\right ].
\end{equation*}
On traite la s{\'e}rie $_5F_4$ {\`a} l'aide de la sommation
\eqref{eq:5F4} et on arrive {\`a}
$$P_n(0,1)=-\lim_{\ep\to0}
{\frac {({ \textstyle 1 - 2 \ep}) _{n} \,({ \textstyle n+2}) _{n} }
   {{n!^2}}}=\binom {2n+1}n,$$
d'o{\`u} $P_n(0,1)=-p_n(0,1)$.
\end{proof}

\begin{Lemma} \label{lem:7}
On a $P_n(1,1)=(-1)^{n+1}p_n(1,1)$.
\end{Lemma}

\begin{proof}[D{\'e}monstration]
La s{\'e}rie $P_n(1,1)$ est {\'e}gale {\`a}
$$\lim_{\ep\to0}\frac {1} {(-\ep)}
\sum _{j=0} ^{n}\(\frac n2-j+\ep\)
\frac {n!} {j!\,({ \textstyle
  1 + 2 \ep}) _{n-j} }\binom {n+j-2\ep} {n}\binom {2n-j} {n}.
$$
En termes hyperg{\'e}om{\'e}triques, cette
somme s'{\'e}crit
\begin{equation*}
\lim_{\ep\to0}\frac {1} {(-\ep)}
{\frac {\left( n+2 \ep  \right)  \,
     ({ \textstyle 1 - 2 \ep}) _{n} \,({ \textstyle n+1}) _{n} }
   {2\,n! \,({ \textstyle 1 + 2 \ep}) _{n} }}
     {} _{4} F _{3} \!\left [ \begin{matrix} { -2 \ep - n, 1 - \ep
      - {\frac n 2}, 1 - 2 \ep + n, -n}\\ { -\ep - {\frac n 2},
      -2 n, 1 - 2 \ep}\end{matrix} ; {\displaystyle -1}\right ].
\end{equation*}
On traite la s{\'e}rie $_4F_3$ {\`a} l'aide de la sommation
\eqref{4F3[-1]} et on arrive {\`a}
$$P_n(1,1)=-\lim_{\ep\to0}
{\frac {({ \textstyle 1 - 2 \ep}) _{n} } {n! }}=-1,$$
d'o{\`u} $P_n(1,1)=(-1)^{n+1}p_n(1,1)$.
\end{proof}

\begin{Lemma} \label{lem:8}
On a $P_n(2,1)=-p_n(2,1)$.
\end{Lemma}

\begin{proof}[D{\'e}monstration]
Compte-tenu de la remarque juste apr{\`e}s la Proposition~\ref{cor:A2},
on a $p_n(2,1)=(-1)^n$.
D'autre part, la s{\'e}rie $P_n(2,1)$ est {\'e}gale {\`a}
$$\lim_{\ep\to0}\frac {1} {(-\ep)}
\sum _{j=0} ^{n}\(\frac n2-j+\ep\)
\binom {n} {j}\frac {n!} {(1-2\ep)_j\,(1+2\ep)_{n-j}}
\binom {n+j-2\ep}n \binom {2n-j}n.
$$
En termes hyperg{\'e}om{\'e}triques, cette
somme s'{\'e}crit
\begin{equation*}
\lim_{\ep\to0}\frac {1} {(-\ep)}
\frac{\left(\frac {n} {2}+ \ep  \right) \,
    ({ \textstyle 1 - 2 \ep}) _{n} \,(2n)! }{
    ({ \textstyle 1+2\ep}) _{n} \,n!^2}
    {} _{5} F _{4} \!\left [ \begin{matrix} { -2 \ep - n, 1 - \ep -
\frac{n}{2}, 1 - 2 \ep + n, -n, -n}\\
{ -\ep - \frac{n}{2}, -2 n, 1 - 2 \ep, 1 \
- 2 \ep}\end{matrix} ; {\displaystyle 1}\right ]  .
\end{equation*}
On traite la s{\'e}rie $_5F_4$ {\`a} l'aide de la sommation
\eqref{eq:5F4} et  on arrive {\`a}
$$P_n(2,1)=-\lim_{\ep\to0}  (-1)^n=(-1)^{n+1},$$
d'o{\`u} $P_n(2,1)=-p_n(2,1)$.
\end{proof}

\begin{Lemma} \label{lem:9}
On a $P_n(3,1)=(-1)^{n+1}p_n(3,1)$.
\end{Lemma}
\begin{proof}[D{\'e}monstration]
Nous commen{\c c}ons avec $p_n(3,1)$, qui, par
d{\'e}finition, vaut
$$
\sum _{i_1=0} ^{n}{\binom n{i_1}}^2\binom {n+i_1}n.
$$
On {\'e}crit cette s{\'e}rie comme la limite
$$\lim_{\ep\to0}\sum _{i_1=0} ^ {\infty}{\binom n
{i_1}} {\binom {i_1 + n} n} \frac{ n!}{(1+2\ep)_{i_1}\,
(1-\ep)_{n-i_1}},$$
soit en termes hyperg{\'e}om{\'e}triques~:
$$\lim_{\ep\to0}\frac{  n! }
{({ \textstyle 1-\ep}) _{n} }
{} _{3} F _{2} \!\left [ \begin{matrix} { 1 + n, \ep - n, -n}\\ { 1,
1 + 2 \ep}\end{matrix} ; {\displaystyle 1}\right ].$$
On applique {\`a} cette s{\'e}rie $_3F_2$ la transformation
(voir
\cite[(3.10.4), $q\to 1$, renvers{\'e}]{GaRaAA}~:
\begin{multline} \label{T3237}
{} _{3} F _{2} \!\left [ \begin{matrix} { x, y, -N}\\ { b, a}\end{matrix} ;
   {\displaystyle 1}\right ]  =
\frac {({ \textstyle a - x, a - y}) _{N}} {({ \textstyle a, a
- x - y}) _{N}}\\
\times
{} _{6} F _{5} \!\left [ \begin{matrix} { -a - N + x + y, 1 -
\frac{a}{2} - \frac{N}{2} + \frac{x}{2} + \frac{y}{2}, 1 - a - b - N + x +
y, x, y, -N}\\ { -\frac{a}{2} - \frac{N}{2} + \frac{x}{2} + \frac{y}{2}, b,
1 - a - N + y, 1 - a - N + x, 1 - a + x + y}\end{matrix} ; {\displaystyle
-1}\right ] ,
\end{multline}
o{\`u} $N$ est un entier positif.
On obtient ainsi l'expression suivante pour $p_n(3,1)$~:
\begin{multline} \label{eq:2a}
\lim_{\ep\to0}\frac{
    n! \,({ \textstyle 2 \ep - n}) _{n} \,
    ({ \textstyle 1 + \ep + n}) _{n} }{({ \textstyle 1-\ep}) _{n} \,
    ({ \textstyle 1+2\ep}) _{n} \,({ \textstyle \ep}) _{n} }
{} _{6} F _{5} \!\left [ \begin{matrix} { -\ep - n, 1 - \frac{\ep}{2} -
\frac{n}{2}, -\ep - n, 1 + n, \ep - n, -n}\\ { -\frac{\ep}{2} - \frac{n}{2}, 1,
-\ep - 2 n, 1 - 2 \ep, 1 - \ep}\end{matrix} ; {\displaystyle -1}\right
]  \\
=
{\left( -1 \right) }^{n}
\lim_{\ep\to0}\frac {2} {\ep}
\frac{({ \textstyle 1 - 2 \ep}) _{n} \,
   }{
    ({ \textstyle 1+2\ep}) _{ n} \,
     }
\sum _{j=0} ^{n}\({ \frac{n}{2} - j +\frac{\ep}{2}}\)
\binom nj
\frac { ({ \textstyle n - j-\ep  + 1}) _{j} } {({ \textstyle 1 - \ep}) _{j} }
\frac { {({ \textstyle n + \ep - j + 1}) _{j} }}
{({ \textstyle 1 - 2 \ep}) _{j} }\\
\cdot
\frac { ({ \textstyle 1 + j}) _{n} }
{ ({ \textstyle 1-\ep}) _{n} }
\frac {  ({ \textstyle 1 + \ep + n}) _{n-j} }
{  ({ \textstyle 1 + \ep}) _{n-j}
}.
\end{multline}
La somme sur $j$ est {\'e}gale {\`a} 0 pour $\ep=0$.
Donc en appliquant le th{\'e}or{\`e}me de l'H{\^o}pital,~on~a
\begin{multline*}
p_n(3,1)=(-1)^{n+1}
\sum _{j=0} ^{n}\(\frac {n}
{2}-j\){\binom nj}^3\binom {n+j}n\binom {2n-j}n\\
\cdot
\(2H_{n-j}-6H_j-2H_{2n-j}-\frac {1} {\frac {n} {2}-j}+C_1(n)\),
\end{multline*}
o{\`u} $C_1(n)$ est ind{\'e}pendant de $j$. Comme le terme
$$\(\frac {n}
{2}-j\){\binom nj}^3\binom {n+j}n\binom {2n-j}n$$
change de signe quand on remplace $j$ par $n-j$ (mais reste identique
autrement), cette derni{\`e}re
expression est aussi {\'e}gale {\`a}
$$
(-1)^{n+1}
\sum _{j=0} ^{n}\(\frac {n}
{2}-j\){\binom nj}^3\binom {n+j}n\binom {2n-j}n
\(4H_{n-j}-4H_j+H_{n+j}-H_{2n-j}-\frac {1} {\frac {n} {2}-j}\),
$$
ce qui est exactement la d{\'e}finition de $P_n(3,1)$, au signe pr{\`e}s.
\end{proof}

\begin{proof}[D{\'e}monstration de la Proposition \ref{cor:B1}]
Dans le Th{\'e}or{\`e}me~\ref{thm:1} on pose $m=B-1$ et on sp{\'e}cialise
$a=2\ep-n$,
$b_1=b_2=\dots=b_{B-1}=n+\ep+1$,
$c_1=c_2=\dots=c_{B-1}=\ep-n$, $b_{B}=\ep+1$, $c_{B}=n+\ep+1$.
On simplifie ensuite les deux membres
de \eqref{eq:main} et on fait tendre $\ep$ vers 0,
ce qui produit les nombres
harmoniques {\`a} dans le sommande de la s{\'e}rie {\`a} gauche de
\eqref{eq:main}.
\end{proof}

\begin{proof}[D{\'e}monstration de la Proposition \ref{cor:B2}]
Dans le Th{\'e}or{\`e}me~\ref{thm:1} on pose encore $m=B-1$ et
on sp{\'e}cialise
$a=2\ep-n$,
$b_1=b_2=\dots=b_{B-1}=n+\ep+1$,
$c_1=c_2=\dots=c_{B-1}=\ep-n$, $c_B=n+\ep+1$,
et finalement on fait tendre $b_B$
vers $\infty$.
On simplifie ensuite les deux membres
de l'identit{\'e} obtenue et on fait finalement
tendre $\ep$ vers 0,
ce qui produit les nombres
harmoniques dans le sommande de la s{\'e}rie {\`a} gauche de
\eqref{eq:main}.
\end{proof}

\begin{proof}[D{\'e}monstration de la Proposition \ref{cor:A4}]
Consid{\'e}rons d'abord le cas o{\`u} $A$ est pair, $A=2m+2$.
Dans le Th{\'e}or{\`e}me~\ref{thm:1}, on remplace $m$ par $m+B-1$.
Si $m\ge1$ on sp{\'e}cialise alors
$a=2\ep-n$,
$b_{1}=b_{2}=\dots=b_{B-1}=\ep-n$,
$c_{1}=c_{2}=\dots=c_{B-1}=n+\ep+1$,
$b_B=b_{B+1}=\dots=b_{m+B-1}=
c_B=c_{B+1}=\dots=c_{m+B-1}=\ep-n$,
$b_{m+B}=\ep-n$, $c_{m+B}=n+\ep+1$,
et si $m=0$
on sp{\'e}cialise
$a=2\ep-n$,
$b_{1}=b_{2}=\dots=b_{B}=n+\ep+1$,
$c_{1}=c_{2}=\dots=c_{B}=\ep-n$.
On simplifie ensuite les deux membres
de \eqref{eq:main} et on fait tendre $\ep$ vers 0,
ce qui produit les nombres
harmoniques dans le sommande de la s{\'e}rie {\`a} gauche de
\eqref{eq:main}.

Si $A$ est impair, $A=2m+1$,
dans le Th{\'e}or{\`e}me~\ref{thm:1}, on remplace encore $m$ par $m+B-1$,
on fait tendre $b_{m+B}$ vers $\infty$,
puis on sp{\'e}cialise comme ci-dessus
$a_1=2\ep-n$, $c_{m+B}=n+\ep+1$,
$b_{1}=b_{2}=\dots=b_{B-1}=\ep-n$,
$c_{1}=c_{2}=\dots=c_{B-1}=n+\ep+1$,
$b_B=b_{B+1}=\dots=b_{m+B-1}=
c_B=c_{B+1}=\dots=c_{m+B-1}=\ep-n$,
$c_{m+B}=n+\ep+1$.
On simplifie ensuite les deux membres de \eqref{eq:main} et
on fait finalement tendre $\ep$ vers 0.
\end{proof}


\section{Corollaires au Th{\'e}or{\`e}me \ref{thm:A1}}\label{CorollairesA}

Les corollaires suivants nous sont n{\'e}cessaires
pour d{\'e}montrer les parties des r{\'e}sultats principaux donn{\'e}s
au paragraphe~\ref{resultats} qui concernent le coefficient
{\og constant\fg}, c'est-{\`a}-dire les coefficients
${\bf p}_{0,C,n}\((-1)^A\)$, $p_{0,E,n}$ et ${\bf v}_n$ (voir
les paragraphes~\ref{sec:demothm1ii} 
{\`a} \ref{sec:Phi}).

\begin{Corollary} \label{cor:10}
Soient $A,B,r$ des nombres entiers tels que $A$ 
soit pair,
$A\ge2$, $B\ge0$ et $r\ge0$. Alors, on a
\begin{multline}
\sum_{j=k}^n
\(\frac {n} {2}-j+\ep\)
\(\frac {n!} { (1-\ep)_j \, (1+\ep)_{n-j}}\)^{A}
{\binom {  r n+j-\ep}{rn}}^{B}
{\binom { (r+1) n-j+\ep}{rn}}^{B}\\
=-\frac {k-\ep} {2}
\sum _{0\le i_1\le\dots\le i_{A/2+B}\le n-k} ^{}(-1)^{i_{A/2+B-1}}
\frac
{i_{A/2+B}!} {i_1!\,(i_2-i_1)!\cdots (i_{A/2+B}-i_{A/2+B-1})!}\\
\cdot
\(\prod _{j=1} ^{B}\frac {(1-\ep)_{rn+k+i_{j-1}}}
{(rn)!\,(1-\ep)_{k+i_{j-1}}}
\frac {(1+\ep)_{(r+1)n-k-i_{j}}}
{(rn-i_{j}+i_{j-1})!\,(1+\ep)_{n-k-i_{j-1}}}\)
\frac {n!}
{(1-\ep)_{k+i_{B}}\,(1+\ep)_{n-k-i_{B}}}\\
\cdot
\(\prod _{j=B+1} ^{A/2+B-1}\frac {n!}
{(1-\ep)_{k+i_{j}}\,(1+\ep)_{n-k-i_{j}}}
\frac {(n+i_j-i_{j-1})!}
{(1-\ep)_{k+i_{j}}\,(1+\ep)_{n-k-i_{j-1}}}\)\\
\times
\frac {n!\,(\ep)_{i_{A/2+B}-i_{A/2+B-1}}\,(1-2\ep)_{k+i_{A/2+B-1}}
\,(1-\ep)_{n-i_{A/2+B}-1}}
{(1-\ep)_n\,(1-2\ep)_{k-1}\,(1-\ep)_{k+i_{A/2+B}}\,
(1+\ep)_{n-k-i_{A/2+B}}},
\label{eq:dernier}
\end{multline}
o{\`u}, par d{\'e}finition, $i_0=0$ et o{\`u}, dans le cas $A=2$ ou $B=0$,
un produit vide doit {\^e}tre interpr{\'e}t{\'e}
comme valant 1.
\end{Corollary}
\begin{proof}[D{\'e}monstration]
En utilisant la notation hyperg{\'e}om{\'e}trique,
on {\'e}crit le membre de gauche
de la fa{\c c}on suivante~:
\begin{multline*}
\(\frac {n} {2}-k+\ep\)
\(\frac {n!} { (1-\ep)_k \, (1+\ep)_{n-k}}\)^{A}
{\binom {  r n+k-\ep}{rn}}^{B}
{\binom { (r+1) n-k+\ep}{rn}}^{B}\\
\times\lim_{\de\to0}\Bigg(
{} _{A+2B+5} F _{A+2B+4} \!\left [ \begin{matrix} -n+2k-2\ep,
-\frac {n} {2}+k-\ep+1,-n+k-\ep,\dots,\\
-\frac {n} {2}+k-\ep,k-\ep+1,\dots,\end{matrix} \right .\\
\left . \begin{matrix} rn+k-\ep+1,\dots,1,k-2\ep-\de+1,-n+k\\
-(r+1)n+k-\ep,\dots,-n+2k-2\ep,-n+k+\de,k-2\ep+1\end{matrix} ; {\displaystyle
       1}\right ]\Bigg),
\end{multline*}
o{\`u} $-n+k-\ep$ et $k-\ep+1$ apparaissent respectivement
$A+B$ fois, o{\`u} $rn+k-\ep+1$ et $-(r+1)n+k-\ep$
apparaissent respectivement
$B$ fois. On remarque l'apparition des termes {\og artificiels\fg}
$k-2\ep-\de+1,-n+k$ en haut et $-n+k+\de,k-2\ep+1$ en bas qui
disparaissent lorsque $\de\to0$: ils garantissent 
(plus pr{\'e}cis{\'e}ment, le terme $-n+k$ en haut)
que la s{\'e}rie est une somme finie.

On pose $m=A/2+B$, $a=-n+2k-2\ep$,
$b_1=b_2=\dots=b_B=-n+k-\ep$,
$c_1=c_2=\dots=c_B=rn+k-\ep+1$, 
$b_{B+1}=b_{B+2}=\dots=b_{A/2+B-1}=
c_{B+1}=c_{B+2}=\dots=c_{A/2+B-1}=-n+k-\ep$, 
$b_{A/2+B}=-n+k-\ep$, $c_{A/2+B}=k-2\ep-\de+1$,
$b_{A/2+B+1}=-n+k-\ep$, $c_{A/2+B+1}=1$, $N=n-k$
dans le Th{\'e}or{\`e}me~\ref{thm:A1}. Le r{\'e}sultat en d{\'e}coule
apr{\`e}s quelques manipulations imm{\'e}diates.
\end{proof}

\begin{Corollary} \label{cor:11}
Soient $A,B,r$ des nombres entiers tels que $A$ 
soit impair,
$A\ge 3$, $B\ge0$ et $r\ge0$. Alors, on a
\begin{multline}
\sum_{j=k}^n
\(\frac {n} {2}-j+\ep\)
\(\frac {n!} { (1-\ep)_j \, (1+\ep)_{n-j}}\)^{A}
{\binom {  r n+j-\ep}{rn}}^{B}
{\binom { (r+1) n-j+\ep}{rn}}^{B}\\
=-\frac {k-\ep} {2}
\sum _{0\le i_1\le\dots\le i_{(A+1)/2+B}\le n-k} ^{}(-1)^{i_{(A-1)/2+B}}
\frac
{i_{(A+1)/2+B}!} {i_1!\,(i_2-i_1)!\cdots (i_{(A+1)/2+B}-i_{(A-1)/2+B})!}\\
\cdot
\frac {n!} { (1-\ep)_k \, (1+\ep)_{n-k}}\\
\cdot
\(\prod _{j=2} ^{B+1}\frac {(1-\ep)_{rn+k+i_{j-1}}}
{(rn)!\,(1-\ep)_{k+i_{j-1}}}
\frac {(1+\ep)_{(r+1)n-k-i_{j}}}
{(rn-i_{j}+i_{j-1})!\,(1+\ep)_{n-k-i_{j-1}}}\)
\frac {n!}
{(1-\ep)_{k+i_{B+1}}\,(1+\ep)_{n-k-i_{B+1}}}\\
\cdot
\(\prod _{j=B+2} ^{(A-1)/2+B}\frac {n!}
{(1-\ep)_{k+i_{j}}\,(1+\ep)_{n-k-i_{j}}}
\frac {(n+i_j-i_{j-1})!}
{(1-\ep)_{k+i_{j}}\,(1+\ep)_{n-k-i_{j-1}}}\)\\
\times
\frac {n!\,(\ep)_{i_{(A+1)/2+B}-i_{(A-1)/2+B}}\,(1-2\ep)_{k+i_{(A-1)/2+B}}
\,(1-\ep)_{n-i_{(A+1)/2+B}-1}}
{(1-\ep)_n\,(1-2\ep)_{k-1}\,(1-\ep)_{k+i_{(A+1)/2+B}}\,
(1+\ep)_{n-k-i_{(A+1)/2+B}}},
\label{eq:dernier2}
\end{multline}
o{\`u}, par d{\'e}finition, $i_0=0$ et o{\`u}, dans le cas $A=3$ ou $B=0$,
un produit vide doit {\^e}tre interpr{\'e}t{\'e}
comme valant 1.
\end{Corollary}
\begin{proof}[D{\'e}monstration]
Cette d{\'e}monstration est compl{\`e}tement analogue {\`a} la
d{\'e}monstration pr{\'e}\-c{\'e}\-dente. 
Cette fois-ci,
on pose $m=(A+1)/2+B$, $a=-n+2k-2\ep$,
$b_1=b_2=\dots=b_{B+1}=-n+k-\ep$,
$c_2=\dots=c_{B+1}=rn+k-\ep+1$, 
$b_{B+2}=b_{B+3}=\dots=b_{(A-1)/2+B}=
c_{B+2}=c_{B+3}=\dots=c_{(A-1)/2+B}=-n+k-\ep$, 
$b_{(A+1)/2+B}=-n+k-\ep$, $c_{(A+1)/2+B}=k-2\ep-\de+1$,
$b_{(A+3)/2+B}=-n+k-\ep$, $c_{(A+3)/2+B}=1$, $N=n-k$
dans le Th{\'e}or{\`e}me~\ref{thm:A1},
et finalement on fait tendre $c_1$ vers $\infty$. 
Ensuite le r{\'e}sultat en d{\'e}coule
apr{\`e}s quelques manipulations imm{\'e}diates.
\end{proof}


\section{Corollaires au Th{\'e}or{\`e}me \ref{thm:1}}\label{Corollaires}

Les corollaires suivants nous sont n{\'e}cessaires
pour d{\'e}montrer les parties des r{\'e}sultats principaux donn{\'e}s
au paragraphe~\ref{resultats} qui concernent les coefficients
${\bf p}_{l,n}\((-1)^A\)$, $p_{l,E,n}$, pour $l\ge1$, et ${\bf u}_n$,
(voir les paragraphes~\ref{sec:demothm1i} {\`a} \ref{sec:demothmvaszud}).

\begin{Corollary} \label{cor:1}
Soient $A,B,r$ des nombres entiers tels que $A$ 
soit pair, $A\ge2$, $B\ge2$ et $r\ge0$.
Consid{\'e}rons
$$
S_{A,B,r}(n)=
\sum _{j=0} ^{n}
 \left( \frac{n}{2} - j + \ep\right)
   \left(
      \frac{n!}{( 1-\ep)_j\,
         (1+\ep)_{n-j}} \right)^{A}
   { \binom { r n+j-\ep} {r n} }^{B}\,
   { \binom {\left( r +1\right)n-j+\ep } {r n} }^{B}.
$$
Soit {\'e}galement
\begin{multline*}
s_{A,B,r}(n)=
\sum _{0\le i_1\le i_2\le\dots\le i_{A/2+B}\le n} ^{}
{\left( -1 \right) }^{i_{A/2+B-1} + i_{ A/2 + B}}
\frac{n!\,
     (1 + \ep )_{rn}\,(1 - \ep )_{rn}}{     (1 - \ep )_n\,
    {\left( r n \right) !}^2}\\
\cdot
\binom{(r+1)n + \ep + 1}{n - i_{A/2 +B}}
\binom{r n + i_{ A/2 + B}+\ep}{
 (r-1)  n  + i_{ A/2 + B}}
\\
\cdot
\frac{(r n + i_{A/2+B}+1)!}
   {( 1- \ep)_{r n + i_{A/2+B}+1}}
\binom{r n -\ep+ i_{A/2+B}- i_{A/2+B-1} }{
     i_{A/2+B}-i_{A/2+B-1} }\\
\cdot
\binom{(r+1)n - \ep + 1 }{
   r n + i_{A/2 + B-1}+1}
 \frac{ i_{A/2+B-1}!\,
     (1+2\ep)_{i_{A/2+B-1}}}{(1+ \ep  )_{i_{A/2+B-1}}\,
     (1+\ep  )_{i_{A/2+B-1}}} \binom{n + i_{A/2+B-1}- i_{A/2+B-2} }{
    i_{A/2+B-1}-i_{A/2+B-2} }\\
\cdot
\Bigg( \prod_{k = B}^{A/2 + B-2}
\binom{n +i_k- i_{k-1} }{
        i_k-i_{k-1} }
    \frac{{n!}}{    ( 1+ \ep )_{i_k}\,
   ( 1-\ep)_{n  - i_k}}
    \frac{{n!}}{
       ( 1+ \ep )_{i_k}\,
       ( 1-\ep)_{n  - i_k}}\Bigg)\\
\cdot
\binom{r n}{i_{B-1}-i_{B-2} }
\frac {n!} {    (1+ \ep )_{i_{B-1}}\,    (1-\ep)_{n  - i_{B-1}}}
 \frac{(1-\ep)_{(r+1) n - i_{B-1}}}
   {\left( r n \right) !\,( 1-\ep)_{n  - i_{B-1}}}
\\
\cdot
\left( \prod_{k = 1}^{B-2}\binom{r n}{i_k-i_{k-1}}
     \binom{(r+1)n  -
       i_k-\ep}{r n}\binom{r n +
i_k+\ep}{r n}\right),
\end{multline*}
o{\`u}, par d{\'e}finition, $i_0=0$ et o{\`u}, dans les cas $A=2$ ou $B=2$,
les produits vides doivent {\^e}tre interpr{\'e}t{\'e}s
comme valant 1.
Alors,
$$
S_{A,B,r}(n)=\ep\cdot(-1)^ns_{A,B,r}(n).
$$
\end{Corollary}

\begin{proof}[D{\'e}monstration]
On sp{\'e}cialise $m=A/2+B$,
$a=2\ep-n$,
$b_{1}=b_{2}=\dots=b_{B-1}=\ep-n$,
$c_{1}=c_{2}=\dots=c_{B-1}=rn+\ep+1$,
$b_B=b_{B+1}=\dots=b_{A/2+B-1}=
c_B=c_{B+1}=\dots=c_{A/2+B-1}=\ep-n$,
$b_{A/2+B}=2\ep+1$, $c_{A/2+B}=1$,
$b_{A/2+B+1}=\ep-n$ et $c_{A/2+B+1}=rn+\ep+1$ dans le
Th{\'e}or{\`e}me~\ref{thm:1}.
On {\'e}crit la s{\'e}rie hyperg{\'e}om{\'e}trique {\`a}
gauche de \eqref{eq:main} comme une somme sur $j$, et ensuite on
inverse l'ordre de sommation, c'est-{\`a}-dire, on remplace $j$ par $n-j$.
Apr{\`e}s
l'{\'e}limination de  certains facteurs redondants et
des simplifications des deux membres, on trouve
l'identit{\'e} annonc{\'e}e.
\end{proof}

\begin{Corollary} \label{cor:2}
Soient $A,B,r$ des nombres entiers tels que $A$ 
soit pair, $A\ge3$, $B\ge2$ et $r\ge0$.
Soient $c,d$ deux entiers tels que $c\ge 2$ et $d\ge1$.
Consid{\'e}rons
$$
S_{A,B,r}(n)=
\sum _{j=0} ^{n}
 \left( \frac{n}{2} - j + \ep\right)
   \left(
      \frac{n!}{( 1-\ep )_j\,
         (1+\ep )_{n-j}} \right)^{A}
   { \binom { r n+j-\ep} {r n} }^{B}\,
   { \binom {\left( r +1\right)n-j+\ep } {r n} }^{B},
$$
et soit
\begin{multline*}
s_{A,B,r}(n)=
\sum _{0\le i_1\le i_2\le\dots\le i_{(A+1)/2+B}\le n} ^{}
{\left( -1 \right) }^{i_{(A-1)/2+B} + i_{ (A+1)/2 + B}}
\frac{n!\,
     (1 + \ep )_{rn}\,(1 - \ep )_{rn}}{     (1 - \ep )_n\,
    {\left( r n \right) !}^2}\\
\cdot
\binom{(r+1)n + \ep + 1}{n - i_{(A+1)/2 +B}}
\binom{r n + i_{ (A+1)/2 + B}+\ep}{
 (r-1)  n  + i_{ (A+1)/2 + B}}\\
\cdot
\frac{(r n + i_{(A+1)/2+B}+1)!}
   {( 1- \ep)_{r n + i_{(A+1)/2+B}+1}}
\binom{r n -\ep+ i_{(A+1)/2+B}- i_{(A-1)/2+B} }{
     i_{(A+1)/2+B}-i_{(A-1)/2+B} }\\
\cdot
\binom{(r+1)n - \ep + 1 }{
   r n + i_{(A-1)/2 + B}+1}
 \frac{ i_{(A-1)/2+B}!\,
     (1+2\ep)_{i_{(A-1)/2+B}}}{(1+ \ep  )_{i_{(A-1)/2+B}}\,
     (1+\ep  )_{i_{(A-1)/2+B}}} \binom{n + i_{(A-1)/2+B}- i_{(A-3)/2+B} }{
    i_{(A-1)/2+B}-i_{(A-3)/2+B} }\\
\cdot
\Bigg( \prod_{k = B+1}^{(A-3)/2 + B}
\binom{n  +i_k- i_{k-1} }{
        i_k-i_{k-1} }
    \frac{{n!}}{    ( 1+ \ep )_{i_k}\,
   ( 1-\ep)_{n  - i_k}}
    \frac{{n!}}{
       ( 1+ \ep )_{i_k}\,
       ( 1-\ep)_{n  - i_k}}\Bigg)\\
\cdot
    \frac{{n!}}{    ( 1+ \ep )_{i_{B}}\,
   ( 1-\ep)_{n  - i_{B}}}
\binom{r n}{i_{B-1}-i_{B-2} }
\frac {n!} {    (1+ \ep )_{i_{B-1}}\,    (1-\ep)_{n  - i_{B-1}}}
 \frac{(1-\ep)_{(r+1) n - i_{B-1}}}
   {\left( r n \right) !\,(        i_{B}-i_{B-1})!
         ( 1-\ep)_{n  - i_{B}}}
\\
\cdot
\left( \prod_{k = 1}^{B-2}\binom{r n}{i_k-i_{k-1}}
     \binom{(r+1)n  -
       i_k-\ep}{r n}\binom{r n +
i_k+\ep}{r n}\right),
\end{multline*}
o{\`u}, par d{\'e}finition, $i_0=0$ et o{\`u}, dans les cas $A=3$ ou $B=2$,
les produits vides doivent {\^e}tre interpr{\'e}t{\'e}s
comme valant 1.
Alors,
$$
S_{A,B,r}(n)=\ep\cdot s_{A,B,r}(n).
$$
\end{Corollary}

\begin{proof}[D{\'e}monstration]
Dans le Th{\'e}or{\`e}me~\ref{thm:1}
on pose $m=(A+1)/2+B$, on fait tendre $c_B$ vers $\infty$,
et puis on sp{\'e}cialise
$a=2\ep-n$,
$b_{1}=b_{2}=\dots=b_{B-1}=\ep-n$,
$c_{1}=c_{2}=\dots=c_{B-1}=rn+\ep+1$,
$b_B=b_{B+1}=\dots=b_{(A-1)/2+B}=
c_{B+1}=\dots=c_{(A-1)/2+B}=\ep-n$,
$b_{(A+1)/2+B}=2\ep+1$, $c_{(A+1)/2+B}=1$,
$b_{(A+3)/2+B}=\ep-n$ et $c_{(A+3)/2+B}=rn+\ep+1$.
On {\'e}crit la s{\'e}rie hyperg{\'e}om{\'e}trique {\`a}
gauche de \eqref{eq:main} comme une somme sur $j$, et ensuite on
inverse l'ordre de sommation, c'est-{\`a}-dire, on remplace $j$ par $n-j$.
Apr{\`e}s
l'{\'e}limination de  certains facteurs redondants
et des simplifications des deux membres, on trouve alors
l'identit{\'e} annonc{\'e}e.
\end{proof}

\begin{Corollary} \label{cor:1a}
Soit $A$ et $r$ deux nombres entiers tels que
$A$ soit pair, $A\ge 2$ et $r\ge0$, et soient
\begin{equation*}
S_{A,1,r}(n)=
\sum _{j=0} ^{n}
 \left( \frac{n}{2} - j + \ep\right)
   \left(
      \frac{n!}{( 1-\ep )_j\,
         (1+\ep )_{n-j}} \right)^{A}
   { \binom { r n+j-\ep} {r n} }\,
   { \binom {\left( r +1\right)n-j+\ep } {r n} }
\end{equation*}
et
\begin{multline*}
s_{A,1,r}(n)=
\sum _{0\le i_1\le i_2\le\dots\le i_{A/2+1}\le n} ^{}
{\left( -1 \right) }^{i_{A/2} + i_{ A/2 +1}}
\frac{n!^2\,(1 - \ep )_{rn}}{     (1 - \ep )_n\,
(1 - \ep )_n\,
    {\left( r n \right) !}}\\
\cdot
\binom{(r+1)n + \ep + 1}{n - i_{A/2 +1}}
\binom{r n + i_{ A/2 + 1}+\ep}{
 (r-1)  n  + i_{ A/2 + 1}}
\frac{(r n + i_{A/2+1}+1)!}
   {( 1- \ep)_{r n + i_{A/2+1}+1}}
\binom{r n -\ep+ i_{A/2+1}- i_{A/2} }{
     i_{A/2+1}-i_{A/2} }\\
\cdot
\binom{(r+1)n - \ep + 1 }{
   r n + i_{A/2 }+1}
 \frac{ i_{A/2}!\,
     (1+2\ep)_{i_{A/2}}}{(1+ \ep  )_{i_{A/2}}\,
     (1+\ep  )_{i_{A/2}}} \binom{n + i_{A/2}- i_{A/2-1} }{
    i_{A/2}-i_{A/2-1} }\\
\cdot
\Bigg( \prod_{k = 1}^{A/2 -1}
\binom{n +i_k- i_{k-1} }{
        i_k-i_{k-1} }
    \frac{{n!}}{    ( 1+ \ep )_{i_k}\,
   ( 1-\ep)_{n  - i_k}}
    \frac{{n!}}{
       ( 1+ \ep )_{i_k}\,
       ( 1-\ep)_{n  - i_k}}\Bigg),
\end{multline*}
o{\`u}, par d{\'e}finition, $i_0=0$ et o{\`u}, dans le cas $A=2$,
le produit vide doit {\^e}tre interpr{\'e}t{\'e}
comme valant 1.
Alors,
$$
S_{A,1,r}(n)=\ep\cdot(-1)^ns_{A,1,r}(n).
$$
\end{Corollary}

\begin{proof}[D{\'e}monstration]
On sp{\'e}cialise $m=A/2+1$,
$a=2\ep-n$,
$b_1=b_2=\dots=b_{A/2}=
c_1=c_2=\dots=c_{A/2}=\ep-n$,
$b_{A/2+1}=2\ep+1$, $c_{A/2+1}=1$,
$b_{A/2+2}=\ep-n$, $c_{A/2+2}=rn+\ep+1$
dans le
Th{\'e}or{\`e}me~\ref{thm:1}. Puis, on simplifie les deux
membres
de l'identit{\'e} \eqref{eq:main} comme auparavant.
\end{proof}

\begin{Corollary} \label{cor:2a}
Soit $A$ et $r$ deux nombres entiers tels que
$A$ soit impair, $A\ge3$ et $r\ge0$, et soient
\begin{multline*}
S_{A,1,r}(n)=
\sum _{j=0} ^{n}
 \left( \frac{n}{2} - j + \ep\right)
   \left(
      \frac{n!}{( 1-\ep )_j\,
         (1+\ep )_{n-j}} \right)^{A}
   { \binom { r n+j-\ep} {r n} }\,
   { \binom {\left( r +1\right)n-j+\ep } {r n} }
\end{multline*}
et
\begin{multline*}
s_{A,1,r}(n)=
\sum _{0\le i_1\le i_2\le\dots\le i_{(A+3)/2}\le n} ^{}
{\left( -1 \right) }^{i_{(A+1)/2} + i_{ (A+3)/2}}
\frac{n!\,(1 - \ep )_{rn}}{     (1 - \ep )_n\,
    {\left( r n \right) !}}\\
\cdot
\binom{(r+1)n + \ep + 1}{n - i_{(A+3)/2}}
\binom{r n + i_{ (A+3)/2}+\ep}{
 (r-1)  n  + i_{ (A+3)/2}}\\
\cdot
\frac{(r n + i_{(A+3)/2}+1)!}
   {( 1- \ep)_{r n + i_{(A+3)/2}+1}}
\binom{r n -\ep+ i_{(A+3)/2}- i_{(A+1)/2} }{
     i_{(A+3)/2}-i_{(A+1)/2} }\\
\cdot
\binom{(r+1)n - \ep + 1 }{
   r n + i_{(A+1)/2 }+1}
 \frac{ i_{(A+1)/2}!\,
     (1+2\ep)_{i_{(A+1)/2}}}{(1+ \ep  )_{i_{(A+1)/2}}\,
     (1+\ep )_{i_{(A+1)/2}}} \binom{n + i_{(A+1)/2}- i_{(A-1)/2} }{
    i_{(A+1)/2}-i_{(A-1)/2} }\\
\cdot
\Bigg( \prod_{k = 2}^{(A-1)/2}
\binom{n  +i_k- i_{k-1} }{
        i_k-i_{k-1} }
    \frac{{n!}}{    ( 1+ \ep )_{i_k}\,
   ( 1-\ep)_{n  - i_k}}
    \frac{{n!}}{
       ( 1+ \ep )_{i_k}\,
       ( 1-\ep)_{n  - i_k}}\Bigg)\\
\cdot
    \frac{{n!}}{    ( 1+ \ep )_{i_1}\,
   ( 1-\ep)_{n  - i_1}}
    \frac{{n!}}{i_1!
       ( 1-\ep)_{n  - i_1}},
\end{multline*}
o{\`u}, par d{\'e}finition, $i_0=0$ et o{\`u}, dans le cas $A=3$,
le produit vide doit {\^e}tre interpr{\'e}t{\'e}
comme valant 1.
Alors,
$$
S_{A,1,r}(n)=\ep\cdot s_{A,1,r}(n).
$$
\end{Corollary}

\begin{proof}[D{\'e}monstration]
Dans le
Th{\'e}or{\`e}me~\ref{thm:1} on pose $m=(A+3)/2$,
on fait tendre $c_1$ vers $\infty$, puis on sp{\'e}cialise
$a=2\ep-n$,
$b_1=b_2=\dots=b_{(A+1)/2}=
c_1=c_2=\dots=c_{(A+1)/2}=\ep-n$,
$b_{(A+3)/2}=2\ep+1$, $c_{(A+3)/2}=1$,
$b_{(A+5)/2}=\ep-n$ et $c_{(A+5)/2}=rn+\ep+1$.
Finalement, on simplifie les deux
membres
de l'identit{\'e} \eqref{eq:main} comme auparavant.
\end{proof}

\section{Lemmes arithm{\'e}tiques}\label{sec:LemmesArithm}

Dans ce paragraphe, nous compilons plusieurs lemmes arithm{\'e}thiques
dont nous avons besoin 
au cours  des d{\'e}monstrations des
Th{\'e}or{\`e}mes~\ref{thm:2} {\`a} \ref{thm:C=3}, aux
paragraphes~\ref{sec:demothm1i} {\`a} \ref{sec:Phi}.
Les d{\'e}mons\-tra\-tions de ces lemmes sont tr{\`e}s similaires,
{\`a} quelques variations pr{\`e}s\footnote{Il serait possible
d'{\'e}noncer un th{\'e}or{\`e}me g{\'e}n{\'e}ral dont ces lemmes seraient
alors des corollaires, mais cela n'ajouterait rien {\`a} la transparence
des d{\'e}monstrations.}.
Pour un nombre premier $p$ donn{\'e},
notons $v_p(N)$ la valuation $p$-adique du nombre entier $N$.
On utilisera constamment le fait que
\begin{equation} \label{eq:vp} 
v_p(s!)=\sum_{j=1}^{\ii} \left[s/p^{\,j}\,\right]
\quad \textup{et}\quad v_p(\textup{d}_n)=
\left[\log_p(n)\right].
\end{equation}

Nous commen{\c c}ons avec un lemme de divisibilit{\'e} tr{\`e}s
{\'e}l{\'e}mentaire o{\`u} intervient le nombre
$$
\Phi_n=\underset{\{n/p\}\in[2/3,1[}{\prod_{p\text{ premier}}} p\
$$
de la Conjecture~\ref{conj2}.

\begin{Lemma}\label{lem:arith2}
Pour tout entier $n>0$ et tout entier $j$, $0\le j\le n$,
le nombre $\Phi_n$ divise
$$\binom{n+j}{n}\binom{2n-j}{n}.$$
\end{Lemma}

\begin{proof}[D{\'e}monstration]
On note $N=\{n/p\}$ et $J=\{j/p\}$
les parties fractionnaires de $n/p$ et $j/p$.
En posant
$$V_p=v_p\bigg(\binom{n+j}{n}\binom{2n-j}{n}\bigg),$$
on a donc
{\openup6pt
\begin{align} \label{eq:Vp}
V_p&\ge
\[\frac{n+j}{p}\]-\[\frac np\]-\[\frac jp\]
+\[\frac{2n-j}{p}\]-\[\frac np\]-\[\frac{n-j}{p}\]
\\
\notag
&\ge\[N+J\]+\[2N-J\]-\[N-J\].
\end{align}}%
Posons $U=\[N+J\]+\[2N-J\]-\[N-J\]$.
Comme clairement
$U\ge 0$, il suffit de montrer que, pour tout $J$,
si $N\ge 2/3$, alors $U\ge 1$.

Fixons $J$ quelconque et
supposons
que $N\ge 2/3$. On a toujours $-[N-J]\ge 0$ et $[N+J]\ge 0$~; comme
$N\ge 2/3$, on a aussi $[2N-J]\ge 0$. Pour que $U=0$, il faut alors
n{\'e}cessairement que $[2N-J]=0$ et $[N+J]=0$, donc que
$2N-J<1$ et $N+J<1$,
dont on d{\'e}duit que $N<2/3$~: contradiction. Donc pour tout $J$,
l'hypoth{\`e}se $N\ge2/3$ implique que  $U\ge 1$. Le lemme en d{\'e}coule.
\end{proof}

Le deuxi{\`e}me lemme concerne
les {\it briques {\'e}l{\'e}mentaires}, introduites par Zudilin dans
\cite{zu3} (et attribu{\'e}es {\`a} Nesterenko).
{\'E}tant donn{\'e}s des entiers positifs $\al$ et $\be$,
une brique {\'e}l{\'e}mentaire est
la fraction rationnelle (en $t$)
$$R(\al,\be;t)=\begin{cases} \dfrac {(t+\be)_{\al-\be}} {(\al-\be)!}
&\text{si }\al\ge \be,\\
\dfrac {(\be-\al-1)!} {(t+\al)_{\be-\al}}&\text{si }\al<\be.\end{cases}
$$
Les briques {\'e}l{\'e}mentaires satisfont aux propri{\'e}t{\'e}s
arithm{\'e}tiques suivantes
(voir \cite[paragraphe~7]{zu3} ou \cite[paragraphe~2]{nesbriques}).

\begin{Lemma} \label{lem:briques}
 Si $\al\ge \be$, alors pour tout entier $H\ge 0$,
\begin{equation*} 
\textup{d}_{\al-\be}^H \cdot\frac {1} {H!}\,
\frac {\partial ^H} {\partial t^H}R(\al,\be;t)\Big\vert_{t=-k}
\end{equation*}
est un nombre entier pour tout nombre entier $k$, et si
$\al_0\le \al<\be\le
\be_0$, alors
\begin{equation*} 
\textup{d}_{\be_0-\al_0-1}^H\cdot \frac {1} {H!}\,
\frac {\partial ^H} {\partial t^H}R(\al,\be;t)(t+k)\Big\vert_{t=-k}
\end{equation*}
est un nombre entier pour $k\in\{\al_0,\ldots, \be_0-1\}$.
\end{Lemma}

Les prochains lemmes concernent des {\it briques sp{\'e}ciales},
not{\'e}es $R_1(\dots)$, 
$R_2(\dots)$, 
\dots, $R_6(\dots)$.

Le lemme suivant est utilis{\'e} dans les d{\'e}monstrations
du Th{\'e}or{\`e}me~\ref{thm:2} (voir la
Propositions~\ref{cor:3} au paragraphe~\ref{sec:demothm1i} et
la Proposition~\ref{conj:4} au paragraphe~\ref{sec:demothm1ii}) et du
Th{\'e}or{\`e}me~\ref{prop:Phi} (au paragraphe~\ref{sec:demothmvaszud}).

\begin{Lemma} \label{lem:Faktor}
Pour des entiers $i, j, n,r\ge 0$, posons
\begin{equation*} \label{eq:R1}
R_1(n,i,j;\ep)=\frac {(rn)!} {n!^r}\binom {rn+i+\ep}{(r-1)n+j}
\end{equation*}
et
\begin{equation*} 
R_2(n,i,j;\ep)=\frac {(rn)!} {n!^r}
\frac {(1-\ep)_{rn}} {(rn)!}
\frac{(r n + j+1)!}
   {( 1- \ep)_{r n + j+1}}
\binom{r n -\ep+ j- i }{
     j-i }
\binom{(r+1)n - \ep + 1 }{
   r n + i+1}.
\end{equation*}
Alors, pour tout entier
$H\ge0$ et pour $0\le i,j\le n$, le nombre
\begin{equation*} 
\textup{d}_{n}^H\cdot \frac {1} {H!}\,
\frac {\partial ^H} {\partial \ep^H}
R_1(n,i,j;\ep)\Big\vert_{\ep=0}
\end{equation*}
est un nombre entier, et pour tout entier $H\ge0$
et pour $0\le i\le j\le n$, le nombre
\begin{equation} \label{eq:F2}
\textup{d}_{n}^H\cdot \frac {1} {H!}\,
\frac {\partial ^H} {\partial \ep^H}
R_2(n,i,j;\ep)\Big\vert_{\ep=0}
\end{equation}
est un nombre entier.
\end{Lemma}

\begin{proof}[D{\'e}monstration]
Nous explicitons la d{\'e}monstration pour $R_2(n,i,j;\ep)$, celle
pour $R_1(n,i,j;\ep)$ {\'e}tant compl{\`e}tement similaire.

On arrange les termes dans \eqref{eq:F2} pour obtenir l'expression
{\'e}quivalente
\begin{equation} \label{eq:F3}
\textup{d}_{n}^H\cdot \frac {1} {H!}\,
\frac {\partial ^H} {\partial \ep^H}\Bigg(
\frac {1} {n!^r}
\binom {rn+j+1}{j-i}
{(n-i+1-\ep)_{(r-1)n+j}} \,
{(rn+j+2-\ep)_{n-j} }
\Bigg)\Bigg\vert_{\ep=0}.
\end{equation}
Nous allons montrer que pour tous nombres entiers $f_1<f_2<\dots<f_H$
appartenant {\`a} $[n-i+1,rn+j-i]\cup[rn+j+2,(r+1)n+1]$, le nombre
\begin{equation} \label{eq:F4}
\textup{d}_{n}^H\cdot
\frac {1} {n!^r}
\binom {rn+j+1}{j-i}
{(n-i+1)_{(r-1)n+j}} \,
{(rn+j+2)_{n-j} }\frac {1} {f_1f_2\cdots f_H}
\end{equation}
est un nombre entier. {\'E}videmment, cela implique alors que
\eqref{eq:F3}, et donc aussi \eqref{eq:F2}, est un nombre entier.

Selon \eqref{eq:vp}, la valuation $p$-adique du nombre \eqref{eq:F4} est 
{\'e}gale {\`a}
\begin{multline} \label{eq:F5}
H\cdot[\log_p(n)]+
\sum _{\ell=1} ^{\infty}\bigg(
\left[\frac {(r+1)n+1} {p^\ell}\right]
-\left[\frac {rn+i+1} {p^\ell}\right]
-\left[\frac {n-i} {p^\ell}\right]\\
+\left[\frac {rn+j-i} {p^\ell}\right]
-\left[\frac {j-i} {p^\ell}\right]
-r\left[\frac {n} {p^\ell}\right]
\bigg)
-
\sum _{h=1} ^{H}v_p(f_h).
\end{multline}
Il nous faut d{\'e}montrer que cette quantit{\'e} est positive.

On {\'e}crit l'expression \eqref{eq:F5} sous la forme
\begin{multline} \label{eq:F7}
\sum _{\ell=1} ^{[\log_p(n)]}\bigg(
\left[\frac {(r+1)n+1} {p^\ell}\right]
-\left[\frac {rn+i+1} {p^\ell}\right]
-\left[\frac {n-i} {p^\ell}\right]
+\left[\frac {rn+j-i} {p^\ell}\right]
-\left[\frac {j-i} {p^\ell}\right]
-r\left[\frac {n} {p^\ell}\right]
\bigg)\\+
\sum _{\ell=[\log_p(n)]+1} ^{\infty}\bigg(
\left[\frac {(r+1)n+1} {p^\ell}\right]
-\left[\frac {rn+i+1} {p^\ell}\right]
+\left[\frac {rn+j-i} {p^\ell}\right]
\bigg)\\
-
\sum _{h=1} ^{H}\(v_p(f_h)-[\log_p(n)]\).
\end{multline}
La somme dans la premi{\`e}re ligne de \eqref{eq:F7} est {\'e}videmment
positive. La somme dans la troisi{\`e}me ligne sera maximale si
l'ensemble
des $f_i$ comprend tous le nombres divisibles par $p^{[\log_p(n)]}$
de l'ensemble de base
$E=[n-i+1,rn+j-i]\cup[rn+j+2,(r+1)n+1]$. Par cons{\'e}quent, on a
\begin{align*}
\sum _{h=1} ^{H}\(v_p(f_h)-[\log_p(n)]\)&\le
\sum _{f\in E,\,p^{[\log_p(n)]}\mid f} \(v_p(f)-[\log_p(n)]\)\\
&\le
\sum _{\ell=[\log_p(n)]+1} ^{\infty}\bigg(
\left[\frac {(r+1)n+1} {p^\ell}\right]
-\left[\frac {rn+j+1} {p^\ell}\right]
+\left[\frac {rn+j-i} {p^\ell}\right]
\bigg).
\end{align*}
Lorsque $i\le j$, la diff{\'e}rence entre la deuxi{\`e}me et la
troisi{\`e}me ligne de \eqref{eq:F7}, et ainsi l'expression compl{\`e}te
\eqref{eq:F7}, est bien positive.
\end{proof}

Le lemme suivant est utilis{\'e} dans la d{\'e}monstration de la
partie ii) du  Th{\'e}or{\`e}me \ref{thm:2}, voir
la Proposition~\ref{conj:4} au paragraphe~\ref{sec:demothm1ii},
et du Th{\'e}or{\`e}me~\ref{prop:Phi} au paragraphe~\ref{sec:demothmvaszud}.

\begin{Lemma} \label{lem:briques1}
Soit
$$R_3(n,k,m_1,m_2;\ep)=\frac {n!\,(1+\ep)_{m_1-m_2-1}\,(1-2\ep)_{k+m_{2}}
\,(1-\ep)_{n-m_{1}-1}}
{(1-\ep)_n\,(1-2\ep)_{k-1}\,(1-\ep)_{k+m_{1}}\,
(1+\ep)_{n-k-m_{1}}}.$$
Alors, pour tout $H\ge0$ et $0\le m_2<m_1\le n-k$, le nombre
\begin{equation} \label{eq:H0}
\textup{d}_n^{H+1}\cdot \frac {1} {H!}\,
\frac {\partial ^H} {\partial \ep^H}
R_3(n,k,m_1,m_2;\ep)\Big\vert_{\ep=0}
\end{equation}
est un nombre entier.
\end{Lemma}

\begin{proof}[D{\'e}monstration]
On arrange les termes dans \eqref{eq:H0} pour obtenir l'expression
{\'e}quivalente
\begin{equation} \label{eq:H1}
\textup{d}_n^{H+1}\cdot \frac {1} {H!}\,
\frac {\partial ^H} {\partial \ep^H}\bigg(
\frac {n!\,(k-2\ep)_{m_{2}+1}\,(1+\ep)_{m_1-m_2-1}}
{(n-m_1-\ep)_{m_1+1}\,(1-\ep)_{k+m_{1}}\,
(1+\ep)_{n-k-m_{1}}}\bigg)\Bigg\vert_{\ep=0}.
\end{equation}
Nous allons montrer que, pour tous nombres entiers $1\le f_1\le
f_2\le\dots\le f_H\le n$, le nombre
\begin{equation} \label{eq:H2}
\textup{d}_n^{H+1}
\frac {n!\,(k-2\ep)_{m_{2}+1}\,(1+\ep)_{m_1-m_2-1}}
{(n-m_1-\ep)_{m_1+1}\,(1-\ep)_{k+m_{1}}\,
(1+\ep)_{n-k-m_{1}}}
\frac {1} {f_1f_2\cdots f_H}
\end{equation}
est un nombre entier. {\'E}videmment, cela implique alors que
\eqref{eq:H1}, et donc aussi \eqref{eq:H0}, est un nombre entier.

On d{\'e}montre cet {\'e}nonc{\'e} en v{\'e}rifiant que la valuation
$p$-adique de \eqref{eq:H2} est positive pour tout nombre premier
$p$. Cette valuation $p$-adique est
{\'e}gale {\`a}
\begin{multline} \label{eq:H3}
(H+1)\cdot[\log_p(n)]+
\sum _{\ell=1} ^{\infty}\bigg(
\left[\frac {k+m_2} {p^\ell}\right]
+\left[\frac {m_1-m_2-1} {p^\ell}\right]
+\left[\frac {n-m_1-1} {p^\ell}\right]\\
-\left[\frac {k-1} {p^\ell}\right]
-\left[\frac {k+m_1} {p^\ell}\right]
-\left[\frac {n-k-m_1} {p^\ell}\right]
\bigg)
-
\sum _{h=1} ^{H}v_p(f_h).
\end{multline}
Si $p>n$, il est {\'evident} qu'elle est positive
parce que toutes les quantit{\'e}s qui
apparaissent sont nulles. D{\'e}sormais nous allons donc supposer que
$p\le n$.

En fait, les conditions sur $k,n,m_1,m_2$ impliquent que les termes
de la s{\'e}rie infinie dans \eqref{eq:H3} sont nuls
pour $\ell>[\log_p(n)]$. On peut donc
{\'e}crire l'expression \eqref{eq:H3} sous la forme
\begin{multline} \label{eq:H4}
[\log_p(n)]+\sum _{\ell=1} ^{[\log_p(n)]}\bigg(
\left[\frac {k+m_2} {p^\ell}\right]
+\left[\frac {m_1-m_2-1} {p^\ell}\right]
+\left[\frac {n-m_1-1} {p^\ell}\right]\\
-\left[\frac {k-1} {p^\ell}\right]
-\left[\frac {k+m_1} {p^\ell}\right]
-\left[\frac {n-k-m_1} {p^\ell}\right]
\bigg)
-
\sum _{h=1} ^{H}\(v_p(f_h)-[\log_p(n)]\).
\end{multline}
Comme, par d{\'e}finition, $1\le f_h\le n$ pour tout $h$, les termes
de la somme sur $h$ sont tous n{\'e}gatifs ou nuls. Par
cons{\'e}quent, il nous suffit {\`a} d{\'e}montrer que les sommandes de
la somme sur $\ell$ sont tous $\ge-1$.

Pour ce faire, on note $N=\{n/p^\ell\}$, $K=\{k/p^\ell\}$,
$M_1=\{m_{1}/p^\ell\}$, $M_2=\{m_{2}/p^\ell\}$ les parties
fractionnaire de $n/p^\ell$, $k/p^\ell$, $m_{1}/p^\ell$ et $m_{2}/p^\ell$.
Avec ces notations, le sommande devient
\begin{multline} \label{eq:H5}
\left[ K+M_2\right]
+\left[ M_1-M_2-\frac {1} {p^\ell}\right]
+\left[N-M_1-\frac {1} {p^\ell}\right]\\
-\left[K-\frac {1} {p^\ell}\right]
-\left[K+M_1\right]
-\left[N-K-M_1\right].
\end{multline}
Supposons que cette expression soit $\le-2$. Comme $-\left[K+M_1\right]
-\left[N-K-M_1\right]\ge0$, la deuxi{\`e}me ligne dans \eqref{eq:H5}
est toujours positive. Compte-tenu du fait que les rationnels $N$ et
$M_1$ ont comme d{\'e}nominateur $p^\ell$, le terme $\left[N-M_1-\frac
{1} {p^\ell}\right]$ est au moins $-1$. Donc, si l'expression
\eqref{eq:H5} est $\le-2$, alors on a forc{\'e}ment
$[K+M_2]=0$, $\left[M_1-M_2-\frac {1} {p^\ell}\right]=
\left[N-M_1-\frac {1} {p^\ell}\right]=-1$
et $\left[K-\frac {1} {p^\ell}\right]=0$, c'est-{\`a}-dire
\begin{align} \label{eq:H6a}
K+M_2&<1,\\
M_1-M_2-\frac {1} {p^\ell}&<0,\label{eq:H6b}\\
N-M_1-\frac {1} {p^\ell}&<0,\label{eq:H6d}\\
K-\frac {1} {p^\ell}&\ge0.\label{eq:H6c}
\end{align}
Nous distinguons deux cas~: si $K+M_1\ge1$, alors la condition
\eqref{eq:H6a} et le fait que les d{\'e}nominateurs de $M_1$ et $M_2$
sont $p^\ell$ impliquent que $M_1\ge M_2+\frac {1} {p^\ell}$. C'est
une contradiction avec \eqref{eq:H6b}. D'autre part, si $K+M_1<1$,
alors $[K+M_1]=0$ et donc $N-K-M_1\ge0$ (sinon
$-[N-K-M_1]=1$ et alors l'expression \eqref{eq:H5} est $\ge-1$).
En combinant avec \eqref{eq:H6d}, on obtient $K<\frac {1} {p^\ell}$,
ce qui contredit \eqref{eq:H6c}. L'expression \eqref{eq:H5} est donc
$\ge-1$, ce qui implique que \eqref{eq:H3} est bien positive.
\end{proof}

Pour d{\'e}montrer le Th{\'e}or{\`e}me~\ref{prop:Phi},
nous avons besoin de deux lemmes portant sur deux autres types de
briques sp{\'e}ciales, qui mettent en jeu le nombre
$$
\tilde\Phi_n=\underset{\{n/p\}\in[2/3,1[}{\prod_{p\text{
premier},\,\,p<n}} p.
$$

\begin{Lemma} \label{lem:binom1}
Pour tout nombre entier
$H\ge0$ et pour $0\le i\le n$, le nombre
\begin{equation} \label{eq:tildePhi}
\tilde\Phi_n^{-1}\textup{d}_n^{H}\cdot \frac {1} {H!}\,
\frac {\partial ^H} {\partial \ep^H}R_4(n,i;\ep)\Big\vert_{\ep=0},
\end{equation}
est un nombre entier, o{\`u}
$$
R_4(n,i;\ep)= \binom{ n +
i+\ep}{ n}\binom{2n -i-\ep}{ n}.
$$
\end{Lemma}

\begin{proof}[D{\'e}monstration]
On suit la d{\'e}marche de la
d{\'e}monstration du Lemme~\ref{lem:Faktor}. Il suffit de montrer
que, pour tous nombres entiers $1\le f_1\le f_2\le\dots\le f_H$ tels que
le multi-ensemble\footnote{Un {\it multi-ensemble} est un ensemble o{\`u}
l'on autorise des r{\'e}p{\'e}titions d'{\'e}l{\'e}ments.
Un multi-ensemble $\mathcal A$ est contenu dans un multi-ensemble
$\mathcal B$ si
pour tout $x\in \mathcal B$ le nombre de r{\'e}p{\'e}titions de $x$
dans $\mathcal A$
est au plus le nombre de r{\'e}p{\'e}titions de $x$ dans $\mathcal B$.}
$\{f_1,f_2,
\dots,f_H\}$ soit contenu dans le multi-ensemble
$[i+1,n+i]\cup[n-i+1,2n-i]$, le nombre
\begin{equation} \label{eq:E1}
\tilde\Phi_n^{-1}\textup{d}_{n}^H\cdot
\binom{ n +
i}{ n}\binom{2n  - i}{ n}\frac {1} {f_1f_2\cdots f_H}
\end{equation}
est un nombre entier. Pour $p$ premier avec $p<n$ et $\{n/p\}\in[2/3,1[$,
la valuation $p$-adique
du nombre \eqref{eq:E1} est
\begin{multline} \label{eq:E2}
-1+H\cdot[\log_p(n)]+
\sum _{\ell=1} ^{\infty}\bigg(
\left[\frac {n+i} {p^\ell}\right]
-\left[\frac {i} {p^\ell}\right]
-\left[\frac {n} {p^\ell}\right]\\
\kern5cm+
\left[\frac {2n-i} {p^\ell}\right]
-\left[\frac {n-i} {p^\ell}\right]
-\left[\frac {n} {p^\ell}\right]\bigg)-
\sum _{h=1} ^{H}v_p(f_h)\\
=-1+
\sum _{\ell=1} ^{[\log_p(n)]}\bigg(
\left[\frac {n+i} {p^\ell}\right]
-\left[\frac {i} {p^\ell}\right]
-\left[\frac {n} {p^\ell}\right]
+
\left[\frac {2n-i} {p^\ell}\right]
-\left[\frac {n-i} {p^\ell}\right]
-\left[\frac {n} {p^\ell}\right]\bigg)\\
+\sum _{\ell=[\log_p(n)]+1} ^{\infty}\bigg(
\left[\frac {n+i} {p^\ell}\right]
+
\left[\frac {2n-i} {p^\ell}\right]
\bigg)-
\sum _{h=1} ^{H}(v_p(f_h)-[\log_p(n)]).
\end{multline}
En proc{\'e}dant comme dans la d{\'e}monstration du Lemme~\ref{lem:Faktor},
on montre
que la deuxi{\`e}me ligne du membre de droite de \eqref{eq:E2} est bien
positive. D'autre part, le Lemme~\ref{lem:arith2} montre que le terme
pour $\ell=1$ de la somme sur la premi{\`e}re ligne est au moins 1. La
valuation $p$-adique \eqref{eq:E2} est donc bien positive.
Si $p\ge n$ ou
$\{n/p\}\notin[2/3,1[$, on proc{\`e}de de la m{\^e}me fa{\c c}on,
et, dans ce cas, les arguments du Lemme~\ref{lem:Faktor} suffisent,
c'est-{\`a}-dire,
qu'il n'est pas n{\'e}cessaire de les compl{\'e}ter par un lemme additionnel
puisque  $v_p(\tilde\Phi_n)=0$ dans ce cas.
\end{proof}

Le lemme suivant est utilis{\'e} dans la d{\'e}monstration du
Th{\'e}or{\`e}me~\ref{prop:Phi} au paragraphe~\ref{sec:demothmvaszud}.

\begin{Lemma} \label{lem:binom2}
Pour tous entiers
$n,c,H\ge 0$ et  $0\le j_0\le j_{1}\le \ldots\le j_{c+1} \le n$,
le nombre
\begin{equation} \label{eq:G5}
\tilde\Phi_n^{-1}\textup{d}_n^{H}\cdot
\frac {1} {H!}
\frac {\partial^{H}}
{\partial \ep^H }
R_5(n,j_{0},j_{1},\dots,j_{c+1};\ep)
\Big\vert_{\ep=0}
\end{equation}
est un nombre entier, o{\`u}
\begin{multline*}
R_5(n,j_{0},j_{1},\dots,j_{c+1};\ep)=\binom{2n + \ep + 1}{n - j_{c +1}}
\binom{ n + j_{ c +1}+\ep}{
  j_{ c +1}}
\frac{( n + j_{c+1}+1)!}
   {( 1- \ep)_{ n + j_{c+1}+1}}\\
\times
\binom{ n -\ep+ j_{c+1}- j_{c} }{
     j_{c+1}-j_{c} }
\binom{2n - \ep + 1 }{
    n + j_{c }+1}
 \binom{n +  j_{c}- j_{c-1} }{
    j_{c}-j_{c-1} }\\
\times
\Bigg( \prod_{k = 1}^{c -1}
\binom{n  +j_{k}- j_{k-1} }{
        j_{k}-j_{k-1} }
    \frac{{n!}}{
       ( 1+ \ep )_{j_{k}}\,
       ( 1-\ep)_{n  - j_{k}}}\Bigg)
 \frac{(1-\ep)_{2 n - j_0}}
   {n  !\,( 1-\ep)_{n  - j_0}}.
\end{multline*}
\end{Lemma}

\begin{proof}[D{\'e}monstration]
On suit de nouveau la d{\'e}marche de la
d{\'e}monstration du Lemme~\ref{lem:Faktor}.
Il suffit en fait de montrer
que, pour tous nombres entiers $1\le f_1\le f_2\le \dots\le f_{H}$
tels que le multi-ensemble
$\{f_h:f_h>n\}$ soit
contenu dans le
multi-ensemble
\begin{multline*}
E=[n+j_{c+1}+2,2n+1]\cup[n+1, n + j_{ c +1}]
\cup[n+1, n + j_{c+1}- j_{c}]\\
\cup[n+j_{c+1}+2,2n+1]
\cup[n+1,2n-j_0],
\end{multline*}
le nombre
\begin{multline} \label{eq:G1}
\tilde\Phi_n^{-1}\textup{d}_n^{H}\cdot
\binom{2n +  1}{n - j_{c +1}}
\binom{ n + j_{ c +1}}{
  j_{ c +1}}
\binom{ n + j_{c+1}- j_{c} }{
     j_{c+1}-j_{c} }
\binom{2n + 1 }{
    n + j_{c }+1}\\
\cdot
 \binom{n + j_{c}- j_{c-1} }{
    j_{c}-j_{c-1} }
\Bigg( \prod_{k = 1}^{c -1}
\binom{n +j_k- j_{k-1} }{
        j_k-j_{k-1} }\binom n{j_k} \Bigg)
 \binom{2 n - j_0}n
\frac {1} {f_1f_2\cdots f_{H}}\\[10pt]
=\tilde\Phi_n^{-1}\textup{d}_n^{H}\cdot
\binom{2n +  1}{n - j_{c +1}}
\binom{ n + j_{ c +1}}{
  j_{ c +1}}
\binom{ n + j_{c+1}- j_{c} }{
     j_{c+1}-j_{c} }
\binom{2n + 1 }{
    n + j_{c }+1}\kern2.6cm\\
\cdot
 \binom{n + j_{c}- j_{c-1} }{
    j_{c}-j_{c-1} }
\binom {2n-j_{c-1}}n
\Bigg( \prod_{k = 1}^{c -1}
\binom{n +j_k- j_{k-1} }{
        j_k }\binom {2n-j_{k-1}}{2n-j_k} \Bigg)
\frac {1} {f_1f_2\cdots f_{H}}
\end{multline}
est un nombre entier.
Pour $p$ premier avec $p<n$ et $\{n/p\}\in[2/3,1[$,
la valuation $p$-adique
du nombre \eqref{eq:G1} est
\begin{multline} \label{eq:G2}
-1+H\cdot[\log_p(n)]+
\sum _{\ell=1} ^{\infty}\Bigg(
\left[\frac {2n+1} {p^\ell}\right]
-\left[\frac {n+j_{c+1}+1} {p^\ell}\right]
-\left[\frac {n-j_{c+1}} {p^\ell}\right]\\
+\left[\frac {n+j_{c+1}} {p^\ell}\right]
-\left[\frac {j_{c+1}} {p^\ell}\right]
-\left[\frac {n} {p^\ell}\right]
+\left[\frac {n + j_{c+1}- j_{c} } {p^\ell}\right]
-\left[\frac {j_{c+1}-j_{c}} {p^\ell}\right]
-\left[\frac {n} {p^\ell}\right]\\
+\left[\frac {2n+1} {p^\ell}\right]
-\left[\frac {n + j_{c }+1} {p^\ell}\right]
-\left[\frac {n-j_{c}} {p^\ell}\right]
+\left[\frac {n + j_{c}- j_{c-1}} {p^\ell}\right]
-\left[\frac {j_{c}-j_{c-1} } {p^\ell}\right]
-\left[\frac {n} {p^\ell}\right]\\
+\left[\frac {2n-j_{c-1}} {p^\ell}\right]
-\left[\frac {n} {p^\ell}\right]
-\left[\frac {n-j_{c-1}} {p^\ell}\right]
+
\sum _{h=1} ^{c-1}\bigg(\left[\frac {n +j_h- j_{h-1}} {p^\ell}\right]
-\left[\frac {j_h} {p^\ell}\right]
-\left[\frac {n-j_{h-1}} {p^\ell}\right]\\
+\left[\frac {2n-j_{h-1}} {p^\ell}\right]
-\left[\frac {2n-j_h} {p^\ell}\right]
-\left[\frac {j_h-j_{h-1}} {p^\ell}\right]\bigg)\Bigg)
-
\sum _{h=1} ^{H}v_p(f_h)\\
=-1+
\sum _{\ell=1} ^{[\log_p(n)]}U(\ell)
+\sum _{\ell=[\log_p(n)]+1} ^{\infty}U(\ell)-
\sum _{h=1} ^{H}(v_p(f_h)-[\log_p(n)]),
\end{multline}
o{\`u} $U(\ell)$ d{\'e}note le sommande de la somme sur $\ell$. Notons
que, si le sommande
$v_p(f_h)-[\log_p(n)]$ de la somme sur $h$ est strictement positif,
alors forcement $f_h>n$~: comme
 le multi-ensemble $\{f_h:f_h>n\}$ est contenu
dans le multi-ensemble $E$, des arguments analogues {\`a} ceux du
Lemme~\ref{lem:Faktor} montrent alors que
\begin{multline*}
\sum _{h=1} ^{H}(v_p(f_h)-[\log_p(n)])\le
\sum _{\ell=[\log_p(n)]+1} ^{\infty}\Bigg(
2\left[\frac {2n+1} {p^\ell}\right]
-2\left[\frac {n+j_{c+1}+1} {p^\ell}\right]\\
\kern5cm
+\left[\frac {n + j_{ c +1}} {p^\ell}\right]
+\left[\frac { n + j_{c+1}- j_{c}} {p^\ell}\right]
+\left[\frac {2n-j_0} {p^\ell}\right]
\Bigg)\\
\le\sum _{\ell=[\log_p(n)]+1} ^{\infty}U(\ell).\kern7.8cm
\end{multline*}
Il nous reste donc {\`a} d{\'e}montrer que $\sum _{\ell=1}^{[\log_p(n)]}U(\ell)\ge1$.
Compte-tenu de la restriction $p<n$, cette somme est non-vide~:
comme tous les $U(\ell)$ sont visiblement positifs,
il suffit donc simplement de d{\'e}montrer que $U(1)\ge 1$.
En oubliant plusieurs termes qui sont toujours positifs, nous allons
d{\'e}montrer l'in{\'e}galit{\'e} plus forte
\begin{multline} \label{eq:G3}
\left[\frac {n+j_{c+1}} {p}\right]
-\left[\frac {j_{c+1}} {p}\right]
-\left[\frac {n} {p}\right]
+\left[\frac {n + j_{c+1}- j_{c} } {p}\right]
-\left[\frac {j_{c+1}-j_{c}} {p}\right]
-\left[\frac {n} {p}\right]\\
+\left[\frac {2n+1} {p}\right]
-\left[\frac {n + j_{c }+1} {p}\right]
-\left[\frac {n-j_{c}} {p}\right]
+\left[\frac {n + j_{c}- j_{c-1}} {p}\right]
-\left[\frac {j_{c}-j_{c-1} } {p}\right]
-\left[\frac {n} {p}\right]\\
+\left[\frac {2n-j_{c-1}} {p}\right]
-\left[\frac {n} {p}\right]
-\left[\frac {n-j_{c-1}} {p}\right]
\ge1.
\end{multline}
Pour ce faire, on d{\'e}finit $N=\{n/p\}$, $J_1=\{j_{c-1}/p\}$,
$J_2=\{j_{c}/p\}$, $J_3=\{j_{c+1}/p\}$.
L'expression {\`a} gauche de \eqref{eq:G3} devient alors
\begin{multline} \label{eq:G4}
\left[N+J_3\right]
+\left[N+J_3-J_2\right]
-\left[J_3-J_2\right]
+\left[2N+\tfrac {1} {p}\right]
-\left[N+J_2+\tfrac {1} {p}\right]
-\left[N-J_2\right]\\
+\left[N+J_2-J_1\right]
-\left[J_2-J_1\right]
+\left[2N-J_1\right]
-\left[N-J_1\right].
\end{multline}
Supposons que cette expression soit nulle.
Puisque, par hypoth{\`e}se,  $N\ge 2/3$, on a
$[2N+\tfrac {1} {p}]=1$  et, de plus,
\begin{multline*}
\left[N+J_3\right]\ge0,\quad\!
\left[N+J_3-J_2\right]
-\left[J_3-J_2\right]\ge0,\quad\!
-\left[N+J_2+\tfrac {1} {p}\right]\ge-1,\quad\!
-\left[N-J_2\right]\ge0,\\
\left[N+J_2-J_1\right]
-\left[J_2-J_1\right]\ge0,\quad
\left[2N-J_1\right]\ge0\quad \text {et}\quad
-\left[N-J_1\right]\ge0.
\end{multline*}
Pour que l'expression \eqref{eq:G4} soit nulle, il faut donc que
toutes ces in{\'e}galit{\'e}s soient en fait des {\'e}galit{\'e}s.
En particulier, on a n{\'e}cessairement
\begin{equation} \label{eq:G6}
N+J_3<1,\quad N\ge J_2, \quad 2N-J_1<1 \quad \text {et}\quad N\ge J_1.
\end{equation}
L'expression $\left[N+J_3-J_2\right]
-\left[J_3-J_2\right]$ est nulle si et seulement si soit $J_3\ge J_2$
(et $N+J_3-J_2<1$)
soit $N+J_3-J_2<0$. De m{\^e}me, l'expression $\left[N+J_2-J_1\right]
-\left[J_2-J_1\right]$ est nulle si et seulement si soit $J_2\ge J_1$
(et $N+J_2-J_1<1$)
soit $N+J_2-J_1<0$.

Si $N+J_2-J_1<0$, il s'ensuit de \eqref{eq:G6} que $N+J_2-J_1\ge
J_2\ge0$~: contradiction.
Le cas $N+J_3-J_2<0$ m{\`e}ne {\`a} une contradiction analogue {\`a} cause
de l'in{\'e}galit{\'e} $J_2\le N$ dans \eqref{eq:G6}.

Il reste la possibilit{\'e} que $J_3\ge J_2\ge J_1$. Dans ce cas,
en utilisant de nouveau les in{\'e}galit{\'e}s dans \eqref{eq:G6}, on a
$$2N-1<J_1\le J_2\le J_3<1-N,$$
c'est-{\`a}-dire, $N<2/3$, ce qui est de nouveau contradictoire.
On a donc bien $U(1)\ge 1$.

Si $p\ge n$ ou
$\{n/p\}\notin[2/3,1[$, on proc{\`e}de de la m{\^e}me fa{\c c}on,
et, dans ce cas, de nouveau les arguments du Lemme~\ref{lem:Faktor} suffisent,
c'est-{\`a}-dire qu'il n'est pas n{\'e}cessaire de les compl{\'e}ter par des
consid{\'e}rations additionnelles (comme, par exemple,
celles ci-dessus pour d{\'e}montrer
que $U(1)\ge1$),
puisque  $v_p(\tilde\Phi_n)=0$ dans ce cas.
\end{proof}

Le lemme suivant aussi est utilis{\'e} dans la d{\'e}monstration du
Th{\'e}or{\`e}me~\ref{prop:Phi} au paragraphe~\ref{sec:demothmvaszud}.

\begin{Lemma} \label{lem:binom3}
Pour tous entiers
$H\ge 0$, $B\ge0$, $0\le k\le n$, et\/
$0\le i_1\le i_2 \le\ldots\le i_B \le n-k$,
le nombre
\begin{equation} \label{eq:G10}
\tilde\Phi_n^{-B+1}\textup{d}_n^{H}\cdot
\frac {1} {H!}
\frac {\partial^{H}}
{\partial \ep^H }
R_6(n,i_{1},i_{2},\dots,i_{B};\ep)
\Big\vert_{\ep=0}
\end{equation}
est un nombre entier, o{\`u}
\begin{multline*}
R_6(n,i_{1},i_{2},\dots,i_{B};\ep)=
\frac
{i_{B}!} {i_1!\,(i_2-i_1)!\cdots (i_{B}-i_{B-1})!}\\
\times
\prod _{j=1} ^{B}\frac {(1-\ep)_{n+k+i_{j-1}}}
{n!\,(1-\ep)_{k+i_{j-1}}}
\frac {(1+\ep)_{2n-k-i_{j}}}
{(n-i_{j}+i_{j-1})!\,(1+\ep)_{n-k-i_{j-1}}}
\end{multline*}
et $i_0=0$.
\end{Lemma}

\begin{proof}[D{\'e}monstration]
De nouveau, on suit la d{\'e}marche de la
d{\'e}monstration du Lemme~\ref{lem:Faktor}. Il suffit de montrer
que, pour tous nombres entiers $1\le f_1\le f_2\le\dots\le 
f_H$ tels que
le multi-ensemble $\{f_1,f_2,
\dots,
f_H\}$ soit contenu dans le multi-ensemble
$$
E=\bigcup _{j=1} ^{B}\Big(\left[k+i_{j-1}+1,n+k+i_{j-1}\right]\cup
\left[n-k-i_{j-1}+1,2n-k-i_{j}\right]\Big)
$$
le nombre
\begin{multline} \label{eq:R7A}
\tilde\Phi_n^{-B+1}\textup{d}_n^{H}\,
\frac
{i_{B}!} {i_1!\,(i_2-i_1)!\cdots (i_{B}-i_{B-1})!}\\
\times
\left(
\prod _{j=1} ^{B}\binom {n+k+i_{j-1}}n 
\binom {2n-k-i_{j}}{n-i_{j}+i_{j-1}}\right)
\frac {1} {f_1f_2\cdots f_H}
\end{multline}
est un nombre entier.

Soit $p$ un nombre premier.
En utilisant la notation $\chi(\mathcal A)=1$ si $\mathcal A$ est vrai
et $\chi(\mathcal A)=0$ sinon,
la valuation $p$-adique de
\eqref{eq:R7A} s'{\'e}crit sous la forme
\begin{multline} \label{eq:R7vala}
-(B-1)\cdot\chi\big(p<n\text{ et }\{n/p\}\in[2/3,1[\big)+H\cdot[\log_p(n)]\\
+
\sum _{\ell=1} ^{\infty}\Bigg(
\left[\frac {i_{B}} {p^\ell}\right]
+
\sum _{j=1} ^{B}\bigg(
-\left[\frac  {i_j-i_{j-1}}  {p^\ell}\right]
+\left[\frac {n+k+i_{j-1}} {p^\ell}\right]
-\left[\frac {n} {p^\ell}\right]
-\left[\frac {k+i_{j-1}} {p^\ell}\right]
\\
+\left[\frac {2n-k-i_{j}} {p^\ell}\right]
-\left[\frac {n-i_{j}+i_{j-1}} {p^\ell}\right]
-\left[\frac {n-k-i_{j-1}} {p^\ell}\right]\bigg)\Bigg)
-
\sum _{h=1} ^{H}v_p(f_h)\\
=-(B-1)\cdot\chi\big(p<n\text{ et }\{n/p\}\in[2/3,1[\big)
\kern6cm\\
+
\sum _{\ell=1} ^{[\log_p(n)]}U(\ell)
+\sum _{\ell=[\log_p(n)]+1} ^{\infty}U(\ell)
-
\sum _{h=1} ^{H}(v_p(f_h)-[\log_p(n)]),
\end{multline}
o{\`u} $U(\ell)$ 
d{\'e}signe le sommande de la somme sur $\ell$.

De nouveau, comme le multi-ensemble $\{f_1,f_2,\dots,f_H\}$ est contenu
dans le multi-ensemble $E$, des arguments analogues {\`a} ceux du
Lemme~\ref{lem:Faktor} montrent alors que
\begin{equation} \label{eq:Hilfe}
\sum _{h=1} ^{H}(v_p(f_h)-[\log_p(n)])\le
\sum _{\ell=[\log_p(n)]+1} ^{\infty}
\sum _{j=1} ^{B}\bigg(
\left[\frac {n+k+i_{j-1}} {p^\ell}\right]
+\left[\frac {2n-k-i_{j}} {p^\ell}\right]\bigg).
\end{equation}

Le sommande $U(\ell)$ est {\'e}videmment positif pour tout $\ell$.
De plus,
en vertu de la condition $0\le i_1\le i_2 \le\ldots\le i_B \le n-k$,
on a l'{\'e}galit{\'e}
$$
\sum _{\ell=[\log_p(n)]+1} ^{\infty}U(\ell)
=\sum _{\ell=[\log_p(n)]+1} ^{\infty}
\sum _{j=1} ^{B}\bigg(
\left[\frac {n+k+i_{j-1}} {p^\ell}\right]
+\left[\frac {2n-k-i_{j}} {p^\ell}\right]\bigg).
$$
Si l'on utilise cette {\'e}galit{\'e} et \eqref{eq:Hilfe} dans
\eqref{eq:R7vala}, on voit que l'expression \eqref{eq:R7vala}
%
serait positive si l'on pourrait d{\'e}montrer que
$$
-(B-1)\cdot\chi\big(p<n\text{ et }\{n/p\}\in[2/3,1[\big)+
\sum _{\ell=1} ^{[\log_p(n)]}U(\ell)
\ge0.
$$
Cette in{\'e}galit{\'e} est s{\^u}rement satisfaite si $p\ge n$
ou $\{n/p\}\notin[2/3,1[$ puisque $U(\ell)$ est positif pour tout $\ell$.
Si $p<n$ et $\{n/p\}\in[2/3,1[$, en utilisant encore que $U(\ell)$ est
positif, on voit qu'il suffit de d{\'e}montrer que
$U(1)\ge B-1$. En r{\'e}arrangeant les termes dans la d{\'e}finition
de $U(1)$, on peut {\'e}crire
\begin{multline*}
U(1)=\left[\frac {n+k} {p}\right]
-\left[\frac n {p}\right]
-\left[\frac k {p}\right]\\
+\left[\frac {2n-k-i_B} {p}\right]
+\left[\frac {i_B} {p}\right]
-\left[\frac {n-k} {p}\right]
-\left[\frac {i_B-i_{B-1}} {p}\right]
-\left[\frac {n-i_B+i_{B-1}} {p}\right]\\
+\sum _{j=1} ^{B-1}\bigg(
\left[\frac {n+k+i_j} {p}\right]
-\left[\frac n {p}\right]
-\left[\frac {k+i_j} {p}\right]
+\left[\frac {2n-k-i_j} {p}\right]
-\left[\frac n {p}\right]
-\left[\frac {n-k-i_j} {p}\right]\\
+\left[\frac n {p}\right]
-\left[\frac {i_j-i_{j-1}} {p}\right]
-\left[\frac {n-i_j+i_{j-1}} {p}\right]\bigg)
\end{multline*}
Le lemme sera donc d{\'e}montr{\'e} si l'on r{\'e}ussit de
v{\'e}rifier que
\begin{multline} \label{eq:ungl1}
\left[\frac {n+k} {p}\right]
-\left[\frac n {p}\right]
-\left[\frac k {p}\right]\\
+\left[\frac {2n-k-i_B} {p}\right]
+\left[\frac {i_B} {p}\right]
-\left[\frac {n-k} {p}\right]
-\left[\frac {i_B-i_{B-1}} {p}\right]
-\left[\frac {n-i_B+i_{B-1}} {p}\right]
\ge0
\end{multline}
et que
\begin{multline} \label{eq:ungl2}
\sum _{j=1} ^{B-1}\bigg(
\left[\frac {n+k+i_j} {p}\right]
-\left[\frac n {p}\right]
-\left[\frac {k+i_j} {p}\right]
+\left[\frac {2n-k-i_j} {p}\right]
-\left[\frac n {p}\right]
-\left[\frac {n-k-i_j} {p}\right]\\
+\left[\frac n {p}\right]
-\left[\frac {i_j-i_{j-1}} {p}\right]
-\left[\frac {n-i_j+i_{j-1}} {p}\right]\bigg)
\ge B-1
\end{multline}
pour tout $\ell$. Implicitement, le dernier a d{\'e}j{\`a} \'et\'e
fait. Dans la d{\'e}monstration du Lemme~\ref{lem:arith2}, on a
d{\'e}montr{\'e} que l'expression \eqref{eq:Vp} est $\ge1$. Cela
implique que la partie du sommande dans la premi{\`e}re ligne de
\eqref{eq:ungl2} est $\ge1$ pour $\ell=1$.
Comme la partie du sommande dans la
deuxi{\`e}me ligne est {\'e}videmment positive, l'in{\'e}galit{\'e}
\eqref{eq:ungl2} est v{\'e}rifi{\'e}e.

Pour d{\'e}montrer l'in{\'e}galit{\'e} \eqref{eq:ungl1},
on note $N=\{n/p\}$, $K=\{k/p\}$, $I_1=\{(k+i_{B-1})/p\}$ et
$I_2=\{(k+i_{B})/p\}$
pour les parties fractionnaires de
$n/p$, $k/p$, $(k+i_{B-1})/p$ et $(k+i_{B})/p$.
Avec ces notations,
le membre de gauche de \eqref{eq:ungl1} 
devient
\begin{equation} \label{eq:NKI}
\left[ {N+K} \right]
+\left[ {2N-I_2} \right]
+\left[ {I_2-K} \right]
-\left[ {N-K} \right]
-\left[ {I_2-I_1} \right]
-\left[ {N-I_2+I_1} \right].
\end{equation}
Supposons que l'expression \eqref{eq:NKI} 
soit strictement n{\'e}gative.
Comme on a $N\ge 2/3$, et comme $-\left[ {I_2-I_1} \right]
-\left[ {N-I_2+I_1} \right]\ge0$, cela implique que
$$\left[ {N+K} \right]=\left[ {2N-I_2} \right]
=\left[ {N-K} \right]=0\quad \text{et}\quad
\left[ {I_2-K} \right]=-1.$$
En particulier,
\begin{equation} \label{eq:Bed1}
N+K<1,\quad 
2N-I_2<1,\quad \text{et}\quad 
I_2-K<0.
\end{equation}
En 
vertu de \eqref{eq:Bed1}, on a
$$3N-1<N+I_2<N+K<1,$$
ce qui implique $N<2/3$, une contradiction. L'in{\'e}galit{\'e}
\eqref{eq:ungl1} est donc v{\'e}rifi{\'e}e, ce qui ach{\`e}ve la
d{\'e}monstration du lemme.
\end{proof}


\newbox\ibox
\setbox\ibox\hbox{i)}

\section{D{\'e}monstration du Th{\'e}or{\`e}me \ref{thm:2},
partie \box\ibox}\label{sec:demothm1i}
On se place dans le cadre de la Conjecture~\ref{conj1} du
paragraphe~\ref{sec:conj1}. Rappelons que
\begin{multline*}
{\bf S}_{n,A,B,C,r}\((-1)^A\)
=n!^{A-2Br}
\sum_{k=1}^{\ii}(-1)^{kA}\frac{1}{C!}\frac{\partial^C }{\partial k^C}\bigg(
\(k+\frac n2\)\frac{(k-rn)_{rn}^B(k+n+1)_{rn}^B}{(k)_{n+1}^A}\bigg)
\\={\bf p}_{0,C,n}\((-1)^A\)+
(-1)^C\sum_{l=1}^A\binom{C+l-1}{l-1} {\bf p}_{l,n}
\((-1)^A\)\Li_{C+l}\((-1)^A\).
\end{multline*}
Les ${\bf p}_{l,n}\((-1)^A\)$ ($l\ge 1$) ne
d{\'e}pendent pas de $C$ et leur expression est donn{\'e}e par les deux
{\'e}quations~\eqref{eq:expressionp_mn}
et~\eqref{eq:expressionp_mnplusexplicite}, {\`a} savoir
\begin{multline*}
(-1)^{Brn}{{\bf p}_{l,n}\((-1)^{A}\)}=\frac {(rn)!^{2B}} {n!^{2rB}}
\sum_{j=0}^n
\frac{1}{(A-l)!}\frac {\partial^{A-l}} {\partial\ep^{A-l}}
\(\frac {n} {2}-j+\ep\)\\
\cdot
\(\frac {n!} { (1-\ep)_j \, (1+\ep)_{n-j}}\)^{A}
{\binom {  r n+j-\ep}{rn}}^{B}
{\binom { (r+1) n-j+\ep}{rn}}^{B}\Bigg\vert_{\ep=0}.
\end{multline*}

La partie~i) du Th{\'e}or{\`e}me~\ref{thm:2}
d{\'e}coulent donc de la
Proposition~\ref{cor:3} ci-dessous, o{\`u} la restriction  analytique
$0\le 2Br<A$ 
du paragraphe~\ref{ssec:conjdenom} n'intervient pas.

\begin{Proposition} \label{cor:3}
Soient $A,B,r$ des entiers  tels que $A\ge2$, $B\ge1$ et
$r\ge 0$.
Alors,
pour tout entier $h\ge1$, la
quantit{\'e}
\begin{multline*}
\textup{d}_n^{h-1}\frac {(rn)!^{2B}} {n!^{2rB}}
\sum _{j=0} ^{n}\frac {1} {h!}\frac {\partial^h} {\partial\ep^h}
\Bigg(    \left( \frac{n}{2} - j + \ep \right)
\(\frac {n!} { (1-\ep)_j \, (1+\ep)_{n-j}}\)^{A}\\
\cdot
{\binom {   r n+j-\ep}{rn}}^{B}
{\binom { (r+1) n-j+\ep}{rn}}^{B}\Bigg)\Bigg\vert_{\ep=0}
\end{multline*}
est un nombre entier.
\end{Proposition}

\begin{proof}[D{\'e}monstration]
Notons que, en utilisant la notation des
Corollaires~\ref{cor:1} {\`a} \ref{cor:2a}, il nous faut d{\'e}montrer que
$$\textup{d}_n^{h-1}\frac {(rn)!^{2B}} {n!^{2rB}}
\frac {1} {h!}\frac {\partial^h} {\partial\ep^h}
S_{A,B,r}(n)\Big\vert_{\ep=0}
$$
est un nombre entier.

Nous allons appliquer les
Corollaires~\ref{cor:1} et \ref{cor:2}
pour traiter le cas $A\ge 2$ et $B\ge 2$. Le cas $A\ge 2$
et $B=1$ se
d{\'e}montrent de la m{\^e}me mani\` ere au moyen des
Corollaires~\ref{cor:1a} et \ref{cor:2a}.

Le Corollaire~\ref{cor:1} montre que
$S_{A,B,r}(n)=\ep \cdot (-1)^ns_{A,B,r}(n)$ pour $A\ge2$ pair, tandis que
le Corollaire~\ref{cor:2} donne la m{\^eme} identit{\'e} pour $A\ge3$ impair.
Par la suite, on se concentre sur le cas que $A$ est pair,
car les arguments sont compl{\`e}tement 
analogues pour l'autre cas.
Il s'ensuit que
\begin{equation} \label{eq:3}
\frac {\partial^h} {\partial\ep^h}
S_{A,B,r}(n)\Big\vert_{\ep=0}=(-1)^nh\cdot
\frac {\partial^{h-1}} {\partial\ep^{h-1}} s_{A,B,r}(n)\Big\vert_{\ep=0}.
\end{equation}
Arrangeons les termes qui contiennent $\ep$ dans le sommande de
$s_{A,B,r}(n)$ donn{\'e} au Corollaire~\ref{cor:1}
en termes des briques (voir les Lemmes~\ref{lem:briques} et
\ref{lem:Faktor} au
paragraphe~\ref{sec:LemmesArithm} pour les d{\'e}finitions des briques
$R(\dots)$, $R_1(\dots)$ et $R_2(\dots)$)~: on a
\begin{multline} \label{eq:F8}
\frac {(rn)!^{2B}} {n!^{r(2B)}}s_{A,B,r}(n)\\=
\sum _{0\le i_1\le i_2\le\dots\le i_{A/2+B}\le n} ^{}
{\left( -1 \right) }^{i_{A/2+B-1} + i_{ A/2 + B}}
R(0,n+1;-\ep)\cdot(-\ep)
\bigg(\prod _{q=0} ^{r-1}R(n,0;nq+\ep+1)\bigg)\\
\cdot
R(n-i_{A/2+B},0;rn+i_{A/2+B}+\ep+2)\,
R_1(n,i_{A/2+B},i_{A/2+B};\ep)\,
R_2(n,i_{A/2+B-1},i_{A/2+B};\ep)\\
\cdot
R(0,i_{A/2+B-1}+1;\ep )\cdot\ep \cdot
R(0,i_{A/2+B-1}+1;\ep )\cdot\ep \\
\cdot
R(i_{A/2+B-1},0;1+2\ep)\,
R(i_{A/2+B-1}-i_{A/2+B-2},0;n +  1)\\
\cdot
\Bigg( \prod_{k = B}^{A/2 + B-2}
R( i_k-i_{k-1},0;n  +1)
\binom n {i_k}
R(0,i_k+1;\ep )\cdot \ep
\cdot
R(0,n-i_k+1;-\ep)\cdot(-\ep)\\
\cdot
\binom n {i_k}
R(0,i_k+1;\ep )\cdot\ep
\cdot
R(0,n-i_k+1;-\ep)\cdot(-\ep)\bigg)\\
\cdot
\binom{r n}{i_{B-1}-i_{B-2} }
\binom n {i_{B-1}}
R(0,i_{B-1}+1;\ep)\cdot\ep\cdot
R(0,n-i_{B-1}+1;-\ep)\cdot(-\ep)\cdot
R(n-i_{B-1},0;1-\ep)\\
\cdot
\bigg(\prod _{q=1} ^{r}R(n,0;nq-i_{B-1}-\ep+1)\bigg)
R(0,n-i_{B-1}+1;-\ep)\cdot(-\ep)
\\
\cdot
\left( \prod_{k = 1}^{B-2}\binom{r n}{i_k-i_{k-1}}
\bigg(\prod _{q=1} ^{r}R(n,0;nq-i_k-\ep+1)\bigg)
\bigg(\prod _{q=0} ^{r-1}R(n,0;nq+i_k+\ep+1)\bigg)\right),
\end{multline}
On peut donc r{\'e}{\'e}crire la somme $s_{A,B,r}(n)$ comme
\begin{equation*}
\frac  {(rn)!^{2B}}{n!^{r(2B)}}s_{A,B,r}(n)
=
\sum _{0\le i_1\le i_2\le\dots\le i_{A/2+B}\le n} ^{}
C_1(i_1,\dots,i_{A/2+B})
\prod _{k=1} ^{M}t_k(i_1,\dots,i_{A/2+B}),
\end{equation*}
o{\`u} chaque $C_1(i_1,\dots,i_{A/2+B})$ est un nombre entier et
o{\`u} chaque $t_k$ est une brique {\'e}l{\'e}mentaire $R(\al,\be;\pm \ep+K)$
avec $\al\ge \be$, ou
une brique {\'e}l{\'e}mentaire $R(\al,\be;\pm \ep)$
multipli{\'e}e
par $\ep$ avec $\al< \be $, ou bien encore une des deux
briques sp{\'e}ciales $R_1(n,i_{A/2+B},i_{A/2+B};\ep)$ et
$R_2(n,i_{A/2+B-1},i_{A/2+B};\ep)$.

En vertu de la formule de Leibniz, on en d{\'e}duit que
\begin{align} \notag
\frac {(rn)!^{2B}} {n!^{r(2B)}}&\frac {1} {h!}
\frac {\partial^h} {\partial\ep^h}
S_{A,B,r}(n)\Big\vert_{\ep=0}=
\frac {(-1)^n} {(h-1)!}
\sum _{\ell_1+\dots+\ell_M=h-1} ^{}\frac {(h-1)!}
{\ell_1!\,\ell_2!\cdots\ell_M!}\\
\notag
&\kern2.5cm
\cdot
\sum _{0\le i_1\le i_2\le\dots\le i_{A/2+B}\le n} ^{}C_1(i_1,\dots,i_{A/2+B})
\prod _{k=1} ^{M}\frac {\partial^{\ell_k}} {\partial\ep^{\ell_k}}
t_k(i_1,\dots,i_{A/2+B})\Big\vert_{\ep=0}\\
&=(-1)^n
\underset {\ell_1+\dots+\ell_M=h-1}
{\sum _{0\le i_1\le i_2\le\dots\le i_{A/2+B}\le n}
^{}}C_1(i_1,\dots,i_{A/2+B})
\prod _{k=1} ^{M}\frac {1} {\ell_k!}
\frac {\partial^{\ell_k}} {\partial\ep^{\ell_k}}
t_k(i_1,\dots,i_{A/2+B})\Big\vert_{\ep=0}.
\label{eq:6}
\end{align}
Gr{\^a}ce aux Lemmes~\ref{lem:briques} et \ref{lem:Faktor}, la quantit{\'e}
$$\textup{d}_n^{\ell_k}\frac {1} {\ell_k!}
\frac {\partial^{\ell_k}} {\partial\ep^{\ell_k}}
t_k(i_1,\dots,i_{A/2+B})\Big\vert_{\ep=0}$$
est un nombre entier pour tout $k$. Le membre de droite de \eqref{eq:6}
multipli{\'e} par $\textup{d}_n^{h-1}$  est donc lui aussi entier,
ce qui ach{\`e}ve la d{\'e}monstration.
\end{proof}


\newbox\iibox
\setbox\iibox\hbox{ii)}

\section{D{\'e}monstration du Th{\'e}or{\`e}me \ref{thm:2},
partie \box\iibox} \label{sec:demothm1ii}

Nous commen{\c c}ons par {\'e}crire \eqref{eq:p0} sous la forme
\begin{multline} \label{eq:p0def}
{\bf p}_{0,C,n}(X)\\
=-\sum _{j=1} ^{n}
\sum _{e=1} ^{A}
(-1)^{C+Aj+Brn} 
\binom {C+e-1}{e-1}
\bigg(\frac {1} {(A-e)!}
\frac {\partial^{A-e}}
{\partial \ep^{A-e}}
T_{n,A,B,r}(j;\ep)\Big\vert_{\ep=0}\bigg)
\sum _{k=1} ^{j}\frac {X^{j-k}} {k^{e+C}},
\end{multline}
avec
\begin{multline*}
T_{n,A,B,r}(j;\ep)
=\frac {(rn)!^{2B}} {n!^{2rB}}
\(\frac {n} {2}-j+\ep\)
\(\frac {n!} { (1-\ep)_j \, (1+\ep)_{n-j}}\)^{A}\\
\cdot
{\binom {  r n+j-\ep}{rn}}^{B}
{\binom { (r+1) n-j+\ep}{rn}}^{B}.
\end{multline*}
En 
sp{\'e}cialisant $X=(-1)^A$ et en \'echangeant les sommations, on a donc
\begin{equation} \label{p0alter}
{\bf p}_{0,C,n}\((-1)^A\)=-\sum _{e=1} ^{A} (-1)^{C+Brn}
\binom {C+e-1}{e-1} \sum_{k=1}^n \frac {(-1)^{Ak}}
{k^{e+C}}\,{\bf q}_{k,n,e,A,B,r}\((-1)^A\),
\end{equation}
avec
$$
{\bf q}_{k,n,e,A,B,r}\((-1)^A\)=\sum_{j=k}^n \left(
\frac {1} {(A-e)!}
\frac {\partial^{A-e}}
{\partial \ep^{A-e}}T_{n,A,B,r}(j;\ep)\Big\vert_{\ep=0}\right).
$$

Compte-tenu de l'expression \eqref{p0alter} pour ${\bf p}_{0,C,n}\((-1)^A\)$,
la partie ii) du Th{\'e}or{\`e}me~\ref{thm:2} d{\'e}coule
imm{\'e}diatement de la proposition
suivante. Sa d{\'e}monstration est bas{\'e}e sur les Corollaires~\ref{cor:10}
et \ref{cor:11} du
paragraphe~\ref{CorollairesA} et,
en outre des Lemmes~\ref{lem:briques} et \ref{lem:Faktor},
le Lemme~\ref{lem:briques1} du paragraphe~\ref{sec:LemmesArithm}.

\begin{Proposition} \label{conj:4} Fixons les entiers
 $A\ge 2$, $B\ge 0$, $r\ge 0$ et $n\ge 0$. Pour
  tout entier $k\in\{1, \ldots, n\}$ et tout entier
$e\in\{1, \ldots,  A\}$, le nombre
$$
2\dd_n^{A-e}{\bf q}_{k,n,e,A,B,r}\((-1)^A\) .
$$
est un nombre entier qui est divisible par $k$.
\end{Proposition}
\begin{proof}
Soit tout d'abord $A$ pair, $A\ge2$.
On remarque que
$$
\sum _{j=k} ^{n}T_{n,A,B,r}(j;\ep)$$
est {\'e}gal au membre de gauche de \eqref{eq:dernier} multipli{\'e}
par $(rn)!^{2B}/n!^{2rB}$. En {\'e}crivant le membre de droite de
\eqref{eq:dernier} en termes de briques {\'e}l{\'e}mentaires et
sp{\'e}ciales (voir les Lemmes~\ref{lem:briques} et \ref{lem:Faktor}
au paragraphe~\ref{sec:LemmesArithm}), on voit
que ${\bf q}_{k,n,e,A,B,r}\(1\)$ est
{\'e}gal {\`a}
\begin{multline} \label{eq:briques}
-\frac {1} {(A-e)!}
\frac {\partial^{A-e}}
{\partial \ep^{A-e}}\Bigg\{
\frac {k-\ep} {2}
\sum _{0\le i_1\le\dots\le i_{A/2+B}\le n-k} ^{}
\frac
{(-1)^{i_{A/2+B-1}}\,(i_{A/2+B})!} {i_1!\,(i_2-i_1)!\cdots (i_{A/2+B}-i_{A/2+B-1})!}\\
\cdot
\(\prod _{j=1} ^{B}
\bigg(\prod _{q=0} ^{r-1}R(n,0;k+i_{j-1}+qn+1-\ep)\bigg)
R_1(n,n-k-i_j,n-i_j+i_{j-1};\ep)\)\\
\cdot(-1)^{k+i_B}\cdot\ep\cdot R(-k-i_{B},n-k-i_{B}+1;\ep)\\
\cdot
\Bigg(\prod _{j=B+1} ^{A/2+B-1}
\ep\cdot R(-k-i_{j},n-k-i_{j}+1;\ep)\kern7cm
\\
\cdot
\binom {n+i_j-i_{j-1}}{k+i_j}\cdot(-\ep)\cdot R(0,k+i_j+1;-\ep)
\cdot\ep\cdot R(0,n-k-i_{j-1}+1;\ep)\Bigg)\\
\cdot\frac {n!\,(\ep)_{i_{A/2+B}-i_{A/2+B-1}}\,(1-2\ep)_{k+i_{A/2+B-1}}
\,(1-\ep)_{n-i_{A/2+B}-1}}
{(1-\ep)_n\,(1-2\ep)_{k-1}\,(1-\ep)_{k+i_{A/2+B}}\,
(1+\ep)_{n-k-i_{A/2+B}}}\Bigg\}\Bigg\vert_{\ep=0}
\end{multline}
Pour simplifier, notons $R_7(n,k,i_{A/2+B},i_{A/2+B-1};\ep)$ la fraction
dans la derni{\`e}re ligne de \eqref{eq:briques}.

On s'interesse {\`a} la quantit{\'e}
$2\dd_n^{A-e}{\bf q}_{k,n,e,A,B,r}\(1\)/k$, c'est-{\`a}-dire, il faut
multiplier l'expression \eqref{eq:briques} par $2\dd_n^{A-e}/k$~:
de fa{\c c}on similaire {\`a} la d{\'e}monstration de la
Proposition~\ref{cor:3},
on {\'e}crit cette nouvelle expression sous la forme
\begin{multline*}
-\frac {2\dd_n^{A-e}} {k}\frac {1} {(A-e)!}
\frac {\partial^{A-e}}
{\partial \ep^{A-e}}\Bigg\{
\frac {k-\ep} {2}\sum _{0\le i_1\le\dots\le i_{A/2+B}\le n-k} ^{}
C_2(i_1,\dots,i_{A/2+B})\\
\cdot
R_7(n,k,i_{A/2+B},i_{A/2+B-1};\ep)
\prod _{h=1} ^{M}t_h(i_1,\dots,i_{A/2+B})
\Bigg\}\Bigg\vert_{\ep=0},
\end{multline*}
o{\`u} chaque $C_2(i_1,\dots,i_{A/2+B})$ est un nombre entier et
o{\`u} chaque $t_h$ est une brique {\'e}l{\'e}mentaire $R(\al,\be;\pm \ep+K)$
avec $\al\ge \be$, ou
une brique {\'e}l{\'e}mentaire $R(\al,\be;\pm \ep)$
multipli{\'e}e
par $\ep$ avec $\al< \be $, ou bien encore une
brique sp{\'e}ciale $R_1(n,n-k-i_j,n-i_j+i_{j-1};\ep)$.

En vertu de la formule de Leibniz, cette derni{\`e}re expression peut
s'{\'e}crire comme
{\refstepcounter{equation}\label{eq:final1}}
\alphaeqn
\begin{align} \notag
-\dd_n^{A-e}&\Bigg\{
\sum _{\ell_0+\dots+\ell_M=A-e} ^{}\frac {1}
{\ell_0!\,\ell_1!\cdots\ell_M!}
\sum _{0\le i_1\le\dots\le i_{A/2+B}\le n-k} ^{}
C_2(i_1,\dots,i_{A/2+B})\\
\label{eq:final1a}
&\kern2cm
\cdot
\frac {\partial^{\ell_0}} {\partial\ep^{\ell_0}}
R_7(n,k,i_{A/2+B},i_{A/2+B-1};\ep)
\prod _{h=1} ^{M}\frac {\partial^{\ell_h}} {\partial\ep^{\ell_h}}
t_h(i_1,\dots,i_{A/2+B})
\Bigg\}\Bigg\vert_{\ep=0}\\
\notag
&+\frac {\dd_n} {k}\dd_n^{A-e-1}\Bigg\{
\sum _{\ell_0+\dots+\ell_M=A-e-1} ^{}\frac {1}
{\ell_0!\,\ell_1!\cdots\ell_M!}
\sum _{0\le i_1\le\dots\le i_{A/2+B}\le n-k} ^{}
C_2(i_1,\dots,i_{A/2+B})\\
&\kern2cm
\cdot
\frac {\partial^{\ell_0}} {\partial\ep^{\ell_0}}
R_7(n,k,i_{A/2+B},i_{A/2+B-1};\ep)
\prod _{h=1} ^{M}\frac {\partial^{\ell_h}} {\partial\ep^{\ell_h}}
t_h(i_1,\dots,i_{A/2+B})
\Bigg\}\Bigg\vert_{\ep=0}.
\label{eq:final1b}
\end{align}
\reseteqn

On doit distinguer deux cas.
Si $i_{A/2+B}=i_{A/2+B-1}$, alors $R_7(n,k,i_{A/2+B},i_{A/2+B-1};\ep)$ se
d{\'e}compose en briques {\'e}l{\'e}mentaires,
\begin{align*}
R_7(n,k,&i_{A/2+B},i_{A/2+B-1};\ep)
=R_7(n,k,i_{A/2+B},i_{A/2+B};\ep)\\
&=
\frac {n!\,(1-2\ep)_{k+i_{A/2+B}}
\,(1-\ep)_{n-i_{A/2+B}-1}}
{(1-\ep)_n\,(1-2\ep)_{k-1}\,(1-\ep)_{k+i_{A/2+B}}\,
(1+\ep)_{n-k-i_{A/2+B}}}\\
&=(-\ep)\cdot R(0,n+1;-\ep)\cdot R(k+i_{A/2+B},0;1-2\ep)\\
&\kern1cm
\cdot(-\ep)\cdot R(0,k+i_{A/2+B}+1;-\ep)
\cdot R(n-i_{A/2+B}-1,0;1-\ep)\\
&\kern1cm
\cdot (-1)^{k-1}\cdot \ep\cdot R(-k+1,n-k-i_{A/2+B}+1;\ep)
\cdot R(k-1,0;1-\ep)\\
&\kern1cm
\cdot (-2\ep)\cdot R(0,k;-2\ep).
\end{align*}
Comme dans la d{\'e}monstration de la Proposition~\ref{cor:3},
on utilise les Lemmes~\ref{lem:briques} et \ref{lem:Faktor} et le fait
que $k$ divise $\dd_n$ pour en
d{\'e}duire que \eqref{eq:final1a} et \eqref{eq:final1b} restreintes
{\`a} $i_{A/2+B}=i_{A/2+B-1}$ sont des entiers.

Lorsque $i_{A/2+B}>i_{A/2+B-1}$, on remarque que
$$R_7(n,k,i_{A/2+B},i_{A/2+B-1};\ep)=\ep\cdot
R_3(n,k,i_{A/2+B},i_{A/2+B-1};\ep),$$
o{\`u} $R_3(\dots)$ est 
la brique sp{\'e}ciale d{\'e}finie au Lemme~\ref{lem:briques1}
du paragraphe~\ref{sec:LemmesArithm}.
Par con\-s{\'e}\-quent, pour $\ell_0\ge1$ on a
\begin{multline*}
\frac {1} {\ell_0!}\frac {\partial^{\ell_0}} {\partial\ep^{\ell_0}}
R_7(n,k,i_{A/2+B},i_{A/2+B-1};\ep)\Bigg\vert_{\ep=0}\\
=\frac {1} {(\ell_0-1)!}\frac {\partial^{\ell_0-1}} {\partial\ep^{\ell_0-1}}
R_3(n,k,i_{A/2+B},i_{A/2+B-1};\ep)\Bigg\vert_{\ep=0}.
\end{multline*}
Le Lemme~\ref{lem:briques1} avec
$m_1=i_{A/2+B}$ et $m_2=i_{A/2+B-1}$ joint aux
Lemmes~\ref{lem:briques} et \ref{lem:Faktor} nous permet alors de
conclure que \eqref{eq:final1a} et \eqref{eq:final1b} restreintes
{\`a} $i_{A/2+B}>i_{A/2+B-1}$ sont aussi des entiers.

La d{\'e}monstration de la proposition dans le cas que $A$ est impair
est compl{\`e}tement analogue. La seule diff{\'e}rence 
est qu'il faut utiliser le Corollaire~\ref{cor:11}, au lieu 
du Corollaire~\ref{cor:10}.
\end{proof}

\begin{Remark} Pour le lecteur avide d'identit\'es hyperg\'eom\'etriques, 
indiquons
les identit\'es suivantes, qui correspondent \`a certains des cas pour $A=e$ de
la Proposition~\ref{conj:4} ci-dessus~:
\begin{align*}
{\bf q}_{k,n,0,0,0,r}(1)&=-\sum_{j=k}^n \left(\frac n2 -j\right)= \frac k2 (n+k-1), \\
{\bf q}_{k,n,1,1,0,r}(-1)&=(-1)^{k+1}\sum_{j=k}^n \left(\frac n2 -j\right)\binom{n}{j}= 
(-1)^k\frac k2 \binom{n}{k},\\
{\bf q}_{k,n,2,2,0,r}(1)&=-\sum_{j=k}^n \left(\frac n2
  -j\right)\binom{n}{j}^2=
\frac k2 \binom{n}{k}\binom{n-1}{k-1} 
\end{align*}
et
$$
{\bf q}_{k,n,0,0,1,1}(1)=-\sum_{j=k}^n \left(\frac n2
  -j\right)\binom{n+j}{n}\binom{2n-j}{n}=
\frac k2 \binom{n+k}{n}\binom{2n-k+1}{n+1}.
$$
La premi\`ere est \'evidemment tr\`es facile \`a montrer et nous avons obtenu
les trois autres  avec {\it Maple}, qui produit imm\'ediatement les membres de droite lorsqu'on lui demande
de calculer formellement les sommes sur $j$ (cela d\'ecoule en effet
de formules de sommation \eqref{4F3[-1]} et \eqref{eq:5F4} 
pour des s{\'e}ries hyperg\'eom\'etriques tr{\`e}s bien {\'e}quilibr{\'e}es, 
connues de {\em Maple}). Nous avons r\'eussi \`a
g\'en\'eraliser la pr\'esence de ce facteur $k$, qui \'etait une voie vers
la d\'emonstration compl\`ete de la partie ii) du Th\'eor\`eme~\ref{thm:2}, ce que nous n'avions
pas pu enti\`erement faire dans~\cite{KR1}~: l'identit\'e d'Andrews s'est av\'er\'ee
\^etre un tr\`es  bon guide sur ce chemin.
\end{Remark}


\section{D{\'e}monstration du Th{\'e}or{\`e}me
\ref{thm:3}, partie {\rm i)}, et des Th{\'e}or{\`e}mes
\ref{prop7gene} et \ref{prop:Phi}}\label{sec:demothmvaszud}

Bien que le Th{\'e}or{\`e}me~\ref{thm:3}, partie~i)
soit un cas particulier du Th{\'e}or{\`e}me~\ref{prop7gene},
nous allons en donner une d{\'e}monstration ind{\'e}pendante, qui donnera l'id{\'e}e du cas g{\'e}n{\'e}ral.
Pour cela, rappelons que la s{\'e}rie
$$
{\bf S}_{n,4,2,1,1}(1)=\sum_{k=1}^{\ii}\frac{\partial}{\partial k} \(
\(k+\frac n2\)\frac{(k-n)_{n}^{2}\,(k+n+1)_{n}^2}{(k)_{n+1}^4}\)
={\bf u}_n\zeta(4)-{\bf v}_n
$$
permet de constuire une suite d'approximations rationnelles
de $\zeta(4)$, o{\`u}
$$
{\bf u}_n=\sum_{j=0}^n \frac{\dd}{\dd j}\(\frac n2-j\)
\binom{n}{j}^4\binom{n+j}{n}^2\binom{2n-j}{n}^2.
$$
La pr{\'e}sence implicite des nombres harmoniques ne permet pas de voir
imm{\'e}diatement  que
ces nombres sont des entiers, ce qui est beaucoup plus visible
gr{\^a}ce {\`a}
l'identit{\'e} suivante, qui d{\'e}coule de la
Proposition~\ref{cor:A4} avec $A=4$ et $B=2$~:
$$
{\bf u}_n=(-1)^{n+1}\sum_{0\le i\le j\le n} \binom{n}{j}^2\binom{n}{i}^2
\binom{n+j}{n}\binom{n+j-i}{n}\binom{2n-i}{n}.
$$

Il nous faut maintenant prouver que $\Phi_n^{-1}{\bf u}_n$ est encore entier, o{\`u}
$
\dis\Phi_n=\underset{\{n/p\}\in[2/3,1[}{\prod_{p\text{ premier}}} p.
$
\begin{Lemma} \label{lem:arith1}
Il existe des entiers $C_{j,n}$ $(0\le  j\le n)$ tels que
$$
{\bf u}_n=\sum_{j=0}^n
\binom{n+j}{n}\binom{2n-j}{n}C_{j,n}.
$$
\end{Lemma}
\begin{proof}[D{\'e}monstration] Notons que~:
\begin{equation*}
\binom{n}{j}\binom{n+j-i}{n}\binom{2n-i}{n}
=\binom{2n-j}{n}
\binom{2n-i}{j-i}\binom{n+j-i}{j}.
\end{equation*}
Il suffit donc de poser
$
\displaystyle C_{j,n}=(-1)^{n+1}\binom{n}{j}\sum_{i=0}^j
\binom{n}{i}^2\binom{2n-i}{j-i}\binom{n+j-i}{j}.
$
\end{proof}

Quand on applique le Lemme~\ref{lem:arith2} aux facteurs
$\binom{n+j}{n}\binom{2n-j}{n}$ pour $j\in\{0, \ldots, n\}$,
mis en {\'e}vidence par le lemme pr{\'e}c{\'e}dent,
on arrive bien au r{\'e}sultat attendu~:
$\Phi_n$ divise ${\bf u}_n$,
ce qui est exactement l'{\'e}nonc{\'e} de la partie~i) du
Th{\'e}or{\`e}me~\ref{thm:3}.

En utilisant les expressions binomiales donn{\'e}es par
la Proposition~\ref{cor:A4} au paragraphe~\ref{sec:hyperharm},
on obtient plus g{\'e}n{\'e}ralement
le Th{\'e}or{\`e}me~\ref{prop7gene}.


\begin{proof}[D{\'e}monstration du Th{\'e}or{\`e}me~\ref{prop7gene}]
On peut supposer que
$B\ge 2$, sinon il n'y a d{\'e}j{\`a} plus rien {\`a} montrer.
On s'int{\'e}resse au nombre
$$
{\bf p}_{A-1,n}\((-1)^A\)
=\sum _{j=0} ^{n}\frac {\dd} {\dd j}
 \left( \frac{n}{2} - j\right)
\binom{n}{j}^{A}
{\binom {n+j}{n}}^{B}
{\binom {2n-j}{n}}^{B},
$$
qui est aussi {\'e}gal {\`a} la
quantit{\'e} $P_n(A,B)$ dans la Proposition~\ref{cor:A4}.
Cette derni{\`e}re quantit{\'e} est
elle-m{\^e}me {\'e}gale {\`a} $p_n(A,B)$, dont la d{\'e}finition
est diff{\'e}rente selon que $A$ est pair ou impair. Si $A$ est impair,
par un jeu d'{\'e}criture similaire au Lemme~\ref{lem:arith1},
on peut reformuler l'expression pour $p_n(A,B)$, donn{\'e}e {\`a}
la Proposition~\ref{cor:A4}~:
\begin{multline*}
p_n(A,B)=
\sum _{0\le i_1\le i_2\le\dots\le i_{m+B-1}\le n} ^{}
\kern-9pt(-1)^{i_{m+B-1}}
\binom {n+i_{m+B-1}-i_{m+B-2}}{i_{m+B-1}}
\\
\cdot
{\binom {n+i_{m+B-1}}n}\binom {2n-i_{m+B-1}}n
\\
\cdot
\Bigg(\prod _{k=B} ^{m+B-2}{\binom n{i_k}}\binom{2n-i_k}{2n-i_{k+1}}\binom
{n+i_{k}-i_{k-1}}{i_k}\Bigg)
\binom n{i_{B-1}}\binom {2n-i_{B-1}}{2n-i_B}\binom n{i_{B-1}-i_{B-2}}\\
\cdot
\Bigg(\prod _{k=1} ^{B-2}{\binom {n+i_k}{n}}\binom {2n-i_k}n
\binom n{i_{k}-i_{k-1}}\Bigg),
\end{multline*}
avec une reformulation similaire de $p_n(A,B)$
lorsque $A$ est pair.
Or, dans le sommande, on remarque
la pr{\'e}sence du produit
$$
{\binom {n+i_{m+B-1}}n}\binom {2n-i_{m+B-1}}n\prod_{k = 1}^{B-2}
\binom {n + i_k }{n}
\binom {2n - i_k}{n}.
$$
On applique donc le
Lemme~\ref{lem:arith2} du paragraphe~\ref{sec:LemmesArithm} pour conclure.
\end{proof}

Pour d{\'e}montrer le Th{\'e}or{\`e}me~\ref{prop:Phi}, en plus
des briques $R(\dots)$, $R_1(\dots)$ et $R_2(\dots)$ d{\'e}j{\`a}
utilis{\'e}es,
nous avons aussi besoin des briques sp{\'e}ciales $R_4(\dots)$
et $R_5(\dots)$ qui apparaissent
aux Lemmes~\ref{lem:binom1} et \ref{lem:binom2}.

\begin{proof}[Esquisse de la d{\'e}monstration du
Th{\'e}or{\`e}me~\ref{prop:Phi}]
Il s'agit essentiellement de raffiner celle du
Th{\'e}o\-r{\`e}\-me~\ref{thm:2}, donn{\'e}e
aux paragraphes~\ref{sec:demothm1i} et
\ref{sec:demothm1ii}.
Nous consid{\'e}rons tout d'abord le cas que $A$ est pair.
On peut l{\`a} aussi supposer que
$B\ge 2$, sinon il n'y a d{\'e}j{\`a} plus rien {\`a} montrer.
Rappelons que, dans ce cas, la
cl{\'e} des d{\'e}mon\-stra\-tions aux paragraphes~\ref{sec:demothm1i} et
\ref{sec:demothm1ii} {\'e}tait la d{\'e}composition en briques \eqref{eq:F8}
de la somme multiple $s_{A,B,r}(n)$ du Corollaire~\ref{cor:1}
et la d{\'e}composition en briques \eqref{eq:briques}
de la somme multiple pour ${\bf q}_{k,n,e,A,B,r}\(1\)$
du Corollaire~\ref{cor:10}.

Consid{\'e}rons tout d'abord la d{\'e}composition \eqref{eq:F8}.
Dans notre cas on a $r=1$.
On remarque la pr{\'e}sence du produit
$$
\prod_{k = 1}^{B-2}
     \binom{ n +
i_k+\ep}{ n}\binom{2n  -
       i_k-\ep}{ n}
$$
dans le sommande de $s_{A,B,1}(n)$. Dans \eqref{eq:F8} (avec $r=1$),
pour $k=1,2,\dots,B-2$, on remplace alors
le produit des briques
\begin{equation} \label{eq:E3}
R(n,0;i_k+\ep+1)R(n,0;n-i_k-\ep+1)= \binom{ n +
i_k+\ep}{ n}\binom{2n  -
       i_k-\ep}{ n}
\end{equation}
par la brique sp{\'e}ciale $R_4(n,i_k;\ep)$ d{\'e}finie au
Lemme~\ref{lem:binom1} au paragraphe~\ref{sec:LemmesArithm}.
De m{\^e}me, on remplace le produit
\begin{multline} \label{eq:R4}
R(n-i_{A/2+B},0;rn+i_{A/2+B}+\ep+2)\,
R_1(n,i_{A/2+B},i_{A/2+B};\ep)\,
R_2(n,i_{A/2+B-1},i_{A/2+B};\ep)\\
\times
R(i_{A/2+B-1}-i_{A/2+B-2},0;n  + 1)\\
\times
\Bigg( \prod_{k = B}^{A/2 + B-2}
R( i_k-i_{k-1},0;n  +1)
\binom n {i_k}
R(0,i_k+1;\ep )\cdot\ep
\cdot
R(0,n-i_k+1;-\ep)\cdot(-\ep)\bigg)\\
\kern0cm
\times
R(n-i_{B-1},0;1-\ep)\,R(n,0;n-i_{B-1}-\ep+1)\,
R(0,n-i_{B-1}+1;-\ep)\cdot(-\ep)
\end{multline}
dans \eqref{eq:F8} (avec $r=1$) par
\begin{equation} \label{eq:RR4}
R(n,0;1-\ep)\,R_5(n,i_{B-1},i_{B},\dots,i_{A/2+B};\ep),
\end{equation}
o{\`u} la brique sp{\'e}ciale
$R_5(n,i_{B-1},i_{B},\dots,i_{A/2+B};\ep)$ est d{\'e}finie au
Lemme~\ref{lem:binom2}. 
Puis, on r{\'e}p{\`e}te les m{\^e}mes arguments que ceux du
paragraphe~\ref{sec:demothm1i}, avec
$R_4(n,i_k;\ep)$ au lieu de \eqref{eq:E3} et
\eqref{eq:RR4} au lieu de \eqref{eq:R4},
et en utilisant bien s{\^u}r le fait que les
nombres \eqref{eq:tildePhi} et \eqref{eq:G5} sont des nombres
entiers. La conclusion est
finalement l'{\'e}nonc{\'e} du th{\'e}or{\`e}me pour
${\bf p}_{l,n}\((-1)^A\)$, $l\ge1$, dans le cas que $A\ge 2$ pair.

D'autre part, dans la d{\'e}composition \eqref{eq:briques} (avec
$r=1$) on remplace
\begin{multline} \label{eq:R6A}
\frac {(i_{A/2+B})!} {i_1!\,(i_2-i_1)!\cdots
(i_{A/2+B}-i_{A/2+B-1})!}\\
\times
\(\prod _{j=1} ^{B}
R(n,0;k+i_{j-1}+1-\ep)
R_1(n,n-k-i_j,n-i_j+i_{j-1};\ep)\)
\end{multline}
par
\begin{equation} \label{eq:R6B}
\frac {(i_{A/2+B})!} {i_B!\,(i_{B+1}-i_{B})!\cdots
(i_{A/2+B}-i_{A/2+B-1})!}R_6(n,i_{1},i_{2},\dots,i_{B};\ep),
\end{equation}
o{\`u} $R_6(n,i_{1},i_{2},\dots,i_{B};\ep)$ est la brique
sp{\'e}ciale du Lemme~\ref{lem:binom3}.
Puis, on r{\'e}p{\`e}te les m{\^e}mes arguments que ceux du
paragraphe~\ref{sec:demothm1ii}, avec
\eqref{eq:R6B} au lieu de \eqref{eq:R6A},
et en utilisant bien s{\^u}r le fait que le
nombre \eqref{eq:G10} est un nombre
entier. La conclusion est
finalement l'{\'e}nonc{\'e} du th{\'e}or{\`e}me pour
${\bf p}_{0,C,n}\((-1)^A\)$ dans le cas que $A\ge 2$ pair.

Lorsque $A\ge 3$ est impair,
il faut utiliser les Corollaires~\ref{cor:11} et
\ref{cor:2}, ainsi qu'une variante
du Lemme~\ref{lem:binom2} adapt{\'e}e au Corollaire~\ref{cor:2}.
La d{\'e}monstration de cette
variante {\'e}tant calqu{\'e}e sur celle du Lemme~\ref{lem:binom2},
nous ne rentrons
pas dans les d{\'e}tails afin de ne pas alourdir davantage le texte.
\end{proof}

\section{D{\'e}monstration du Th{\'e}or{\`e}me
\ref{thm:3}, partie {\rm ii)}, et du Th{\'e}or{\`e}me
\ref{thm:C=3}}
\label{sec:Phi}

Nous commen\c{c}ons avec la partie~ii) du Th{\'e}or{\`e}me
\ref{thm:3}.
Le cas sp{\'e}cial $r=1$, $A=4$, $B=2$, $C=1$ du
Th{\'e}or{\`e}me~\ref{prop:Phi} montre d{\'e}j{\`a} que
$2\tilde\Phi_n^{-1}\textup{d}_n^{4}{\bf v}_n$
est un nombre entier. Pour achever la d{\'e}monstration de la
partie~ii) du Th{\'e}or{\`e}me~\ref{thm:3}, il nous reste {\`a}
d{\'e}montrer que tout nombre premier $p$, $n<p\le \frac {3} {2}n$,
divise $2\textup{d}_n^{4}{\bf v}_n$. Comme un tel nombre premier $p$
est premier avec $2\textup{d}_n$, il suffit de d{\'e}montrer que la
valuation $p$-adique $v_p({\bf v}_n)$ est au moins 1.
Pour cela, nous aurons besoin du lemme hyperg{\'e}om{\'e}trique suivant.

\begin{Lemma} \label{lem:10}
Pour tout nombre entier $q\ge 0$, on a
\begin{multline} \label{eq:idzeta4}
\sum _{j=0} ^{q-1}{\binom {q+j-1}j}^4{\binom {q-1}j}^2\frac
{(q+j-1)!^2} {(2q+j-1)!^2}\Bigg(\(\frac {q}
{2}+j\)\(2H_j^{(3)}-2H_{q+j-1}^{(3)}\)\\
\kern1cm
+\(H_j^{(2)}-H_{q+j-1}^{(2)}\)\Big(-1+\(\frac {q}
{2}+j\)(6H_j-6H_{q+j-1}+2H_{2q+j-1}-2H_{q+j-1})\Big)\Bigg)=0,
\end{multline}
o{\`u} les nombres $H_j$ sont les nombres harmoniques et les nombres
$H_j^{(e)}$ sont les nombres harmoniques g{\'e}n{\'e}ralis{\'e}s
d{\'e}finis par
$H_j^{(e)}= \sum_{i=1} ^j 1/i^e.$
\end{Lemma}

\begin{proof}On v\'erifie facilement que la somme {\`a} gauche
de \eqref{eq:idzeta4} s'exprime sous la forme
\begin{multline} \label{eq:der1}
\frac {1} {2}
\frac {\partial} {\partial\ep}\frac {\partial^2} {\partial\et^2}
\Bigg(\sum _{j=0} ^{n}
\left( \frac{q}{2}+ j -\ep    \right)
\binom {q+j-1}j\\
\cdot
\frac {({ \textstyle 1 -
   2 \ep}) _{q+j-1}}
 {({ \textstyle 1 - 2 \ep}) _{j}\,(q-1)!}
\frac {({ \textstyle 1 - \ep -
   \et}) _{q+j-1} } {({ \textstyle 1 - \ep - \et})
   _{j}\,(q-1)!}
\frac {({ \textstyle 1 - \ep +
   \et}) _{q+j-1}} {({ \textstyle 1 - \ep + \et}) _{j}\,(q-1)!}
\\
\kern1cm\cdot
\frac {(q-1)!} {({ \textstyle 1 -
   \ep}) _{j}\,\left( q-j-1 \right) !}
\frac {(q-1)!} {({ \textstyle 1 -
   \ep}) _{j}\,({ \textstyle 1 + 2 \ep}) _{q-j-1} }
\frac{{({ \textstyle 1 - \ep})
   _{q+j-1} ^2} }{\left(
    2 q+ j-1 \right) !\,({
   \textstyle 1 - 2 \ep}) _{2 q+ j -1}
   } \Bigg)\Bigg\vert_{\ep=\et=0},
\end{multline}
soit en notation hyperg{\'e}om{\'e}trique~:
\begin{multline*}
\frac {1} {2}
\frac {\partial} {\partial\ep}\frac {\partial^2} {\partial\et^2}
\Bigg( \frac{
     {({ \textstyle 1 - \ep}) _{q-1} }^2\,
     ({ \textstyle 1 - \ep - \et}) _{q-1} \,
     ({ \textstyle 1 - \ep + \et}) _{q-1} }{4
     \left( q-1 \right) !\,\left( 2 q -1 \right) !\,
     ({ \textstyle 1 + 2 \ep}) _{q-1} \,
     ({ \textstyle 1 - 2 \ep + q}) _{q-1} }\\
\times
{} _{9} F _{8} \!\left [ \begin{matrix} {  q-2 \ep, 1 + \frac{q}{2}-
      \ep ,  1 - 2 \ep - q,q -\ep -
      \et , q-\ep + \et, q,  q-\ep,
      q-\ep,1 - q}\\ {
       \frac{q}{2}-\ep ,2 q,1 -
      \ep + \et, 1 - \ep -
      \et  , 1 -
      2 \ep, 1 - \ep, 1 - \ep,  2 q-2 \ep}\end{matrix} ;
{\displaystyle 1}\right ]  \Bigg)\Bigg\vert_{\ep=\et=0}.
\end{multline*}
Notons $f(\ep,\et)$ le terme entre parenth{\`e}ses.
Nous allons d{\'e}montrer que $f(\ep,\et)$ est une fonction paire de $\ep$,
c'est-{\`a}-dire, que $f(\ep,\et)=f(-\ep,\et)$. Il est clair que cela
implique l'{\'e}nonc{\'e}.

La d{\'e}monstration de cette sym{\'e}trie est bas{\'e}e sur la
transformation de Bailey entre deux s{\'e}ries $_9F_8$ tr{\`e}s bien
{\'e}quilibr{\'e}es (voir \cite[(2.4.4.1)]{SlatAC})~:
\begin{multline} \label{eq:9F8}
 {} _{9} F _{8} \!\left [ \begin{matrix} { a, 1 + \frac{a}{2}, b, c, d, e,
    f,}\\ { \frac{a}{2}, 1 + a - b, 1 + a
    - c, 1 + a - d, 1 + a - e, 1 + a - f,}\end{matrix} \right . \\
\left.\begin{matrix} 2 + 3 a - b - c - d - e - f + N, -N\\
-1 - 2 a + b + c + d + e + f - N, 1
    + a + N\end{matrix} ; {\displaystyle 1}\right ]
\\
 =
   \frac{ ({ \textstyle 1 + a}) _{N} \,
      ({ \textstyle 2 + 2 a - b - c - d - e}) _{N} \,
      ({ \textstyle 2 + 2 a - b - c - d - f}) _{N} \,
      ({ \textstyle 1 + a - e - f}) _{N} }{({ \textstyle 2 + 2 a - b - c -
       d}) _{N} \,({ \textstyle 1 + a - e}) _{N} \,
      ({ \textstyle 1 + a - f}) _{N} \,
      ({ \textstyle 2 + 2 a - b - c - d - e - f}) _{N} }\\
\times
{} _{9} F _{8} \!\left [ \begin{matrix} { 1 + 2 a - b - c - d,
       \frac{3}{2} + a - \frac{b}{2} - \frac{c}{2} - \frac{d}{2}, 1 + a - c -
       d, 1 + a - b - d,}\\ { \frac{1}{2} + a - \frac{b}{2} - \frac{c}{2} - \frac{d}{2}, 1
       + a - b, 1 + a - c,1 + a - d, 2 + 2 a - b - c - d - e,}
     \end{matrix}\right .\\
\left. \begin{matrix} 1 + a - b - c, e, f, 2 + 3 a - b - c - d - e - f +
       N, -N\\
  2 + 2 a - b -
       c - d - f, -a + e + f - N, 2 + 2 a - b - c - d + N
\end{matrix} ;
       {\displaystyle 1}\right ] ,
\end{multline}
o{\`u} $N$ est un entier positif. On applique cette transformation
{\`a} la s{\'e}rie $_9F_8$ dans la d{\'e}finition de $f(\ep,\et)$~: on obtient
l'expression {\'e}quivalente
\begin{multline*}
 \frac{
     {({ \textstyle 1 + \ep}) _{q-1} }^2\,
     ({ \textstyle 1 - \ep - \et}) _{q-1} \,
     ({ \textstyle 1 - \ep + \et}) _{q-1} }{4
     \left( 2 q-1 \right) !\,
     ({ \textstyle 1 - 2 \ep}) _{q-1} \,
     ({ \textstyle 1 + 2 \ep}) _{q-1} \,
     ({ \textstyle q+1}) _{q-1} }\\
\times
{} _{9} F _{8} \!\left [ \begin{matrix} { q, 1 + \frac{q}{2}, 1 - q,
      q-\ep , q-\ep , q+\ep -
      \et , q+\ep + \et , q, 1 - q}\\ {
      \frac{q}{2}, 2 q, 1 + \ep, 1 + \ep, 1 -
      \ep + \et, 1 - \ep - \et, 1,
      2 q}\end{matrix} ; {\displaystyle 1}\right ]
\end{multline*}
pour $f(\ep,\et)$. Si l'on applique la transformation de Bailey une
deuxi{\`e}me fois, on arrive exactement {\`a} $f(-\ep,\et)$.
\end{proof}

\begin{proof}[D{\'e}monstration de la partie {\em ii)} du Th{\'e}or{\`e}me
\ref{thm:3}]
Soit $p$ un nombre premier, $n<p\le \frac {3} {2}n$. Si $n=2$, on peut
v{\'e}rifier l'{\'e}nonc{\'e} directement. Nous pouvons
donc d{\'e}sormais supposer que $p>3$.

En utilisant la formule \eqref{eq:p0def} pour $r=1$, $A=4$, $B=2$ et
$C=1$, on obtient que $\mathbf v_n$ est donn{\'e} par l'expression
\begin{multline} \label{eq:vn}
\mathbf v_n=\frac {1} {2}\sum _{j=0} ^{n}
\sum _{e=1} ^{4}
\frac {e} {(4-e)!}
\frac {\partial^{4-e}}
{\partial \ep^{4-e}}\Bigg(
\(\frac {n} {2}-j+\ep\)
\(\frac {n!} { (1-\ep)_j \, (1+\ep)_{n-j}}\)^{4}\\
\cdot
{\binom {   n+j-\ep}{n}}^{2}
{\binom { 2 n-j+\ep}{n}}^{2}\Bigg)\Bigg\vert_{\ep=0}
\cdot
\(H^{(e+1)}_j+(-1)^{e+1}H^{(e+1)}_{n-j}\).
\end{multline}
Nous voulons d{\'e}montrer que la valuation $p$-adique de cette
expression est au moins 1. Pour faciliter le calcul, on fait deux
remarques~:

\medskip
\begin{itemize}
\item[(R1)] On peut {\`a} volont{\'e} multiplier ou
diviser l'expression \eqref{eq:vn} par
des nombres qui sont premiers avec $p$, et la valuation $p$-adique de
l'expression obtenue sera au moins 1 si et seulement si
c'est vrai pour l'expression originelle.
\item[(R2)] On peut
{\`a} volont{\'e} remplacer un nombre $x$ par un nombre $y$
dans \eqref{eq:vn} si $v_p(x-y)\ge1$,
et la valuation $p$-adique de l'expression obtenue sera au moins 1 si
et seulement si c'est vrai pour l'expression originelle.
Cette remarque vaut puisque, apr\`es un instant de r\'eflexion,
on se convainc qu'aucun d{\'e}nominateur dans
\eqref{eq:vn} n'est divisible par $p$.
\end{itemize}

\medskip
Nous rappelons maintenant que le Lemme~\ref{lem:arith2}
du paragraphe~\ref{sec:LemmesArithm} montre que $p$
divise le produit $\binom {n+j}n\binom {2n-j}n$ pour tout $j$, $0\le
j\le n$. C'est en fait compl{\`e}tement {\'e}vident pour un nombre
premier $p$ entre $n+1$ et $\frac {3} {2}n$ puisque
$v_p\big((2n-j)!)=1$ pour $j\le n/2$ et $v_p\big((n+j)!)=1$ pour $j\ge
n/2$.
Il s'ensuit que
\begin{equation} \label{eq:bin1}
v_p\(\frac {\partial}
{\partial \ep}\Bigg({\binom {   n+j-\ep}{n}}^{2}
{\binom { 2 n-j+\ep}{n}}^{2}
\Bigg)\Bigg\vert_{\ep=0}\)
\ge1,
\end{equation}
\begin{equation} \label{eq:bin2}
\frac {1} {2}\frac {\partial^{2}}
{\partial \ep^{2}}\Bigg({\binom {   n+j-\ep}{n}}^{2}
{\binom { 2 n-j+\ep}{n}}^{2}
\Bigg)\Bigg\vert_{\ep=0}=\frac {1} {p^2}{\binom {   n+j}{n}}^{2}
{\binom { 2 n-j}{n}}^{2}+N_1,
\end{equation}
o{\`u} $v_p(N_1)\ge1$, et
\begin{multline} \label{eq:bin3}
\frac {1} {6}\frac {\partial^{3}}
{\partial \ep^{3}}\Bigg({\binom {   n+j-\ep}{n}}^{2}
{\binom { 2 n-j+\ep}{n}}^{2}
\Bigg)\Bigg\vert_{\ep=0}\\
=\frac {1} {p^2}{\binom {   n+j}{n}}^{2}
{\binom { 2 n-j}{n}}^{2}
\bigg(2H_j-2H_{n-j}+2H_{2n-j}\kern4cm\\
+\frac {2\cdot (\chi(j>n/2)-\chi(j<n/2))} {p}
-2H_{n+j}\bigg)
+N_2,
\end{multline}
o{\`u} $v_p(N_2)\ge1$, et, comme
auparavant, $\chi(\mathcal A)=1$ si $\mathcal A$ est vrai
et $\chi(\mathcal A)=0$ sinon.
En appliquant le Lemme~\ref{lem:arith2}
et 
le fait~\eqref{eq:bin1} dans \eqref{eq:vn},
en combinaison avec la remarque~(R2), on obtient que
les valuations $p$-adiques des sommandes pour $e=4$ et $e=3$ sont
au moins 1. Ensuite on applique la formule de Leibniz pour calculer
les d{\'e}riv{\'e}es des sommandes pour $e=2$ et $e=1$.
En utilisant de nouveau le Lemme~\ref{lem:arith2}, puis les 
faits~\eqref{eq:bin1}--\eqref{eq:bin3} en combinaison avec la
remarque~(R2), on conclut que la valuation $p$-adique de $\mathbf v_n$
est au moins 1 si et seulement si la valuation $p$-adique de
\begin{multline*}
\sum _{j=0} ^{n}
2\(\frac {n} {2}-j\)
{\binom nj}^4
\frac {1} {p^2}{\binom {   n+j}{n}}^{2}
{\binom { 2 n-j}{n}}^{2}
\cdot
\(H^{(3)}_j-H^{(3)}_{n-j}\)\\
+
\sum _{j=0} ^{n}
\(\frac {n} {2}-j\)
{\binom nj}^4
\frac {1} {p^2}
{\binom {   n+j}{n}}^{2}
{\binom { 2 n-j}{n}}^{2}
\cdot
\(H^{(2)}_j+H^{(2)}_{n-j}\)\kern4cm\\
\cdot
\(6H_j-6H_{n-j}+2H_{2n-j}+\frac {2\cdot (\chi(j>n/2)-\chi(j<n/2))} {p}
-2H_{n+j}+\frac {1} {\frac {n} {2}-j}\)
\end{multline*}
est au moins 1. En utilisant le Lemme~\ref{lem:arith2} une
derni{\`e}re fois, on remarque que, dans ces deux sommes, les
sommandes pour $p-n\le j\le 2n-p$ sont divisibles par $p$. De plus,
les sommes partielles $
\sum _{j=0} ^{p-n-1}$ et $
\sum _{j=2n-p} ^{n}$ sont {\'e}gales,
ce que l'on peut v{\'e}rifier en faisant la substitution $j\to n-j$. On
en conclut que $v_p(\mathbf v_n)\ge1$ si et seulement si la valuation
$p$-adique de
\begin{multline} \label{eq:sumq}
\sum _{j=0} ^{p-n-1}
2\(\frac {n} {2}-j\)
{\binom nj}^4
\frac {1} {p^2}{\binom {   n+j}{n}}^{2}
{\binom { 2 n-j}{n}}^{2}
\cdot
\(H^{(3)}_j-H^{(3)}_{n-j}\)\\
+
\sum _{j=0} ^{p-n-1}
\(\frac {n} {2}-j\)
{\binom nj}^4
\frac {1} {p^2}
{\binom {   n+j}{n}}^{2}
{\binom { 2 n-j}{n}}^{2}
\cdot
\(H^{(2)}_j+H^{(2)}_{n-j}\)\\
\cdot
\(6H_j-6H_{n-j}+2H_{2n-j}-\frac {2} {p}
-2H_{n+j}+\frac {1} {\frac {n} {2}-j}\)
\end{multline}
est au moins 1.

{\'E}crivons $p=n+q$. On v{\'e}rifie facilement que
\begin{align} \label{eq:congr1}
\binom nj&\equiv (-1)^j\binom {q+j-1}j\pmod p,\\
\binom {n+j}n=\binom {n+j}j&\equiv (-1)^j\binom {q-1}j\pmod p,
\label{eq:congr2}
\end{align}
et
\begin{align} \notag
\frac {n!} {p}\binom {2n-j}n&=\frac {(n-j+1)(n-j+2)\cdots(2n-j)} {p}\\
\label{eq:congr3}
&\equiv
(-1)^{q+1}(q+j-1)!\,(2q+j-1)!^{-1}\pmod p.
\end{align}
Puis, en groupant par deux les termes dans la d{\'e}finition du nombre
harmonique $H_{p-1}$, on remarque que $v_p(H_{p-1})\ge
1$. Par la m{\^e}me raison on a $v_p\(H_{p-1}^{(3)}\)\ge1$. Nous
pr{\'e}tendons que l'on a aussi $v_p\(H_{p-1}^{(2)}\)\ge1$ si $p>3$.
En effet, on a
$$
\sum _{i=1} ^{p-1}i^2=\frac{\left(p-1 \right) p\left( 2p-1 \right)
}{6},$$
ce qui implique que la somme des restes quadratiques modulo $p$
est divisible par $p$ si $p>3$. Comme l'ensemble
des r{\'e}ciproques des restes quadratiques modulo $p$ est le m{\^e}me que
l'ensemble des restes quadratiques modulo $p$, l'{\'e}nonc{\'e}
$v_p\(H_{p-1}^{(2)}\)\ge1$ en
d{\'e}coule\footnote{Les \'enonc\'es sur $H_{p-1}$ et $H_{p-1}^{(2)}$ sont
bien s\^ur des corollaires du th\'eor\`eme de Wolstenholme
(voir~\cite[Chapter VII]{hw}), qui affirme que
$v_p(H_{p-1})\ge 2$ si $p>3$ et $v_p\(H_{p-1}^{(2)}\)\ge1$ si $p\ge 3$.}.
Pour $0\le j\le q=p-n$, on a donc
\begin{align} \notag
H_{n-j}&=H_{p-1}-\(\frac {1} {p-q-j+1}+\frac {1} {p-q-j+2}+\dots+\frac
{1} {p-1}\)\\
&=H_{q+j-1}+N_3,\\
\notag
H_{n+j}&=H_{p-1}-\(\frac {1} {p-q+j+1}+\frac {1} {p-q+j+2}+\dots+\frac
{1} {p-1}\)\\
&=H_{q-j-1}+N_4,\\
\notag
H_{2n-j}-\frac {1} {p}&=
H_{p-1}+\(\frac {1} {p+1}+\frac {1} {p+2}+\dots+\frac
{1} {2p-2q-j}\)\\
\notag
&=H_{p-1}+\(\frac {1} {p+1}+\frac {1} {p+2}+\dots+\frac
{1} {2p-1}\)\\
\notag
&\kern3cm-
\(\frac {1} {2p-2q-j+1}+\frac {1} {2p-2q-j+2}+\dots+\frac
{1} {2p-1}\)\\
&=H_{2q+j-1}+N_5,\\
\notag
H_{n-j}^{(2)}&=H_{p-1}^{(2)}-\(\frac {1} {(p-q-j+1)^2}+
\frac {1} {(p-q-j+2)^2}+\dots+\frac
{1} {(p-1)^2}\)\\
&=-H_{q+j-1}^{(2)}+N_6,\\
\notag
H_{n-j}^{(3)}&=H_{p-1}^{(3)}-\(\frac {1} {(p-q-j+1)^3}+
\frac {1} {(p-q-j+2)^3}+\dots+\frac
{1} {(p-1)^3}\)\\
&=H_{q+j-1}^{(3)}+N_7,
\label{eq:congr4}
\end{align}
o{\`u} $N_3,N_4,N_5,N_6,N_7$ sont des nombres dont la valuation
$p$-adique est au moins 1.
En substituant \eqref{eq:congr1}--\eqref{eq:congr4} dans
\eqref{eq:sumq}, on obtient que $v_p(\mathbf v_n)\ge1$ si et seulement
si la valuation $p$-adique de
\begin{multline*}
\sum _{j=0} ^{q-1}
2\(-\frac {q} {2}-j\)
{\binom {q+j-1}j}^4
{\binom {   q-1}{j}}^{2}
\frac {(q+j-1)!^2} {(2q+j-1)!^2}
\cdot
\(H^{(3)}_j-H^{(3)}_{q+j-1}\)\\
+
\sum _{j=0} ^{q-1}
\(-\frac {q} {2}-j\)
{\binom {q+j-1}j}^4
{\binom {   q-1}{j}}^{2}
\frac {(q+j-1)!^2} {(2q+j-1)!^2}
\cdot
\(H^{(2)}_j-H^{(2)}_{q+j-1}\)\\
\cdot
\(6H_j-6H_{q+j-1}+2H_{2q+j-1}
-2H_{q-j-1}-\frac {1} {\frac {q} {2}+j}\)
\end{multline*}
est au moins 1. Or cette derni{\`e}re expression est l'oppos{\'e}e du
membre de gauche dans \eqref{eq:idzeta4}, qui, selon le
Lemme~\ref{lem:10}, vaut 0. La valuation $p$-adique de $\mathbf v_n$ est
donc au moins 1.
\end{proof}

Pour la d{\'e}monstration du Th{\'e}or{\`e}me
\ref{thm:C=3}, l'approche est similaire.
Le cas sp{\'e}cial $r=1$, $A=4$, $B=2$, $C=3$ du
Th{\'e}or{\`e}me~\ref{prop:Phi} montre d{\'e}j{\`a} que
$2\tilde\Phi_n^{-1}\textup{d}_n^{6}{\bf p}_{0,3,n}\(1\)$
est un nombre entier. Pour achever la d{\'e}monstration
du Th{\'e}or{\`e}me~\ref{thm:C=3}, il nous reste {\`a}
d{\'e}montrer que pour tout nombre premier $p$, $n<p\le \frac {3} {2}n$,
la valuation $p$-adique de ${\bf p}_{0,3,n}\(1\)$ est au moins 1.
Pour cela, nous avons besoin d'un autre lemme hyperg{\'e}om{\'e}trique.

\begin{Lemma} \label{lem:11}
Pour tout nombre entier $q\ge 0$, on a
\begin{multline} \label{eq:C=3}
\sum _{j=0} ^{q-1}{\binom {q+j-1}j}^4{\binom {q-1}j}^2\frac
{(q+j-1)!^2} {(2q+j-1)!^2}\Bigg(\(\frac {q}
{2}+j\)\(4H_j^{(5)}-4H_{q+j-1}^{(5)}\)\\
\kern1cm
+\(H_j^{(4)}-H_{q+j-1}^{(4)}\)\Big(-1+\(\frac {q}
{2}+j\)(6H_j-6H_{q+j-1}+2H_{2q+j-1}-2H_{q+j-1})\Big)\Bigg)=0.
\end{multline}
\end{Lemma}

\begin{proof}[Esquisse de la d{\'e}monstration]
Soit $f(\ep,\et_1,\et_2)$ la somme donn{\'e}e par
\begin{multline*}
\sum _{j=0} ^{n}
\left( \frac{q}{2}+ j -\ep    \right)
\binom {q+j-1}j\\
\cdot
\frac {({ \textstyle 1 -
   2 \ep}) _{q+j-1}}
 {({ \textstyle 1 - 2 \ep}) _{j}\,(q-1)!}
\frac {({ \textstyle 1 - \ep -
   \et_1}) _{q+j-1} } {({ \textstyle 1 - \ep - \et_1})
   _{j}\,(q-1)!}
\frac {({ \textstyle 1 - \ep +
   \et_1}) _{q+j-1}} {({ \textstyle 1 - \ep + \et_1}) _{j}\,(q-1)!}
\frac {(q-1)!} {({ \textstyle 1 -
   \ep-\et_2}) _{j}\,\left( q-j-1 \right) !}\\
\cdot
\frac {(q-1)!} {({ \textstyle 1 -
   \ep+\et_2}) _{j}\,({ \textstyle 1 + 2 \ep}) _{q-j-1} }
\frac{{({ \textstyle 1 - \ep-\et_2})
   _{q+j-1} } }{\left(
    2 q+ j-1 \right) !}
\frac{{({ \textstyle 1 - \ep+\et_2})
   _{q+j-1} } }{({
   \textstyle 1 - 2 \ep}) _{2 q+ j -1}
   } .
\end{multline*}
On peut v{\'e}rifier que la somme {\`a} gauche de
\eqref{eq:C=3} est {\'e}gale {\`a}
\begin{equation} \label{eq:Abl}
\frac {1} {12}
\frac {\partial} {\partial\ep}\frac {\partial^4} {\partial\et_1^4}
f(\ep,\et_1,\et_2)\Big\vert_{\ep=\et_1=0}-
\frac {1} {4}
\frac {\partial} {\partial\ep}\frac {\partial^2} {\partial\et_1^2}
\frac {\partial^2} {\partial\et_2^2}
f(\ep,\et_1,\et_2)\Big\vert_{\ep=\et_1=\et_2=0}.
\end{equation}
{\'E}videmment, la somme not{\'e}e $f(\ep,\et)$ dans \eqref{eq:der1} n'est
rien d'autre que $f(\ep,\et,0)$. En effet, $f(\ep,\et_1,\et_2)$ peut
s'aussi exprimer comme une s{\'e}rie $_9F_8$ tr{\`e}s bien
{\'e}quilibr{\'e}e. On d{\'e}montre alors
de la m{\^e}me fa{\c c}on que lors de la
d{\'e}monstration du Lemme~\ref{lem:10} que
$f(\ep,\et_1,\et_2)=f(-\ep,\et_1,\et_2)$. De nouveau, on a besoin
d'une double application de la transformation \eqref{eq:9F8} de Bailey
pour le faire. Les termes dans
\eqref{eq:Abl} sont donc nuls, ce qui implique l'{\'e}nonc{\'e}.
\end{proof}

\begin{proof}[Esquisse de la d{\'e}monstration du Th{\'e}or{\`e}me
\ref{thm:C=3}]
On suit la d{\'e}marche de la d{\'e}monstration de la partie~ii)
du Th{\'e}or{\`e}me~\ref{thm:3} ci-dessus. C'est-{\`a}-dire, soit
$p$ un nombre premier, $n<p\le \frac {3} {2}n$. On part
de l'expression \eqref{eq:p0def} (avec $r=1$, $A=4$, $B=2$, $C=3$)
pour ${\bf p}_{0,3,n}\(1\)$, et ensuite on r{\'e}duit l'expression
obtenue en utilisant les remarques 
(R1) et (R2) 
dans la
d{\'e}monstration pr{\'e}c{\'e}dente. Le r{\'e}sultat de la
r{\'e}duction sera le membre de gauche de \eqref{eq:C=3},
au signe pr{\`e}s. On conclut
comme auparavant que la valuation
$p$-adique de ${\bf p}_{0,3,n}\(1\)$ est au moins 1, ce qui 
ach{\`e}ve la d{\'e}monstration du th{\'e}or{\`e}me.
\end{proof}

\begin{Remark}
Le Th{\'e}or{\`e}me~\ref{thm:3}, partie~ii), et
le Th{\'e}or{\`e}me~\ref{thm:C=3} sont extr{\^e}mement sp{\'e}ciaux.
En effet, quand on regarde l'expression \eqref{eq:p0},
il est \'evident que l'on ne peut esp{\'e}rer
une telle am{\'e}loriation du
Th{\'e}or{\`e}me~\ref{prop:Phi} que lorsqu'il n'y a pas trop
de d{\'e}riv{\'e}es, c'est-{\`a}-dire lorsque $A$ est petit
(dans ces deux th{\'e}or{\`e}mes, on a $A=4$). Mais cela ne suffit
pas. Il faut aussi les deux identit{\'e}s
hyperg{\'e}om{\'e}trico-harmoniques aux Lemmes~\ref{lem:10} et
\ref{lem:11}, qui d{\'e}pendent fortement de la possibilit{\'e} d'une
{\'e}criture compacte de la somme en question comme une somme de
d{\'e}riv{\'e}es
(voir \eqref{eq:der1} et \eqref{eq:Abl}), ainsi que de la
co{\"\i}ncidence {\og miraculeuse \fg}
que la fonction {\`a} laquelle on applique les
d{\'e}riv{\'e}es soit paire, ce que l'on a d{\'e}montr{\'e} en utilisant
la transformation \eqref{eq:9F8} de Bailey. Or, si $C$ est pair, une
{\'e}criture comme somme de d{\'e}riv{\'e}es ne semble pas
d'{\^e}tre possible. En particulier,
c'est la raison pour laquelle
cette approche ne marche pas pour $C=2$, et, comme on l'a d{\'e}j{\`a}
mentionn{\'e} au paragraphe~\ref{resultats}, on peut v{\'e}rifier
num\'eriquement qu'une telle am{\'e}lioration  du
Th{\'e}or{\`e}me~\ref{prop:Phi} n'a pas lieu.

D'autre part, pour
$C\ge4$, il serait n{\'e}cessaire d'introduire davantage de
param{\`e}tres auxiliaires (comme $\et$ dans la d{\'e}monstration du
Lemme~\ref{lem:10} et $\et_1$ et $\et_2$ dans la d{\'e}monstration du
Lemme~\ref{lem:11}). Or, il n'y a plus la place de le faire et, l\`a aussi, on
ne sait pas comment exprimer la somme comme somme de
d{\'e}riv{\'e}es. De nouveau, les exp\'eriences indiquent
qu'une telle am{\'e}lioration du
Th{\'e}or{\`e}me~\ref{prop:Phi} n'a pas lieu pour $C\ge4$.
\end{Remark}

\section{Encore un peu d'hyperg{\'e}om{\'e}trie}
\label{sec:hyp}

Nous d{\'e}montrons dans ce paragraphe que l'on peut d{\'e}montrer
l'{\'e}quivalence\footnote{Par {\og {\'e}quivalence\fg} nous entendons
{\it deux} choses~: premi{\`e}rement, que les valeurs de deux s{\'e}ries
sont {\'e}gales et, deuxi{\`e}mement, une v{\'e}rification directe
(qui n'utilise surtout pas l'irrationalit{\'e} de $\zeta(3)$)
que les d{\'e}compositions en combinaisons
lin{\'e}aires de 1 et $\zeta(3)$, que l'on obtient en appliquant la
d{\'e}marche classique, sont identiques.}
de la s{\'e}rie de Beukers, Gutnik et Nesterenko
\eqref{eq:bgn} et de la s{\'e}rie de Ball \eqref{eq:Ball} directement,
{\`a} l'aide des identit{\'e}s classiques pour les s{\'e}ries
hyperg{\'e}om{\'e}triques, et de m{\^e}me
l'{\'e}quivalence des s{\'e}ries
\eqref{eq:Apery2} et \eqref{eq:Ball2} pour $\zeta(2)$. Cela constitue
une alternative attractive
{\`a} la d{\'e}monstration de Zudilin \cite{zu1} qui utilise l'algorithme
de Gosper--Zeilberger.\footnote{Dans \cite[dernier paragraphe]{zu1},
Zudilin donne aussi une autre d{\'e}monstration {\og classique\fg}
de l'{\'e}galit{\'e}
des valeurs de ces deux s{\'e}ries (mais pas de l'{\'e}galit{\'e} des
d{\'e}compositions en combinaisons lin{\'e}aires de 1 et $\zeta(3)$)
en utilisant une identit{\'e}
int{\'e}grale due {\`a} Bailey. Cette d{\'e}monstration est essentiellement
{\'e}quivalente {\`a} celle donn{\'e}e ici, une fois constat{\'e},
apr{\`e}s application du th{\'e}or{\`e}me des
r{\'e}sidus, que
l'int{\'e}grale dans l'identit{\'e} de Bailey peut s'{\'e}crire comme
somme de deux s{\'e}ries $_4F_3$.}

Nous commen\c cons avec la s{\'e}rie de Ball
\begin{multline*}
{\bf B}_n=n!^2\sum_{k=1}^{\ii}\(k+\frac n2\)
\frac{(k-n)_n(k+n+1)_n}{(k)_{n+1}^4}
\\=\frac{n!^7(3n+2)!}{2(2n+1)!^5}\;  _{7}F_6\left[
\begin{array}{cccccc}
3n+2, \frac32 n+2 , n+1    ,n+1,n+1,n+1,  n+1 \\
     \frac32 n+1 ,  2n+2,2n+2,2n+2,2n+2, 2n+2
\end{array}\,;\, 1\right],
\end{multline*}
que l'on {\'e}crit comme la limite suivante~:
\begin{multline*}
{\bf B}_n=\lim_{\ep\to0}\frac{n!^7(3n+2)!}{2(2n+1)!^5}\\
\times  {}_{7}F_6\left[
\begin{matrix}
3n+2\ep+2, \frac32 n+\ep+2 , n+\ep+1,n+\ep+1,n+\ep+1,n+\ep+1,  n+\ep+1 \\
     \frac32 n+\ep+1 ,  2n+\ep+2,2n+\ep+2,2n+\ep+2,2n+\ep+2, 2n+\ep+2
\end{matrix}\,;\, 1\right].
\end{multline*}
{\`A} cette s{\'e}rie, on applique alors la version
non termin{\'e}e de la transformation de Whipple, due {\`a} Bailey
(voir \cite[(2.4.4.3)]{SlatAC})~:
\begin{multline} \label{eq:7643}
 {} _{7} F _{6} \!\left [ \begin{matrix} { a, 1 + \frac{a}{2}, b, c, d, e, f}\\ {
    \frac{a}{2}, 1 + a - b, 1 + a - c, 1 + a - d, 1 + a - e, 1 + a -
    f}\end{matrix} ; {\displaystyle 1}\right ] \\
\kern-4cm =
   \frac{\Gamma({ \textstyle 1 + a - d}) \,
       \Gamma({ \textstyle 1 + a - e}) \,\Gamma({ \textstyle 1 + a - f}) \,
       \Gamma({ \textstyle 1 + a - d - e - f}) }{\Gamma({ \textstyle 1 + a})
        \,\Gamma({ \textstyle 1 + a - d - e}) \,
       \Gamma({ \textstyle 1 + a - d - f}) \,
       \Gamma({ \textstyle 1 + a - e - f}) } \\
\times
{} _{4} F _{3} \!\left [ \begin{matrix} { 1 + a - b - c, d, e, f}\\ { 1 +
        a - b, 1 + a - c, -a + d + e + f}\end{matrix} ; {\displaystyle 1}\right
        ] \\
\kern-2cm+
    \frac{\Gamma({ \textstyle 1 + a - b}) \,
       \Gamma({ \textstyle 1 + a - c}) \,\Gamma({ \textstyle 1 + a - d}) \,
       \Gamma({ \textstyle 1 + a - e}) \,\Gamma({ \textstyle 1 + a - f})
      }{\Gamma({ \textstyle 1 + a})
        \,\Gamma({ \textstyle 1 + a - b - c}) \,\Gamma({ \textstyle d}) \,
       \Gamma({ \textstyle e})
       }\\
\times
\frac {  \Gamma({ \textstyle 2 + 2\,a - b - c - d - e - f}) \,
       \Gamma({ \textstyle -1 - a + d + e + f})}
 { \Gamma({ \textstyle 2 + 2\,a - b - d - e - f}) \,
       \Gamma({ \textstyle 2 + 2\,a - c - d - e - f}) \,
       \Gamma({ \textstyle f})}\\
\times
{} _{4} F _{3} \!\left [ \begin{matrix} { 1 + a - d - e, 1 + a - d - f, 1
        + a - e - f, 2 + 2\,a - b - c - d - e - f}\\ { 2 + 2\,a - b - d - e -
        f, 2 + 2\,a - c - d - e - f, 2 + a - d - e - f}\end{matrix} ;
        {\displaystyle 1}\right ] .
\end{multline}
De la sorte, on obtient l'expression
$$
{\bf B}_n=\lim_{\ep\to0}
\frac {1} {2\ep}\Bigg(
\sum _{k=1} ^{\infty}\frac {(k-n)_n\,(k-n-\ep)_n} {(k)_{n+1}^2}
-\sum _{k=1} ^{\infty}\frac {(k-n)_n\,(k-n+\ep)_n} {(k+\ep)_{n+1}^2}
\Bigg),
$$
soit, en utilisant le Th{\'e}or{\`e}me de l'H{\^o}pital,
$$
{\bf B}_n=-
\frac {1} {2}\sum _{k=1} ^{\infty}\frac {\partial} {\partial\ep}
\frac {(k-n+\ep)_n^2} {(k+\ep)_{n+1}^2}\Bigg\vert_{\ep=0},
$$
ce qui est exactement la moiti{\'e} de la s{\'e}rie de Beukers, Gutnik et
Nesterenko \eqref{eq:bgn}. (Pour un $q$-analogue de ces s{\'e}ries et
de cette d{\'e}monstration, voir \cite[Par.~4.2]{KrRZAA}.)

L'{\'e}galit{\'e} du coefficient $a_n$ d'Ap{\'e}ry (voir \eqref{eq:a_n})
et du coefficient ${\bf a}_n$ de Ball (voir \eqref{eq:bolda_n}) vient
directement en sp{\'e}cialisant la transformation de Whipple
\eqref{Whipple} en $a=n+1$, $b= - n+\ep$, $c= n + \ep + 1$,
$d=  -n$, $e=1$ et $f= 1 + 2 \ep$, et ensuite en faisant tendre $\ep$
vers 0. De la sorte, on obtient $a_n$ {\`a} gauche et, encore en vertu de
la formule de l'H{\^o}pital, le coefficient $\mathbf a_n$ {\`a} droite.
Il d{\'e}coule de ces deux {\'e}galit{\'e}s
que les coefficients $b_n/2$ et $\mathbf
b_n$ sont {\'e}gaux.

\medskip
Consid{\'e}rons maintenant les s{\'e}ries \eqref{eq:Apery2} et
\eqref{eq:Ball2}. En termes hyperg{\'e}om{\'e}triques, la s{\'e}rie
\eqref{eq:Apery2} s'{\'e}crit sous la forme
$${\left( -1 \right) }^n  \frac{
      {n!}^4}{{\left
         ( 2 n +1\right) !}^2}\,
 {} _{3} F _{2} \!\left [ \begin{matrix} { n+1, n+1, n+1}\\ {
      2 n+2, 2n+2}\end{matrix} ; {\displaystyle 1}\right ],$$
alors que la s{\'e}rie \eqref{eq:Ball2} s'{\'e}crit sous la forme
$$- {\left( -1 \right) }^n  \frac{
        {n!}^5\,\left( 3 n +2\right) !  }{2\,
     {\left( 2 n +1\right) !}^4}\,
 {} _{6} F _{5} \!\left [ \begin{matrix} { 3 n+2,
        \frac{3 n}{2}+2, n+1, n+1, n+1, n+1}\\ {  \frac{3 n}{2}+1,
        2 n+2, 2n+2, 2n+2, 2n+2}\end{matrix} ; {\displaystyle
        -1}\right ].$$
En sp{\'e}cialisant $a=3 n+2$, $b=c=d=e=n+1$ dans \eqref{eq:6F5}
(une transformation qui est en fait un cas limite 
de la transformation
\eqref{eq:7643}, ce que l'on voit en restreignant $f$ aux valeurs
enti{\`e}res n{\'e}gatives, et en faisant ensuite tendre $f$ vers $-\infty$),
on constate que l'oppos{\'e}e de \eqref{eq:Ball2} est la moiti{\'e} de
\eqref{eq:Apery2}.

Comme on l'a d{\'e}j{\`a} remarqu{\'e},
l'{\'e}galit{\'e} des coefficients $\al_n$ et $p_n$
est le contenu
de la Proposition~\ref{cor:A2} sp{\'e}cialis{\'e}e en $A=3$.
Il s'ensuit
que les coefficients $-\be_n/2$ et $q_n$ sont {\'e}gaux.

\section{Perspectives}
\label{sec:persp}

Dans ce paragraphe de conclusion, nous pr{\'e}sentons quelques
probl{\`e}mes ouverts, ou en cours de r{\'e}solution,
que l'on peut essayer d'aborder par les m{\'e}thodes introduites
dans cet article.

\subsection{Les s{\'e}ries asym{\'e}triques de Zudilin}\label{ssec:zudconj}

Pour d{\'e}montrer l'irrationalit{\'e} d'au moins un
des nombres $\zeta(5), \zeta(7), \zeta(9), \zeta(11)$, Zudilin \cite{zu3}
a  utilis{\'e} la s{\'e}rie d{\'e}riv{\'e}e suivante, qui est une
perturbation des s{\'e}ries du type \eqref{eq:seriesderiveesformeslineaires}
(voir aussi \cite{fi}, {\`a} qui nous empruntons la formulation ci-dessous)~:
$$
{\bf Z}_n = \prod_{u=1}^{10} \frac{((13+2u)n)!} {(27n)!^6}
\sum_{k=1} ^{\infty} \frac{1}{2}
\frac{\partial^2}{\partial k^2}
\left( \left(k+\frac{37n}{2}\right)
\frac{(k-27n)_{27n} ^3 (k+37n+1)_{27n} ^3 }{\prod_{u=1} ^{10}
(k+(12-u)n)_{(13+2u)n+1}} \right).
$$
Il existe des rationnels ${\bf z}_{j,n}$ ($j=0, 1, \dots, 4$) tels que
${\bf Z}_n = {\bf z}_{0,n} + \sum_{j=1} ^{4}
{\bf z}_{j,n} \zeta(2j+3)$ et
$
2\dd_{35n} ^3 \dd_{34n} \dd_{33n} ^8 {\bf z}_{j,n}\in\,\mathbb{Z}.
$
Num{\'e}ri\-que\-ment, il semble que l'on
ait ici aussi une conjecture des d\'enominateurs~:
\begin{equation}
\label{zuddenom2}
2\dd_{35n} ^3 \dd_{34n} \dd_{33n}^7
{\bf z}_{j,n}\in\,\mathbb{Z},
\end{equation}
c'est-\`a-dire que l'on peut esp\'erer gagner un facteur $\dd_{33n}.$ 
{\em Stricto sensu}, Zudilin ne prouve pas cette conjecture  mais la contourne en  montrant 
qu'il existe un entier $\hat\Phi_n$ tel que\footnote{Le facteur $\hat\Phi_n$ co{\"\i}ncide
essentiellement avec notre $\tilde\Phi_n$ quand on sp{\'e}cialise
la construction de Zudilin au cas des s{\'e}ries {\og sym{\'e}triques\fg}
consid{\'e}r{\'e}es dans cet
article, ce qui donne encore
plus de poids aux observations num{\'e}riques faites au
paragraphe~\ref{resultats} entre les Th{\'e}or{\`e}mes~\ref{prop:Phi}
et \ref{thm:C=3}.} 
\begin{equation}
\label{zuddenom}
2\dd_{35n} ^3 \dd_{34n} \dd_{33n} ^8 \hat\Phi_n^{-1} {\bf z}_{j,n}\in\,\mathbb{Z}.
\end{equation}
Puisque $\dd_{33n}<{\hat{\Phi}}_n$, l'estimation~\eqref{zuddenom}  est  
meilleure que celle pr\'edite par la seule conjecture des 
d\'enominateurs~\eqref{zuddenom2}, 
et elle suffit \`a montrer le th\'eor\`eme envisag\'e.

N\'eanmoins, en oubliant ce facteur $\hat\Phi_n$ qui joue un r{\^o}le {\`a} part, Zudilin {\'e}nonce une
conjecture tr{\`e}s
g{\'e}n{\'e}rale sur le d{\'e}nominateur commun aux
coefficients des combinaisons lin{\'e}aires de valeurs de z{\^e}ta
que l'on peut
construire {\`a} l'aide de ce type de s{\'e}ries
{\og asym{\'e}triques \fg}~: voir \cite[paragraphe~9]{zu3}
pour l'{\'e}nonc{\'e} de cette conjecture.
Comme les coefficients des s{\'e}ries asym{\'e}triques consid{\'e}r{\'e}es par Zudilin
s'expriment encore comme des s{\'e}ries hyperg{\'e}om{\'e}tri\-ques
tr{\`e}s bien {\'e}quilibr{\'e}es, on peut leur appliquer
nos Th{\'e}or{\`e}mes \ref{thm:A1} et \ref{thm:1} et les exprimer
sous forme de sommes multiples.
La difficult\'e majeure pour en d\'eduire
le d{\'e}nominateur conjectur{\'e} par Zudilin
r{\'e}side dans le sommande de la somme multiple ainsi obtenue, qu'il
semble tr\`es compliqu\'e
d'arranger sous forme d'une d{\'e}composition en
briques {\'e}l{\'e}mentaires ou sp{\'e}ciales.

Par ailleurs, sous cette forme, cette conjecture n'est probablement
pas assez forte pour une {\'e}ventuelle
 application diophantienne, puisque, comme remarqu\'e ci-dessus, 
sur\break l'exemple de la s{\'e}rie
${\bf Z}_n$, on a $\dd_{33n}<{\hat{\Phi}}_n$.
Il est donc en fait  essentiel  
de g{\'e}n{\'e}raliser notre Th{\'e}or{\`e}me~\ref{prop:Phi}, c'est-{\`a}-dire
de d{\'e}terminer dans chaque situation des  entiers ${\check{\Phi}}_n$
tels que, en poursuivant sur l'exemple de ce paragraphe, l'on ait
${\hat{\Phi}}_n\le \dd_{33n}{\check{\Phi}}_n$ et
$$
2\dd_{35n} ^3 \dd_{34n} \dd_{33n} ^7
{\check{\Phi}}_n^{-1}
{\bf z}_{j,n}\in\,\mathbb{Z},
$$
ce qui serait alors au moins aussi bon que l'estimation~\eqref{zuddenom}.

\subsection{La conjecture des d{\'e}nominateurs li{\'e}e aux
valeurs de la fonction beta}
\label{ssec:beta}

La fonction beta est d{\'e}finie par la s{\'e}rie de Dirichlet
$$
\beta(s)=\sum_{k=0}^{\infty} \frac{(-1)^k}{(2k+1)^s},
$$
dont il est bien connu que les valeurs aux entiers impairs $j$ sont dans $\mathbb{Q}\pi^j$.
En revanche, rien n'{\'e}tait connu sur
la nature arithm{\'e}tique des valeurs de beta aux entiers pairs,
avant l'article \cite{rivzud} o{\`u} il est d{\'e}montr{\'e} que~:
{\it {\og une infinit{\'e} des nombres $\beta(j)$, $j\ge 2$
pair, sont lin{\'e}airement ind{\'e}pendants sur $\mathbb{Q}$\fg} }
et {\it {\og au moins un des sept nombres
$\beta(2)$, $\beta(4), \ldots, \beta(14)$ est irrationnel \fg} }.
La m{\'e}thode est bas{\'e}e sur la construction de combinaisons lin{\'e}aires en les valeurs
de beta au moyen, de nouveau, de certaines
s{\'e}ries hyperg{\'e}om{\'e}triques tr{\`e}s bien {\'e}quilibr{\'e}es,
{\'e}ventuellement d{\'e}riv{\'e}es.
Par exemple, la s{\'e}rie suivante produit
des approximations rationnelles de la constante de Catalan $G=\beta(2)$~:
$$
{\bf G}_n=n!\sum_{k=1}^{\infty}(-1)^k\bigg(k+\frac{n-1}{2} \bigg) \,
\frac{(k-n)_n(k+n)_n}{\left(k-\frac 12\right)_{n+1}^3}={\bf e}_nG-{\bf f}_n,
$$
o{\`u} $2^{4n}\dd_{2n}{\bf e}_n$ et $2^{4n}\dd_{2n}^3{\bf f}_n$ sont entiers.
Il appara{\^\i}t num{\'e}riquement
que $2^{4n}{\bf e}_n$ et $2^{4n}\dd_{2n}^2{\bf f}_n$
sont d{\'e}j{\`a} entiers (voir \cite[fin du paragraphe~9]{rivzud}), ce qui
a {\'e}t{\'e} confirm{\'e}
par le deuxi{\`e}me auteur dans \cite{rivbeta} par une m{\'e}thode tr{\`e}s
indirecte
d'approximation de Pad{\'e}.
(Voir \cite{zu6} pour une d{\'e}monstration {\og asymptotique\fg}.)
Ce raffinement demeure
n{\'e}anmoins insuffisant pour montrer l'irrationalit{\'e} de $G$.

Bien que cela ne soit mentionn{\'e} ni dans \cite{rivzud}
ni dans \cite{rivbeta}, on peut formuler une conjecture des d{\'e}nominateurs
pour les s{\'e}ries de \cite{rivzud} tout aussi g{\'e}n{\'e}rale
que celle de Zudilin \cite{zu3}.
En utilisant les techniques de cet article,
nous sommes capables de d{\'e}montrer cette conjecture pour les analogues
de nos coefficients
${\bf p}_{l,n}\((-1)^A\)$ dans le cas des s{\'e}ries 
{\og sym{\'e}triques\fg}~:
$$
n!^{A-2Br}\sum_{k=1}^{\infty}\frac{1}{C!}\frac{\partial^C }{\partial k^C}
\left(\bigg(k+\frac{n-1}{2} \bigg) \,
\frac{(k-rn)_{rn}^B(k+n)_{rn}^B}{\left(k-\frac 12\right)_{n+1}^A}\right)(-1)^{kA}.
$$
Il est {\'e}vident que notre approche s'applique aussi {\`a} ces
s{\'e}ries (l{\'e}g{\'e}rement diff{\'e}rentes des s{\'e}ries
${\bf S}_{n,A,B,C,r}\((-1)^A\)$ de notre article) puisqu'elles sont aussi
des s{\'e}ries tr{\`e}s bien {\'e}quilibr{\'e}es.
Si, toutefois, on parvenait
{\`a} d{\'e}montrer la totalit{\'e} de cette conjecture,
on pourrait alors peut-{\^e}tre  d{\'e}montrer l'irrationalit{\'e}
d'au moins un des six nombres
$\beta(2)$, $\beta(4), \ldots, \beta(12)$.
Nous envisageons de revenir {\`a} ces questions dans une publication
future.

\subsection{La $q$-conjecture des d{\'e}nominateurs}
\label{ssec:q}

D{\'e}finissons, pour tous  $q$ et $s$
tels que $\vert q\vert <1$ et $s\ge 1$, la s{\'e}rie~
$$
\zeta_q(s)=\sum_{k=1}^{\infty} k^{s-1}\frac{q^k}{1-q^k},
$$
qui est un $q$-analogue de $\zeta(s)$ au sens suivant :
$$
\lim_{q\to 1 }(1-q)^s \zeta_q(s)=(s-1)!\,\zeta(s).
$$
En dehors du cas $s=1$, la nature arithm{\'e}tique des
valeurs de $\zeta_q(s)$ aux entiers $s$ impairs {\'e}tait
totalement inconnue, jusqu'{\`a} l'article \cite{KrRZAA} o{\`u} les
versions $q$-analogiques de certains des th{\'e}or{\`e}mes
rappel{\'e}s au paragraphe \ref{sec:histoire} ont {\'e}t{\'e}
{\'e}tablies~: pour tout $q\not=\pm 1$ tel que $1/q\in\,\mathbb{Z}$,
{\og{\it une infinit{\'e} des nombres
$\zeta_q(j)$, $j\ge 1$ impair, sont lin{\'e}airement ind{\'e}pendants
sur $\mathbb{Q}$}\fg} et {\og{\it
au moins un des cinq nombres
$\zeta_q(3), \zeta_q(5), \zeta_q(7), \zeta_q(9), \zeta_q(11)$ est
irrationnel}\fg}.

Les d{\'e}monstrations sont bas{\'e}es sur l'{\'e}tude d'une
s{\'e}rie hyperg{\'e}om{\'e}trique basique similaire {\`a}
celles de cet article, mais qui n'est cependant
pas toujours tr{\`e}s bien {\'e}quilibr{\'e}e.
Il est indiqu{\'e} au dernier paragraphe de \cite{KrRZAA} que
pour d{\'e}montrer les m{\^e}mes th{\'e}or{\`e}mes, on aurait tout aussi bien
pu utiliser une autre s{\'e}rie basique tr{\`e}s bien {\'e}quilibr{\'e}e, qui
redonne la s{\'e}rie $\bar{\bf S}_{n,A,r}(1)$ du
paragraphe \ref{ssec:hisinfzeta}, lorsque $q$ tend vers 1.
Avec cette s{\'e}rie alternative, appara{\^\i}t alors une $q$-conjecture des
d{\'e}nominateurs, dont la preuve permettrait de montrer
l'irrationalit{\'e} d'au moins un des quatre
nombres $\zeta_q(3), \zeta_q(5), \zeta_q(7), \zeta_q(9)$.
Or des versions $q$-analogiques de nos Th{\'e}or{\`e}mes~\ref{thm:A1}
et \ref{thm:1} sont connues~: comme nous l'avons mentionn{\'e} au
paragraphe~\ref{gigantesques}, Andrews a en fait
d{\'e}montr{\'e} un $q$-analogue du Th{\'e}or{\`e}me~\ref{thm:A1} dans
\cite{Andr}, et un $q$-analogue de la transformation \eqref{Sears}
(ici n{\'e}cessaire pour d{\'e}duire le
Th{\'e}or{\`e}me~\ref{thm:1} du Th{\'e}or{\`e}me~\ref{thm:A1}) est
aussi connu (voir \cite[(2.10.4)~; 
Appendix (III.15)]{GaRaAA}).
Il est donc possible que
l'on puisse d{\'e}montrer cette $q$-conjecture des
d{\'e}nominateurs plus facilement que les autres conjectures
mentionn{\'e}es dans ce paragraphe.

\subsection{La version non-termin{\'e}e des identit{\'e}s gigantesques}
\label{ssec:nonterm}

Le paragraphe~\ref{sec:hyp} pr{\'e}\-sente des identit{\'e}s pour
des s{\'e}ries infinies
qui relient des s{\'e}ries {\og {\`a} mauvais\fg} d{\'e}nominateurs (plus
pr{\'e}cis{\'e}ment, la
s{\'e}rie de Ball $\bf B_n$ pour $\zeta(3)$ et la
s{\'e}rie~\eqref{eq:Ball2} pour $\zeta(2)$) aux s{\'e}ries
{\og {\`a} bons\fg} d{\'e}nominateurs.
On peut aussi consid{\'e}rer
que le th{\'e}or{\`e}me de Zudilin~\cite[paragraphe~8]{zu2}, liant les int{\'e}grales de Vasilyev {\`a}
une s{\'e}rie hyperg{\'e}om{\'e}trique tr{\`e}s bien {\'e}quilibr{\'e}e,
est une identit{\'e} de ce type~: en effet, comme nous l'avons expliqu{\'e}
au paragraphe~\ref{sec:explications},
les int{\'e}grales de Vasilyev se d{\'e}veloppent naturellement en
une s{\'e}rie multiple infinie et c'est la comparaison des deux
membres de l'identit{\'e} de Zudilin qui nous a permis de deviner
certaines propositions au paragraphe~\ref{sec:hyperharm}.

Il est donc envisageable que ce type d'identit\'es puisse fournir une approche plus simple 
et/ou  globale pour les diverses conjectures des d\'enominateurs. Dans ce but, il 
est donc utile d'obtenir les versions
non-termin{\'e}es de nos Th{\'e}or{\`e}mes~\ref{thm:A1} et~\ref{thm:1},
ou m{\^e}me plus
simplement de certaines de leurs sp{\'e}cialisations. 
De fait, dans le cas du Th\'eor\`eme~\ref{thm:A1}, 
nous avons d\'ej\`a obtenu une telle version non-termin\'ee (voir~\cite{KR}), qui nous a permis de donner 
une nouvelle d\'emonstration purement hyperg\'eom\'etrique du th\'eor\`eme de Zudilin mentionn\'e ci-dessus. 

On ne doit cependant pas n{\'e}gliger le fait que l'approche par les
s{\'e}ries hyperg{\'e}om{\'e}triques multiples
pr{\'e}sente pour le moment de s{\'e}rieuses difficult{\'e}s. Pour les
d{\'e}crire, nous avons besoin  d'introduire les  nombres
polyz{\^e}tas, d{\'e}finis par
$$
\zeta(\underline j)=\sum_{1\le k_s< \cdots< k_1}
\frac{1}{k_1^{j_1}k_2^{j_2}\cdots k_s^{j_s}},
$$
o{\`u} $s\ge 1$ et $\underline{j}=(j_1, j_2, \dots, j_s)$
est un $s$-uplet d'entiers $\ge 1$,
avec $j_1\ge 2$ (voir~\cite{miw} pour une introduction aux propri{\'e}t{\'e}s de ces
nombres). L'int{\'e}grale de Vasilyev de dimension $E\ge 2$, {\'e}crite sous forme
d'int{\'e}grale de Sorokin, se d{\'e}veloppe
en une s{\'e}rie hyperg{\'e}om{\'e}trique multiple ${\bf S}_n$. Son sommande
est une fraction rationnelle que l'on d{\'e}compose en {\'e}l{\'e}ments simples pour obtenir
$$
{\bf S}_n= {\bf p}_{0,n} +
\sum_{\underline{j}\in J} {\bf p}_{\underline{j},n}\,\zeta(\underline j),
$$
o{\`u} $J$ est un certain ensemble fini (ne
d{\'e}pendant pas de $n$) de $s$-uplets avec $s$ variant entre $1$ et $E$, et
${\bf p}_{\underline{j},n}$ sont des rationnels, ainsi que ${\bf p}_{0,n}$.  
Il est extr{\^e}mement difficile de d{\'e}crire {\it pr{\'e}cis{\'e}ment}
l'ensemble $J$, c'est-{\`a}-dire
de d{\'e}terminer les $\underline{j}$ tels que ${\bf
  p}_{\underline{j},n}\neq 0$.
Cependant, il est possible de d{\'e}terminer un d{\'e}nominateur commun aux rationnels ${\bf p}_{0,n}$ et 
${\bf p}_{\underline{j},n}$,
qui, dans ce cas, est $\dd_n^E$, comme on s'y attend {\it a priori}.

Par ailleurs, le th{\'e}or{\`e}me de Zudilin
exprime ${\bf S}_n$ comme
une s{\'e}rie hyperg{\'e}om{\'e}trique tr{\`e}s bien {\'e}quilibr{\'e}e $S_n$, qui v{\'e}rifie~:
$$
S_n= p_{0,n} + \underset{j\equiv \,E \,\textup{(mod 2)}}{\sum_{j=2}^E}
p_{j,n}\,\zeta(j),
$$
o{\`u} les rationnels $p_{j,n}$ admettent $\dd_n^{E+1}$ comme d{\'e}nominateur
{\it a priori}. Ainsi,
l'identit{\'e}  ${\bf S}_n=S_n$ devient~:
\begin{equation}
\label{eq:polyzetazeta}
{\bf p}_{0,n} + \sum_{\underline{j}\in J}
{\bf p}_{\underline{j},n}\,\zeta(\underline j)
=p_{0,n} + \underset{j\equiv \,E \,\textup{(mod 2)}}{\sum_{j=2}^E}
p_{j,n}\,\zeta(j).
\end{equation}
Il est donc tentant de se dire que l'on peut d{\'e}montrer
la conjecture des d{\'e}nominateurs dans ce cas en comparant les
coefficients des deux membres. Malheureusement, lorsque $E$ est impair,
c'est actuellement impossible pour deux raisons.
La premi{\`e}re est que nous ne savons pas encore montrer directement que
$J=\{3, 5, \ldots, E\}$, m{\^e}me dans le cas le plus simple o{\`u} $E=3$.
La deuxi{\`e}me est que  m{\^e}me si on savait lever la premi{\`e}re
difficult{\'e},
on ne pourrait conclure que les coefficients rationnels
co{\"\i}ncident que sous l'hypoth{\`e}se que les valeurs de z{\^e}ta aux entiers impairs
sont toutes lin{\'e}airement ind{\'e}pendantes sur $\mathbb{Q}$, ce qui, 
\'evidemment, 
est loin d'{\^e}tre av{\'e}r{\'e} actuellement. Si $E$ est pair, comme les
valeurs de z{\^e}ta aux entiers pairs
sont en revanche bien lin{\'e}airement ind{\'e}pendantes sur $\mathbb{Q}$, on pourrait
{\'e}ventuellement conclure si l'on savait montrer que $J=\{2, 4, \ldots,
E\}$, mais  c'est pour l'instant tout aussi impossible
(sauf quand $E=2$, voir le paragraphe \ref{sec:hyp}). 
Une mani\`ere de lever ces difficult\'es 
serait d'obtenir un analogue {\em fonctionnel} convenable
de~\eqref{eq:polyzetazeta}, car les propri\'et\'es d'ind\'ependance 
lin\'eaire sur 
$\mathbb{C}(z)$ des polylogarithmes $\Li_s(z)$ et de leurs versions 
multiples sont 
bien connues (voir~\cite[Theorem 4, paragraphe~5]{petitot}). 
Dans cette direction, indiquons que 
Zudilin propose dans~\cite[Theorem 5, p. 224]{zu7} une identit\'e int\'egrale 
qui pourrait donner une version fonctionnelle de~\eqref{eq:polyzetazeta} 
de la forme souhait\'ee.


\section*{Remerciements}
Nous ne saurions assez remercier K. Srinivasa Rao et
Michel Waldschmidt  d'avoir
organis{\'e} l'``Int\-er\-na\-tion\-al
Conference on Special Functions'' {\`a} Chennai (Madras) au mois de  septembre
2002. Sans cette conf{\'e}rence, o{\`u} nous nous sommes
rencontr{\'e}s pour la premi{\`e}re fois et
o{\`u} ce travail a d{\'e}but{\'e}, nous n'aurions
probablement jamais pris contact. Nous remercions {\'e}galement
Wadim Zudilin pour ses commentaires pertinents.


\def\refname{Bibliographie}

\message{ }
\message{ }
\message{ }
\message{ }
\message{ }
\message{ }
\message{ }
\message{ }
\message{ }
\message{ }
\message{ }
\message{ }
\message{ }
\message{ }
\message{ }
\message{ }
\message{ }
\message{ }
\message{ }
\message{ }
\message{ }
\message{ }
\message{ }
\message{ }
\message{ }
\message{ }
\message{ }
\message{ }
\message{ }
\message{ }
\message{ }
\message{ }
\message{ }
\message{ }
\message{ }
\message{ }
\message{ }
\message{ }
\message{ }
\message{ }
\message{ }
\message{ }
\message{ }
\message{ }
\message{ }
\message{ }
\message{ }
\message{ }
\message{ }
\message{ }
\message{ }
\message{ }
\message{ }
\message{ }
\message{ }
\message{ }
\message{ }
\message{ }
\message{ }
\message{ }
\message{ }
\message{ }
\message{ }
\message{ }
\message{ }
\message{ }
\message{ }
\message{ }
\message{ }
\message{ }
\message{ }
\message{ }
\message{ }
\message{!!!!!!!!!!!!!!!!!!!!!!!!!!!!!!!!!!!!!!!!!!!!!!!!}
\message{Attention! Il faut adapter les nombres des pages}
\message{dans la table des matieres a la fin!}
\message{!!!!!!!!!!!!!!!!!!!!!!!!!!!!!!!!!!!!!!!!!!!!!!!!}
\message{ }
\message{ }
\message{ }
\message{ }
\message{ }
\message{ }
\message{ }
\message{ }
\message{ }
\message{ }
\message{ }
\message{ }
\message{ }
\message{ }
\message{ }
\message{ }
\message{ }
\message{ }
\message{ }
\message{ }
\message{ }
\message{ }
\message{ }
\message{ }
\message{ }
\message{ }
\message{ }
\message{ }
\message{ }
\message{ }
\message{ }
\message{ }
\message{ }
\message{ }
\message{ }
\message{ }
\message{ }
\message{ }
\message{ }
\message{ }
\message{ }
\message{ }
\message{ }
\message{ }
\message{ }
\message{ }
\message{ }
\message{ }
\message{ }
\message{ }
\message{ }
\message{ }

\end{document}